\documentclass[11pt]{article}

\usepackage{amssymb}
\usepackage{amsmath}
\usepackage{amscd}
\usepackage[dvips]{graphicx}
\usepackage{mathrsfs}
\usepackage{picins}
\usepackage{setspace}
\numberwithin{equation}{section}

\begin{document}

\title{On the Sojourn Time Distribution in a Finite Population Markovian Processor Sharing Queue
\thanks 
{We dedicate this paper to the memory of John A. Morrison (1927-2013). John was an outstanding applied mathematician who spent his entire career (1956-1995) at Bell Laboratories (and then Alcatel-Lucent), continuing to be very active in research after his retirement. John worked in many different areas of applied mathematics, including wave motion, elastic materials, satellite orbit motion, stochastic differential equations and communications networks. 
Much of his latter career was spent analyzing queueing models, by a combination of classical analysis and singular perturbation methods. In particular he made fundamental contributions to processor sharing queues, including the finite population model treated here. John was an excellent problem solver that could make progress on extremely difficult models.}
\vspace{.5cm}}

\author{
Qiang Zhen\thanks{
Department of Mathematics and Statistics,
University of North Florida, 1 UNF Dr, Jacksonville, FL 32224-7699, USA.
{\em Email:} q.zhen@unf.edu}
\and
Johan S.H. van Leeuwaarden\thanks{
Eindhoven University of Technology and EURANDOM, P.O. Box 513, 5600 MB Eindhoven, The Netherlands.
{\em Email:} j.s.h.v.leeuwaarden@tue.nl
}
\and
Charles Knessl\thanks{
Department of Mathematics, Statistics, and Computer Science,
University of Illinois at Chicago, 851 South Morgan (M/C 249),
Chicago, IL 60607-7045, USA.
{\em Email:} knessl@uic.edu\
\newline\indent\indent{\bf Acknowledgement:} C. Knessl was partly supported by NSA grants H 98230-08-1-0102 and H 98230-11-1-0184. 
}}

\date{ }
\maketitle

\begin{abstract}
\noindent We consider a finite population processor-sharing (PS) queue, with Markovian arrivals and an exponential server. Such a queue can model an interactive computer system consisting of a bank of terminals in series with a central processing unit (CPU). For systems with a large population $N$ and a commensurately rapid service rate, or infrequent arrivals, we obtain various asymptotic results. We analyze the conditional sojourn time distribution of a tagged customer, conditioned on the number $n$ of others in the system at the tagged customer's arrival instant, and also the unconditional distribution. The asymptotics are obtained by a combination of singular perturbation methods and spectral methods. We consider several space/time scales and parameter ranges, which lead to different asymptotic behaviors. We also identify precisely when the finite population model can be approximated by the standard infinite population $M/M/1$-PS queue.  

\end{abstract}

\section{Introduction}

In this paper we present various asymptotic expressions for the sojourn times in a finite population processor-sharing (PS) queue.
The sojourn time is defined as the time it takes for a tagged customer to get through the system (after having obtained the required amount of service). With the PS discipline, first introduced by Kleinrock in the 1960's \cite{KLE_A,KLE_T}, we mean that if there are $\mathcal{N}(t)>0$ customers in the system at time $t$, each gets an equal fraction $1/\mathcal{N}(t)$ of the server's capacity. The PS discipline is often a good model for
computer and communication systems (see \cite{HEY}-\cite{NAB}). The PS discipline has the advantage over, say, first-in first-out (FIFO), in that shorter jobs tend to get through the system more rapidly, due to the egalitarian way of serving. 
PS also serves as  a sensible approximation
for round robin (RR) disciplines used in computer operating systems, or head-of-the-line 
fair-queueing service disciplines used in communication network routers. Our asymptotic expressions serve as predictions for sojourn times that provide valuable information, both for the customers and the system manager. 

 The most classical processor-sharing model is the $M/M/1$-PS queue. Here customers arrive according to a Poisson process with rate parameter $\lambda$ and the service requirements are exponentially distributed with mean $\mu^{-1}$.  The sojourn time is a random variable denoted throughout by $R$. For the simplest $M/M/1$ model, the distribution of $R$ depends on the total service time $\mathcal{X}$ that the customer requests and also on the number of other customers present when the tagged customer enters the system. 
 
 One natural variant of the $M/M/1$-PS model is the finite population model, which puts an upper bound on the number of customers that can be served by the processor. The model assumes that there are a total of $N$ customers, and each customer will enter service in the next $\Delta t$ time units with probability $\lambda \Delta t+o(\Delta t)$. At any time there are $\mathcal{N}(t)$ customers being served and the remaining $N-\mathcal{N}(t)$ customers are in the general population. Hence the total arrival rate is $\lambda [N-\mathcal{N}(t)]$ and we may view the model as a PS queue with a state-dependent arrival rate that decreases linearly to zero. Once a customer finishes service that customer re-enters the general population. The service times are exponentially distributed with mean $1/\mu$ and we define the traffic intensity $\rho$ by $\rho=\lambda N/\mu$. This model may describe, for example, a network of $N$ terminals in series with a processor-shared CPU, which can be viewed as a closed two node queueing network.

The finite population model does not seem amenable to an exact solution. However, various asymptotic studies have been done in the limit $N\to\infty$, so that the total population, or the number of terminals, is large. If $N$ is large it is reasonable to assume either that $\lambda $, the arrival rate of each individual customer, is small, of the order $O(N^{-1})$, or that the service rate $\mu$ is large, of the order $O(N)$. Then $\rho=\lambda N/\mu$ will remain $O(1)$ as $N\to\infty$. Previous studies of the finite population model were carried out by Morrison and Mitra (see \cite{MIT}-\cite{MOR_C88}), in each case for $N\to\infty$. For example, the moments of the sojourn time $R$ conditioned on the service time $\mathcal{X}$ are obtained in \cite{MOR_A}, where it was found that the asymptotics are very different according as $\rho<1$ (called ``normal usage''), $\rho-1=O(N^{-1/2})$ (called ``heavy usage''), or $\rho>1$ (called ``very heavy usage''). In \cite{MOR_M} the unconditional sojourn time distribution is investigated for $N\to\infty$ and the three cases of $\rho$, in \cite{MOR_C87} the author obtains asymptotic results for the conditional sojourn time distribution, conditioned on the service time $\mathcal{X}$, in the very heavy usage case $\rho>1$, and in \cite{MOR_C88} the results of \cite{MOR_C87} are generalized to multiple customer classes (here the population $N$ is divided into several classes, with each class having different arrival and service times). In \cite{MIT} the authors analyze the multiple class model and obtain the unconditional sojourn time moments for $N\to\infty$ in the normal usage case, while in \cite{MOR_H} heavy usage results are obtained.

We shall focus on the sojourn time of customers in the finite population model as $N\to\infty$ with a fixed $\rho=O(1)$. This corresponds to systems with infrequent arrivals (small $\lambda$) or, equivalently, rapid service (large $\mu$).
Let $\mathcal{N}(0^-)$ be the number of other customers in the system immediately before the tagged customer arrives (time $t=0^-$). Thus, if the total population is $N$ we have $0\le \mathcal{N}(0^-)\le N-1$.
The conditional sojourn time density is defined by
\begin{equation*}\label{s2_pn(t)}
p_n(t)dt=\mathrm{Prob}\big[R\in(t,t+dt)\big| \mathcal{N}(0^-)=n\big].
\end{equation*}
Removing the condition on $n$ will yield the unconditional sojourn time density, $p(t)$, as
\begin{equation}\label{s2_p(t)}
p(t)=\sum_{n=0}^{N-1}\pi_n\,p_n(t)
\end{equation}
where
\begin{equation}\label{s2_pin}
\pi_n=\frac{(N-1)!}{(N-n-1)!}\Big(\frac{\rho}{N}\Big)^n\bigg/\sum_{l=0}^{N-1}\frac{(N-1)!}{(N-l-1)!}\Big(\frac{\rho}{N}\Big)^l,
\end{equation}
and we note that the right side of (\ref{s2_pin}) is the steady state distribution of customers in a finite population (or repairman) model with population $N-1$, which is one less than that in our model (see \cite{SEV}). We shall consider $\rho<1$ and $\rho>1$ with $N$ large, and also $\rho\to\infty$ with $N$ fixed. However, we do not treat the transition case where $N\to\infty$ and $\rho\to 1$ simultaneously.

For both the conditional and unconditional sojourn time density we will derive asymptotic expression by using a combination of singular perturbation methods and spectral methods. We consider several space/time scales and parameter ranges, which lead to different asymptotic behaviors. We also identify precisely when, and when not, the finite population model can be approximated by the standard infinite population $M/M/1$-PS queue.

Intuitively, for a fixed $\rho<1$ and $N\to\infty$ we would expect that the finite population model may be approximated by the infinite population $M/M/1$-PS queue. We shall show that this is true for time scales $t=O(1)$ and when there are only a few other customers present when the tagged customer arrives, i.e., $n=O(1)$. But for larger space/time scales this is certainly not true, and we shall show that very different asymptotics come into play for $t=O(N)$ and $n=O(N)$. For example, let us set $\xi=n/N$ and $\tau=t/N$, so that $\xi$ is the fraction of the population that is in the system and $\tau$ is a large time scale. If, say, $\xi=1/2$ we will obtain different asymptotic expansions for the sojourn time density $p_n(t)$ for the three time ranges $0\le \tau<\tau_0(1/2)$, $\tau_0(1/2)<\tau<\tau_*(1/2)$ and $\tau>\tau_*(1/2)$, where $\tau_0(\cdot)$ and $\tau_*(\cdot)$ are defined in (\ref{s33_tau0}) and (\ref{s33_tau*}). In the first time range the density $p_n(t)$ varies slowly, on the $\tau$-time scale, and most of the probability mass accumulates in this range of $\tau$. In the other two ranges the density is roughly exponentially small in $N$, and these correspond to the tail(s) of the distribution. For $\tau_0(1/2)<\tau<\tau_*(1/2)$ the dependence on $\tau$ is intricate, while for $\tau>\tau_*(1/2)$ the dependence on $\tau$ is exponential, and for the latter the dominant eigenvalue begins to determine the asymptotics. Note that the spectrum of the infinite matrix that arises in the analysis of the infinite population model is purely continuous (and leads to the interesting tail behavior in (\ref{s31_p(t)sigma})), but the finite matrix (cf. (\ref{s2_A})) in the finite population model has clearly a purely discrete spectrum. Our analysis for $N\to\infty$ will indicate how the continuous spectrum separates into the discrete eigenvalues, as $\tau$ increases through the critical value $\tau_*(1/2)$. We shall also consider in detail the transition ranges $\tau\approx \tau_0(1/2)$ and $\tau\approx\tau_*(1/2)$, where the approximation to $p_n(t)$ will be given in terms of parabolic cylinder functions and elliptic integrals, respectively. We will also study the spatial scales $n=O(1)$ and $n=O(N^{1/4})$, which lead to yet different asymptotics, and show that then the time scale $t=O(N^{3/4})$ is also important. For the unconditional density $p(t)$, we shall show that the finite population model behaves similarly to the infinite population model as long as $t\ll N^{3/4}$, but a ``phase transition" occurs for $t=O(N^{3/4})$, and for $t\gg N^{3/4}$ the behavior of $p(t)$ becomes purely exponential in time. There is also a critical time range (see (\ref{s31_sigma*scale}) with $\sigma=tN^{-3/4}$) where all of the eigenvalues of the matrix contribute roughly equally to the expansion of $p(t)$, and then we obtain an interesting double-exponential approximation (see (\ref{s31_sigma*})).

We comment that our approach is analytical rather than probabilistic. The singular perturbation analysis does make certain assumptions about the forms of the asymptotic expansions and the asymptotic matching between the different space/time scales. We shall obtain the conditional and unconditional densities for time ranges where most of the probability mass concentrates, but also treat carefully the tails, where the densities are often exponentially small in $N$. In the tails we shall show (see (\ref{s33_R2pn}) and (\ref{s33_R3pn})) that the expansion of $\log\big[p_n(t)\big]$ has terms of order $O(N)$, $O(\log N)$ and $O(1)$, and sometimes additional terms of orders $O(\sqrt{N})$ and $O(N^{1/4})$. Then several terms in the expansion of $\log\big[p_n(t)\big]$ must be computed in order to get the leading term for $p_n(t)$ itself. In contrast, probabilistic approaches based on large deviations theory typically compute only the limit of $N^{-1}\log\big[p_n(t)\big]$ for $N\to\infty$. We also indicate how to compute higher order terms in the various asymptotic expansions, and compute these explicitly for some of the space/time ranges. 

The paper is organized as follows. In section 2 we give a precise statement of the mathematical problem and obtain the basic equations. In section 3 we summarize our main results and briefly discuss them. The derivations appear in sections 4-6. In section 7 we discuss an alternate approach based on spectral results obtained in \cite{ZHE_O12}, while in section 8 we consider the limit of fixed $N$ and large $\rho$.

\section{Statement of the problem}


We easily obtain for $p_n(t)$ the equation(s)
\begin{equation}\label{s2_recu}
p_n'(t)=\frac{\rho}{N}(N-n-1)\big[p_{n+1}(t)-p_n(t)\big]+\frac{n}{n+1}p_{n-1}(t)-p_n(t),\; 0\le n\le N-1,
\end{equation}
and the initial condition is 
\begin{equation}\label{s2_ic}
p_n(0)=\frac{1}{n+1},
\end{equation}
which can be obtained by integrating (\ref{s2_recu}) from $t=0$ to $t=\infty$, and noting that $\int_0^\infty p_n(t)dt=1$ for each $n$. The finite system of ODEs in (\ref{s2_recu}) may also be written in matrix form as
\begin{equation}\label{s2_recuMatrx}
\mathbf{p}'(t)=\mathbf{A}\mathbf{p}(t)
\end{equation}
where $\mathbf{p}(t)$ is the column vector $\mathbf{p}(t)=\big(p_0(t),p_1(t),\cdots,p_{n-1}(t)\big)^T$ ($T$=transpose) and $\mathbf{A}$ is the tridiagonal matrix
\begin{equation}\label{s2_A}
\mathbf{A}=\left[\begin{array}{cccccccc}
-1-\rho\big(\frac{N-1}{N}\big) & \rho\big(\frac{N-1}{N}\big) & 0 & 0 & 0 & \cdots & 0 & 0\\
1/2 & -1-\rho\big(\frac{N-2}{N}\big) & \rho\big(\frac{N-2}{N}\big) & 0 & 0 & \cdots & 0 & 0 \\
0 & 2/3 & -1-\rho\big(\frac{N-3}{N}\big) & \rho\big(\frac{N-3}{N}\big) & 0 & \cdots & 0 & 0 \\
\vdots  & \vdots  & \vdots  & \vdots  & \vdots & \ddots & \vdots  & \vdots \\
0 & 0 & 0 &0 & 0 & \cdots & \frac{N-1}{N}  & -1
\end{array}\right].
\end{equation}

Thus formally we can write the solution as $\mathbf{p}(t)=e^{\mathbf{A}t}\mathbf{p}(0)$, but this does not yield much immediate insight. We can also write the solution in terms of the eigenvalues and eigenvectors of the matrix $\mathbf{A}$. Denoting the eigenvalues by $-\nu_j$ for $j=0,1,\cdots,N-1$ with $0<\nu_0<\nu_1<\cdots<\nu_{N-1}$, and the corresponding eigenvectors by $\phi_j(n)$, we see from (\ref{s2_recu}) that
\begin{equation}\label{s2_recu_eigen}
-\nu_j\phi_j(n)=\rho\Big(1-\frac{n+1}{N}\Big)\big[\phi_j(n+1)-\phi_j(n)\big]+\frac{n}{n+1}\phi_j(n-1)-\phi_j(n),\; 0\le n\le N-1.
\end{equation}
From (\ref{s2_recu_eigen}) we can easily establish the orthogonality relation
\begin{equation}\label{s2_orth}
\sum_{n=0}^{N-1}(n+1)\pi_n\,\phi_j(n)\phi_k(n)=\delta_{jk}\Big[\sum_{n=0}^{N-1}(n+1)\pi_n\,\phi_j^2(n)\Big].
\end{equation}
Thus the spectral representation of the solution to (\ref{s2_recu}) and (\ref{s2_ic}) is given by
\begin{equation}\label{s2_pn(t)eigen}
p_n(t)=\sum_{j=0}^{N-1}e^{-\nu_jt}c_j\phi_j(n),
\end{equation}
where
\begin{eqnarray}\label{s2_cj}
c_j&=&\sum_{n=0}^{N-1}\pi_n\,\phi_j(n)\bigg/\sum_{n=0}^{N-1}(n+1)\pi_n\,\phi_j^2(n)\nonumber\\
&=&\frac{\sum_{n=0}^{N-1}\frac{N!}{(N-n-1)!}\big(\frac{\rho}{N}\big)^n\phi_j(n)}{\sum_{n=0}^{N-1}\frac{N!}{(N-n-1)!}\big(\frac{\rho}{N}\big)^n(n+1)\phi_j^2(n)}.
\end{eqnarray}
Note that (\ref{s2_pn(t)eigen}) with (\ref{s2_cj}) is invariant under $\phi_j(n)\to c\phi_j(n)$ for any constant $c$, so that the eigenvectors need not be normalized. The spectral expansion of the unconditional density is then given by 
\begin{equation}\label{s2_p(t)eigen}
p(t)=\sum_{j=0}^{N-1}d_je^{-\nu_jt},\quad d_j=\frac{\sum_{n=0}^{N-1}\pi_n\,\phi_j^2(n)}{\sum_{n=0}^{N-1}(n+1)\pi_n\,\phi_j^2(n)}.
\end{equation}

Unfortunately it does not seem possible to explicitly obtain the eigenvalues and eigenvectors. In \cite{ZHE_O12} we did a detailed study of these, asymptotically for $N\to\infty$. In particular we showed that if $\rho>1$
\begin{equation}\label{s2_nuj>1}
\nu_j=\frac{1}{N}\Big(\rho\,j+\frac{\rho}{\rho-1}\Big)+O(N^{-2})
\end{equation}
while if $\rho<1$
\begin{eqnarray}\label{s2_nuj<1}
\nu_j&=&(1-\sqrt{\rho})^2+2\sqrt{\rho}\sqrt{1-\sqrt{\rho}}\,N^{-1/2}+(2j+1)\sqrt{\rho}(1-\sqrt{\rho})^{3/4}\,N^{-3/4}\nonumber\\
&&+\bigg[\frac{\sqrt{\rho}(22\sqrt{\rho}-3\rho-15)}{16(1-\sqrt{\rho})}-\frac{3}{8}\sqrt{\rho}(1-\sqrt{\rho})j(j+1)\bigg]N^{-1}+O(N^{-5/4}).
\end{eqnarray}
These results apply to the eigenvalue index $j$ being $O(1)$, though in (\ref{s2_nuj<1}) $j$ is allowed to be slightly large, as long as $j=o(N^{1/4})$. We also obtained in \cite{ZHE_O12} the large $N$ asymptotics of the eigenvectors, where it proved necessary to approximate $\phi_j(n)$ differently in different ranges of $n$, such as $n=O(1)$ and $n=O(N)$. For sufficiently large times $t$, no matter how large $N$ is, the eigenvalue $\nu_0$ dominates and we have
\begin{equation}\label{s2_pasym}
p_n(t)\sim c_0\phi_0(n)e^{-\nu_0t},\quad p(t)\sim d_0e^{-\nu_0t},\; t\to\infty.
\end{equation}
However, for other time ranges the spectral representation in (\ref{s2_pn(t)eigen}) is not useful for understanding the behavior of the sojourn time density. 

Our goal here is to analyze both $p_n(t)=p_n(t;N,\rho)$ and $p(t)=p(t;N,\rho)$ for $N\to\infty$. We shall consider separately the cases $\rho<1$ and $\rho>1$, and in either case it is necessary to analyze several different space-time ranges, as we summarize in section 3. We do not consider the critical case $\rho\approx 1$ here, but comment that the scale $\rho-1=O(N^{-2/3})$ is important in understanding the eigenvalue problem \cite{ZHE_O12}, and to see how the expansion in (\ref{s2_nuj>1}) transitions, as $\rho\downarrow 1$, to that in (\ref{s2_nuj<1}). In \cite{MOR_A} it was shown that the scale $\rho-1=O(N^{-1/2})$ is important for the unconditional density $p(t)$. 

In section 7, we shall give an alternate derivation of some of our asymptotic results, which uses the spectral expansion and the large $N$ results for the $\nu_j$ and $\phi_j(n)$ in \cite{ZHE_O12}. However, this is useful only in time ranges where (\ref{s2_pasym}) applies, or on slightly shorter times where all the eigenvalues $\nu_j$ contribute about equally to the expansions of $p_n(t)$ and $p(t)$. This will be made more precise later. For now we note that if $\rho>1$, (\ref{s2_nuj>1}) shows that all the eigenvalues are small and of order $O(N^{-1})$, and that the leading term depends on the eigenvalue index $j$. This suggests that for times $t=O(N)$ one sees the distinct eigenvalues, and that for times $t\gg O(N)$, $\nu_0$ should dominate. In contrast, if $\rho<1$, (\ref{s2_nuj<1}) shows that the eigenvalue index does not appear until the third ($O(N^{-3/4})$) term in the asymptotic expansion. This suggests that on time scales $t=O(N^{3/4})$ all the eigenvalues contribute roughly equally, but for $t\gg O(N^{3/4})$, $\nu_0$ should dominate. We make these criteria more precise in section 3. Also, we note that the leading term in (\ref{s2_nuj<1}), which is independent of $j$, is precisely the relaxation rate in the standard (infinite population) $M/M/1$ model. 

In addition to using a spectral expansion, we could also try to analyze (\ref{s2_recu}) using transform methods, but there seems to be no standard transform that would be appropriate for projecting out the spatial variable $n$. In Appendix A we briefly discuss using a generating function over $n$, but this leads to a complicated functional-partial differential equation, whose chances of solution are minimal. 

We can solve (\ref{s2_recu}) and (\ref{s2_ic}) for some very small values of $N$. When $N=2$ we obtain
\begin{equation*}\label{s2_N=2}
\left\{ \begin{array}{l}
\displaystyle p_0(t)=\frac{1}{2}\sqrt{\frac{\rho+4}{\rho}}\,\bigg[\frac{1-\nu_0}{\nu_0}e^{-\nu_0t}-\frac{1-\nu_1}{\nu_1}e^{-\nu_1t}\bigg],\\
\\
\displaystyle p_1(t)=\frac{1}{4}\sqrt{\frac{\rho+4}{\rho}}\,\bigg[\frac{1}{\nu_0}e^{-\nu_0t}-\frac{1}{\nu_1}e^{-\nu_1t}\bigg],
\end{array} \right.
\end{equation*}
where 
\begin{equation*}\label{s2_N=2nu}
\nu_0=1+\frac{\rho}{4}-\frac{1}{4}\sqrt{\rho(\rho+4)},\quad \nu_1=1+\frac{\rho}{4}+\frac{1}{4}\sqrt{\rho(\rho+4)}.
\end{equation*}
Note that as $\rho\to\infty$, $\nu_0\to 1/2$ while $\nu_1\sim \rho/2$.

When $N=3$ the result is already less explicit, with
\begin{equation*}\label{s2_N=3}
\left\{ \begin{array}{l}
\displaystyle p_0(t)=\frac{\rho\, C_0}{\nu_0-1-2\rho/3}\,e^{-\nu_0t}+\frac{\rho\, C_1}{\nu_1-1-2\rho/3}\,e^{-\nu_1t}+\frac{\rho\, C_2}{\nu_2-1-2\rho/3}\,e^{-\nu_2t},\\
\\
\displaystyle p_1(t)=-\frac{3}{2}\Big[C_0e^{-\nu_0t}+C_1e^{-\nu_1t}+C_2e^{-\nu_2t}\Big],\\
\\
\displaystyle p_2(t)=\frac{C_0}{\nu_0-1}\,e^{-\nu_0t}+\frac{C_1}{\nu_1-1}\,e^{-\nu_1t}+\frac{C_2}{\nu_2-1}\,e^{-\nu_2t},\\
\end{array} \right.
\end{equation*}
the $C_j$ are obtained by solving
\begin{equation*}\label{s2_N=3Cj}
\left\{ \begin{array}{l}
\displaystyle 1=\frac{\rho\, C_0}{\nu_0-1-2\rho/3}+\frac{\rho\, C_1}{\nu_1-1-2\rho/3}+\frac{\rho\, C_2}{\nu_2-1-2\rho/3},\\
\\
\displaystyle \frac{1}{2}=-\frac{3}{2}\big(C_0+C_1+C_2\big),\\
\\
\displaystyle \frac{1}{3}=\frac{C_0}{\nu_0-1}+\frac{C_1}{\nu_1-1}+\frac{C_2}{\nu_2-1},\\
\end{array} \right.
\end{equation*}
and the $\nu_j$ are the three roots of the cubic polynomial
\begin{equation}\label{s2_N=3nu}
\Big(\nu-1-\frac{2}{3}\rho\Big)\Big[\Big(\nu-1-\frac{1}{3}\rho\Big)(\nu-1)-\frac{2}{9}\rho\Big]=\frac{1}{3}\rho(\nu-1),
\end{equation}
again ordered as $0<\nu_0<\nu_1<\nu_2$. It is easy to show from (\ref{s2_N=3nu}) that $\nu_0\in(0,1)$, $\nu_1\in(1,1+\rho/3)$ and $\nu_2>1+2\rho/3$. As $\rho\to\infty$ we have $\nu_0\to 1/3$, $\nu_1\sim \rho/3$ and $\nu_2\sim 2\rho/3$. 

Our observations for the simple cases $N=2$ and $N=3$ suggest that for any $N$, $\nu_0\to 1/N$ as $\rho\to\infty$, and $\nu_j\sim \rho\,j/N$ for $j\ge 1$. Note that this conjecture, which assumes that we fix $N$ and let $\rho\to\infty$, is consistent with (\ref{s2_nuj>1}), which lets $N\to\infty$ for a fixed $\rho>1$. 

It seems fruitless to pursue any further exact solutions to (\ref{s2_recu}) so we begin the asymptotic analysis.

\section{Summary of results}

Here we shall summarize our main results, and briefly discuss and interpret them, for the cases $\rho\lessgtr 1$ and for both the conditional and unconditional probabilities. We discuss first the conditional sojourn time distribution for $N\to\infty$, and then the unconditional distribution.

\subsection{Conditional distribution for $\rho>1$}

We consider $N\to\infty$ with a fixed $\rho>1$. We use the scaled space/time variables 
\begin{equation}\label{s32_scales}
n=N\xi,\quad t=N\tau
\end{equation}
so that $\xi$ $(\in(0,1))$ is essentially the fraction of the population that is being serviced by the PS server. The asymptotic expansion of the sojourn time distribution will be different for the three scales (i) $n,\,t=O(N)$, (ii) $n=O(1)$, $t=O(N)$ and (iii) $n,\,t=O(1)$. 

For $n,\,t=O(1)$ we have
\begin{equation}\label{s32_O(1)}
p_n(t)=p_n(t;\,M/M/1-\mathrm{PS})+O(N^{-1}),
\end{equation}
where $p_n(t;M/M/1-\textrm{PS})$ refers to the conditional sojourn time distribution in the standard \textit{infinite population} model. The infinite population model was analyzed by Coffman, Muntz and Trotter \cite{COF}, who also conditioned the response time $R$ on the service time that the tagged customer will receive. Also, the model is equivalent to the random order service (ROS) $M/M/1$ queue (see \cite{COH}, \cite{BOR}), which was analyzed already in 1946 by Pollaczek \cite{POL}. The $O(N^{-1})$ correction term may be obtained by analyzing (\ref{s2_recu}), as it will satisfy an inhomogeneous version of the problem for the leading term in (\ref{s32_O(1)}). But here we focus on the behavior of $p_n(t)$ on larger time scales, where the differences between the finite and infinite population models will become apparent. 
If we denote the correction term in (\ref{s32_O(1)}) as $p_n^{(1)}(t)N^{-1}$, then $p_n^{(1)}$ will satisfy
\begin{equation}\label{s32_recu}
\frac{d}{dt}p_n^{(1)}(t)-\rho\big[p_{n+1}^{(1)}(t)-p_n^{(1)}(t)\big]-\frac{n}{n+1}p_{n-1}^{(1)}(t)+p_n^{(1)}(t)=-\rho(n+1)\big[p_{n+1}^{(0)}(t)-p_n^{(0)}(t)\big]
\end{equation}
and the initial condition $p_n^{(1)}(0)=0$. In (\ref{s32_recu}) $p_n^{(0)}$ is the infinite capacity solution, and (\ref{s32_recu}) may be analyzed by transform or spectral methods, but we shall not do so. The leading term $p_n^{(0)}=p_n(t;\,M/M/1-\mathrm{PS})$ can be expressed, for example, as the integral
\begin{equation}\label{s32_p^(0)}
p_n^{(0)}(t)=\frac{1}{2\pi i}\int_{Br}\widehat{p}_n(\theta)\,e^{\theta t}d\theta,
\end{equation}
\begin{equation}\label{s32_phat}
\widehat{p}_n(\theta)=\int_{\mathcal{C}_*}z_+z_-\frac{(1-z_+)^{\alpha_1-1}}{(1-z_-)^{\alpha_1}}\,\frac{e^{\pi i\alpha_1}}{1-e^{2\pi i\alpha_1}}z^n(z_+-z)^{-\alpha_1}(z-z_-)^{\alpha_1-1}dz,
\end{equation}
where
\begin{equation*}\label{s32_zpm}
z_\pm(\theta)=\frac{1}{2\rho}\Big[1+\rho+\theta\pm\sqrt{(1+\rho+\theta)^2-4\rho}\Big],
\end{equation*}
\begin{equation}\label{s32_alpha1}
\alpha_1=\frac{z_+}{z_+-z_-}=\frac{1+\rho+\theta+\sqrt{(1+\rho+\theta)^2-4\rho}}{2\sqrt{(1+\rho+\theta)^2-4\rho}},
\end{equation}
$Br$ is a vertical Bromwich contour in the complex $\theta$-plane on which ℜ$\Re(\theta)>0$,
and $\mathcal{C}_*$ is a closed loop around the branch cut, where $\Im(z)=0$ and $z_-<\Re(z)<z_+$. The branches in (\ref{s32_phat}) are chosen so that the integrand, as a function of $z$, is analytic in the complex $z$-plane, exterior to this cut. In Figure \ref{figure:1} we sketch $\mathcal{C}_*$. 

For $n=O(1)$, corresponding to having only a few customers in the system, and large times $t=N\tau=O(N)$ we shall obtain 
\begin{equation}\label{s32_O(N)}
p_n(t)\sim N^{-\alpha_0}\,\alpha_0^{2\alpha_0-1}\,\Gamma(\alpha_0)\big(1-e^{-\rho\tau}\big)^{-\alpha_0}e^{-\alpha_0\tau}\,\frac{e^{\pi i\alpha_0}}{2\pi i}\oint_\mathcal{C}s^n(1-s)^{-\alpha_0}\Big(s-\frac{1}{\rho}\Big)^{\alpha_0-1}ds,
\end{equation}
where $\alpha_0=\rho/(\rho-1)$
is obtained by setting $\theta=0$ in (\ref{s32_alpha1}), and the closed contour $\mathcal{C}$ encircles the branch points at $s=1/\rho$ and $s=1$, with now the integrand chosen to be analytic exterior to the slit $\Im(s)=0$, $1/\rho<\Re(s)<1$. 

We can show that (\ref{s32_O(1)}), with (\ref{s32_p^(0)}), when expanded for $t\to\infty$, agrees with the expansion of (\ref{s32_O(N)}) as $\tau\to 0$. In the matching region $1\ll t\ll N$, $p_n(t)$ is proportional to $t^{-\alpha_0}$ and thus shows an algebraic behavior, which is also the tail of $p_n(t;\,M/M/1-\mathrm{PS})$ when $\rho>1$. However, for times $t=O(N)$ ($\tau=O(1)$) the time dependence is more complicated, as can be seen from (\ref{s32_O(N)}). For times where $\tau\to\infty$ ($t\gg N$) (\ref{s32_O(N)}) becomes proportional to the exponential $e^{-\alpha_0\tau}$, and this corresponds to $p_n(t)\sim c_0\phi_0(n)e^{-\nu_0t}$, as in (\ref{s2_pasym}), as now the zeroth eigenvalue dominates the asymptotics. Note also that the contour integral in (\ref{s32_O(N)}) may be expressed in terms of hypergeometric functions.

Next we consider the $(\xi,\tau)$ scale, so we are looking at the problem on large space/time scales. Then we obtain 
\begin{equation}\label{s32_xitau}
p_n(t)=N^{-1}P(\xi,\tau)+N^{-2}P^{(1)}(\xi,\tau)+O(N^{-3})
\end{equation} 
where
\begin{equation}\label{s32_Pxitau}
P(\xi,\tau)=\xi^{\frac{1}{\rho-1}}\exp\Big(-\frac{\rho\tau}{\rho-1}\Big)\Big[\Big(\xi-\frac{\rho-1}{\rho}\Big)e^{-\rho\tau}+\frac
{\rho-1}{\rho}\Big]^{-\frac{\rho}{\rho-1}}
\end{equation}
and
\begin{equation}\label{s32_Pxitau1}
P^{(1)}(\xi,\tau)=P_\xi(\xi,\tau)+P(\xi,\tau)C(\xi,\tau),
\end{equation}
\begin{eqnarray}\label{s32_Cxitau}
C(\xi,\tau)&=&\frac{2\rho-1}{2(\rho-1)^2}\frac{e^{\rho\tau}-1}{(\xi+\Delta_*)^2}\Big(e^{\rho\tau}+\rho-\rho\xi\Big)-\frac{\rho(\rho+1)}{(\rho-1)^3}\frac{e^{\rho\tau}-1}{\xi+\Delta_*}+\frac{\rho(3-\rho)}{2(\rho-1)^3}\Big[\frac{1}{\xi}-\frac{e^{\rho\tau}}{\xi+\Delta_*}\Big]\nonumber\\
&&+\frac{\rho}{(\rho-1)^4}\Big[1+2\rho^2\frac{\xi-1+1/\rho}{\xi+\Delta_*}\Big]\log\Big(\frac{\xi+\Delta_*}{\xi}\Big),
\end{eqnarray}
\begin{equation*}\label{s32_Delta*}
\Delta_*=\Delta_*(\tau)=\frac{\rho-1}{\rho}(e^{\rho\tau}-1).
\end{equation*}
As $\tau\to 0^+$ we have $P(\xi,0)=1/\xi$ and $C(\xi,0)=0$, so that (\ref{s32_xitau}) becomes $N^{-1}\big[\xi^{-1}-N^{-1}\xi^{-2}+O(N^{-2})\big]$, which is consistent with the exact initial value $p_n(0)=1/(n+1)=N^{-1}(\xi+N^{-1})^{-1}=N^{-1}\big[\xi^{-1}-N^{-1}\xi^{-2}+O(N^{-2})\big]$. Thus (\ref{s32_xitau}) holds also for shorter time scales, including $t=O(1)$, as long as $\xi>0$. When $\rho>1$ we clearly have $\xi+(1-\rho^{-1})(1-e^{-\rho\tau})>0$ except for $\xi=\tau=0$ (when $\rho<1$, however, this quantity may vanish and this will lead to some very different asymptotics of $p_n(t)$). As $\xi\to 0$, both (\ref{s32_Pxitau}) and (\ref{s32_Cxitau}) have singular behavior, with the leading term in (\ref{s32_xitau}) becoming
\begin{equation}\label{s32_leading}
N^{-1}\xi^{\frac{1}{\rho-1}}\Big(\frac{\rho}{\rho-1}\Big)^{\frac{\rho}{\rho-1}}(1-e^{-\rho\tau})^{-\frac{\rho}{\rho-1}}\exp\Big(-\frac{\rho\tau}{\rho-1}\Big).
\end{equation}
But (\ref{s32_leading}) can also be obtained by expanding (\ref{s32_O(N)}) for $n\to\infty$, using singularity analysis to approximate the contour integral. Noting that $N^{-\alpha_0}n^{1-\alpha_0}=N^{-1}\xi^{1/(\rho-1)}$, this thus verifies the asymptotic matching, for a fixed $\tau>0$, of the expansions in (\ref{s32_O(N)}) and (\ref{s32_xitau}). 

We have thus shown that when $\rho>1$ the expansion of $p_n(t)$ is different in three different space/time ranges, with (\ref{s32_xitau}) applying in the majority of the $(\xi,\tau)$ plane, except near $\xi=0$ where the solution involves the contour integral in (\ref{s32_O(N)}), and near the ``corner" $(\xi,\tau)=(0,0)$, where for $n,\,t=O(1)$ we approximate the present model by the infinite population $M/M/1$-PS queue. If we expand (\ref{s32_xitau}) for $\tau\to\infty$ (then $t\gg N$) the leading term in (\ref{s32_Pxitau}) will be proportional to $\exp\big[-\rho\tau/(\rho-1)\big]$, so then the zeroth eigenvalue dominates. The correction term in (\ref{s32_Pxitau1}) will have, for $\tau\to\infty$, additional terms proportional to $\exp\big[-\rho^2\tau/(\rho-1)\big]$ and $\exp\big[-\rho\tau/(\rho-1)-2\rho\tau\big]$, and such terms correspond to contributions from the first and second eigenvalues ($\nu_1$ and $\nu_2$), in view of the asymptotic result in (\ref{s2_nuj>1}). But, for times $\tau=O(1)$ all of the eigenvalues contribute to even the leading term (\ref{s32_xitau}), as can be seen from (\ref{s32_Pxitau}) and (\ref{s2_nuj>1}).

\subsection{Conditional distribution for $\rho<1$}

This case is much more complicated, and we will need to consider a total of 10 space/time scales where the asymptotic behavior of $p_n(t)$ is different. Viewed on the coarse $(\xi,\tau)$ space/time scale, there will be 3 main regions and various internal transition layers and boundary layers, the latter along $\xi=0$.

We begin by introducing two curves in the $(\xi,\tau)$ plane, defined by 
\begin{equation}\label{s33_tau0}
\tau_0(\xi)=\frac{1}{\rho}\log\Big(1+\frac{\rho\xi}{1-\rho}\Big)
\end{equation}
and
\begin{equation}\label{s33_tau*}
\tau_*(\xi)=\frac{1}{\rho}\log\bigg[1+\frac{\rho\xi+\sqrt{\rho^2\xi^2+4\rho\xi(1-\sqrt{\rho})}}{2(1-\sqrt{\rho})}\bigg].
\end{equation}
We have $\tau_0(0)=0=\tau_*(0)$, and for $\xi>0$,
$0<\tau_0(\xi)<\tau_*(\xi)$.
For large $\xi$ we have 
\begin{equation*}\label{s33_tau0large}
\tau_0(\xi)=\frac{1}{\rho}\log\xi+\frac{1}{\rho}\log\Big(\frac{\rho}{1-\rho}\Big)+o(1),\; \xi\to\infty
\end{equation*}
\begin{equation*}\label{s33_tau*large}
\tau_*(\xi)=\frac{1}{\rho}\log\xi+\frac{1}{\rho}\log\Big(\frac{\rho}{1-\sqrt{\rho}}\Big)+o(1),\; \xi\to\infty
\end{equation*}
and the inverse functions to (\ref{s33_tau0}) and (\ref{s33_tau*}) are given by 
\begin{equation*}\label{s33_xi0}
\xi_0(\tau)=\frac{1-\rho}{\rho}(e^{\rho\tau}-1)
\end{equation*}
and
\begin{equation*}\label{s33_xi*}
\xi_*(\tau)=\frac{1-\sqrt{\rho}}{\rho}e^{\rho\tau}(1-e^{-\rho\tau})^2=\frac{4(1-\sqrt{\rho})}{\rho}\sinh^2\Big(\frac{\rho\tau}{2}\Big).
\end{equation*}
Note that $\xi_0(\tau)$ hits the corner $(\xi,\tau)=(0,0)$ along a line, while $\xi_*(\tau)$ hits the corner along a parabola.

The curvers $\tau_0$ and $\tau_*$ divide the $(\xi,\tau)$ plane into the three regions $R_1$, $R_2$, $R_3$, with
\begin{equation*}\label{s33_Regions}
\begin{array}{l}
\displaystyle R_1=\big\{(\xi,\tau):\;\xi>0,\,0\le \tau<\tau_0(\xi)\big\},\\
\displaystyle R_2=\big\{(\xi,\tau):\;\xi>0,\,\tau_0(\xi)< \tau<\tau_*(\xi)\big\},\\
\displaystyle R_3=\big\{(\xi,\tau):\;\xi>0,\,\tau>\tau_*(\xi)\big\}.
\end{array}
\end{equation*}
The regions are sketched in Figure \ref{figure:2}. Note that the $R_j$ are defined so as to exclude the curves that separate them, and to also exclude $\xi=0$.

In region $R_1$ the expansion (\ref{s32_xitau}) holds. But now that $\rho<1$ the leading term vanishes when $\tau=\tau_0(\xi)$, and becomes complex for $\tau>\tau_0(\xi)$. Also, from (\ref{s32_Cxitau}) we see that the correction term has a singularity as $\xi\uparrow \xi_0(\tau)$ ($=-\Delta_*(\tau)$). Thus (\ref{s32_xitau}) becomes invalid as we exit region $R_1$.

In region $R_2$ we shall show that 
\begin{equation}\label{s33_R2pn}
p_n(t)\sim e^{N\phi(\xi,\tau)}N^{\phi^{(1)}(\xi,\tau)}K(\xi,\tau)
\end{equation}
where
\begin{equation}\label{s33_R2phi}
\phi(\xi,\tau)=\frac{1-Be^{-\rho\tau}}{\rho}\Big(1+e^{\rho\tau}\frac{B+\rho-1}{1-B}\Big)\log(1-Be^{-\rho\tau})-\log(1-B)-\frac{B}{\rho}(1-e^{-\rho\tau})
\end{equation}
and
\begin{equation}\label{s33_R2B}
B=B(\xi,\tau)=\frac{e^{2\rho\tau}-1-\rho\xi e^{\rho\tau}-\sqrt{(\rho\xi e^{\rho\tau}+1-e^{2\rho\tau})^2-4e^{\rho\tau}(e^{\rho\tau}-1)\big[(1-\rho)(e^{\rho\tau}-1)-\rho\xi \big]}}{2(e^{\rho\tau}-1)}.
\end{equation}
Furthermore,
\begin{equation}\label{s33_R2phi(1)}
\phi^{(1)}(\xi,\tau)=-\frac{3(1-B)^2-\rho}{2(1-B)^2-2\rho}
\end{equation}
and
\begin{equation}\label{s33_R2K}
K(\xi,\tau)=K_0(B)e^{5\rho\tau/2}\frac{\big[B^2-B+e^{\rho\tau}(1-\rho-B)\big]^{\frac{\rho}{(1-B)^2-\rho}}}{\sqrt{(1-B)^2+\rho e^{\rho\tau}}\sqrt{e^{\rho\tau}-B}}(e^{\rho\tau}-1)^{\frac{\rho-3(1-B)^2}{2(1-B)^2-2\rho}}, 
\end{equation}
with
\begin{eqnarray}\label{s33_R2K0}
K_0(B)&=&\frac{1}{\sqrt{2\pi}}(1-B)\rho^{\frac{3(1-B)^2-\rho}{2(1-B)^2-2\rho}}B^{\frac{(1-B)^2}{\rho-(1-B)^2}}\Big(1-\frac{\rho}{1-B}\Big)^{\frac{\rho}{(1-B)^2-\rho}}\nonumber\\
&&\times\Big(1-B-\frac{\rho}{1-B}\Big)^{\frac{(1-B)^2+2\rho}{\rho-(1-B)^2}}\Gamma\bigg(\frac{(1-B)^2}{(1-B)^2-\rho}\bigg)
\end{eqnarray}
and $\Gamma(\cdot)$ is the Gamma function. Despite the rather complicated form, we note that $\phi$, $\phi^{(1)}$ and $K$, and thus the right side of (\ref{s33_R2pn}), are explicit functions of the $(\xi,\tau)$ variables. The curve $B=0$ corresponds to $\xi=\xi_0(\tau)$ while the curve $B=1-\sqrt{\rho}$ corresponds to $\xi=\xi_*(\tau)$, which separates regions $R_2$ and $R_3$. In region $R_2$ we have $0<B<1-\sqrt{\rho}$. We can show that in $R_2$, $\phi(\xi,\tau)<0$ so that the conditional density is exponentially small in $N$. The factor $N^{\phi^{(1)}}$ gives an algebraic correction factor, and we note that $\phi^{(1)}$ is constant along curves in the $(\xi,\tau)$ plane which are level curves of $B=B(\xi,\tau)$. When $B=0$ we have $\phi^{(1)}=-(3-\rho)/(2-2\rho)$ and when $B\uparrow 1-\sqrt{\rho}$, $\phi^{(1)}$ becomes singular, as does $K$ in (\ref{s33_R2K}).

In region $R_3$, $p_n(t)$ will again be exponentially small, but now the leading term is
\begin{equation}\label{s33_R3pn}
p_n(t)\sim e^{N\psi (\xi,\tau)}e^{\sqrt{N}\psi^{(1)}(\xi,\tau)}e^{N^{1/4}\psi^{(2)}(\xi,\tau)}N^{-9/8}L(\xi,\tau)
\end{equation}
where
\begin{eqnarray}\label{s33_R3psi}
\psi(\xi,\tau)&=&-(1-\sqrt{\rho})^2\tau+\xi\log 2+\frac{\xi}{2}-\frac{1}{2\rho}\sqrt{\rho^2\xi^2+4\rho\xi(1-\sqrt{\rho})}\nonumber\\
&&-\xi\log\bigg[2\sqrt{\rho}-\rho\xi+\sqrt{\rho^2\xi^2+4\rho\xi(1-\sqrt{\rho})}\bigg]\nonumber\\
&&+\frac{\rho-4\sqrt{\rho}+2}{2\rho}\log\bigg[\frac{2-2\sqrt{\rho}+\rho\xi+\sqrt{\rho^2\xi^2+4\rho\xi(1-\sqrt{\rho})}}{2(1-\sqrt{\rho})}\bigg]\nonumber\\
&&+\frac{1}{2}\log\bigg[\frac{(\rho-2\sqrt{\rho}+2)\xi+2(1-\sqrt{\rho})+(2-\sqrt{\rho})\sqrt{\rho\xi^2+4\xi(1-\sqrt{\rho})}}{2(1-\sqrt{\rho})}\bigg].
\end{eqnarray}
Setting 
\begin{equation}\label{s33_R3B0}
B_0=\frac{e^{\rho\tau}}{2}\bigg[2-2\sqrt{\rho}+\rho\xi-\sqrt{\rho^2\xi^2+4\rho\xi(1-\sqrt{\rho})}\bigg]
\end{equation}
we note that $B_0=1-\sqrt{\rho}$ along $\xi=\xi_*(\tau)$ and $B_0>1-\sqrt{\rho}$ in region $R_3$. Then
{\setlength\arraycolsep{1pt}
\begin{eqnarray}\label{s33_R3psi(1)}
\psi^{(1)}(\xi,\tau)&=&\frac{2\sqrt{1-\sqrt{\rho}}}{\sqrt{\rho}}\log\Big(\frac{1-\sqrt{\rho}}{B_0}\Big)\nonumber\\
&=&-2\sqrt{\rho}\sqrt{1-\sqrt{\rho}}\,\tau-\frac{2\sqrt{1-\sqrt{\rho}}}{\sqrt{\rho}}\log\bigg[\frac{2-2\sqrt{\rho}+\rho\xi-\sqrt{\rho^2\xi^2+4\rho\xi(1-\sqrt{\rho})}}{2(1-\sqrt{\rho})}\bigg],
\end{eqnarray}}
\begin{eqnarray}\label{s33_R3psi(2)}
\psi^{(2)}(\xi,\tau)&=&\frac{(1-\sqrt{\rho})^{3/4}}{\sqrt{\rho}}\log\Big(\frac{1-\sqrt{\rho}}{B_0}\Big)-\frac{8}{3}(1-\sqrt{\rho})^{-1/4}\nonumber\\
&=&-\sqrt{\rho}(1-\sqrt{\rho})^{3/4}\tau-\frac{8}{3}(1-\sqrt{\rho})^{-1/4}\nonumber\\
&&-\frac{(1-\sqrt{\rho})^{3/4}}{\sqrt{\rho}}\log\bigg[\frac{2-2\sqrt{\rho}+\rho\xi-\sqrt{\rho^2\xi^2+4\rho\xi(1-\sqrt{\rho})}}{2(1-\sqrt{\rho})}\bigg],
\end{eqnarray}
and
\begin{eqnarray}\label{s33_R3L}
L(\xi,\tau)&=&L_0(B_0)\sqrt{2}\rho\,(1-\sqrt{\rho})^{-3/4}\bigg[2(1-\sqrt{\rho})+\rho\xi+\sqrt{\rho^2\xi^2+4\rho\xi(1-\sqrt{\rho})}\bigg]^{5/2}\nonumber\\
&&\times\bigg[2\sqrt{\rho}(1-\sqrt{\rho})+\rho\xi+\sqrt{\rho^2\xi^2+4\rho\xi(1-\sqrt{\rho})}\bigg]^{-1/2} \bigg[\rho\xi+\sqrt{\rho^2\xi^2+4\rho\xi(1-\sqrt{\rho})}\bigg]^{-3/2}\nonumber\\
&&\times\bigg[4(1-\sqrt{\rho})+\rho\xi+\sqrt{\rho^2\xi^2+4\rho\xi(1-\sqrt{\rho})}\bigg]^{-1/2} \exp\bigg\{\frac{(\sqrt{\rho}-2)\rho\xi}{2(1-\sqrt{\rho})\sqrt{\rho^2\xi^2+4\rho\xi(1-\sqrt{\rho})}}\bigg\}\nonumber\\
&\equiv& L_0(B_0)J(\xi)
\end{eqnarray}
with
\begin{eqnarray}\label{s33_R3L0}
L_0(B_0)&=& 8(1-\sqrt{\rho})^{-5/8}\exp\Big(\frac{1}{1-\sqrt{\rho}}\Big)e^{-5/2}\exp\bigg\{\frac{22\sqrt{\rho}-3\rho-15}{16\sqrt{\rho}(1-\sqrt{\rho})}\log\Big(\frac{1-\sqrt{\rho}}{B_0}\Big)\bigg\}\nonumber\\
&=&8(1-\sqrt{\rho})^{-5/8}\exp\Big(\frac{1}{1-\sqrt{\rho}}\Big)e^{-5/2}\exp\bigg\{-\frac{\sqrt{\rho}(22\sqrt{\rho}-3\rho-15)}{16(1-\sqrt{\rho})}\tau\bigg\}\nonumber\\
&&\times\bigg[\frac{2-2\sqrt{\rho}+\rho\xi-\sqrt{\rho^2\xi^2+4\rho\xi(1-\sqrt{\rho})}}{2(1-\sqrt{\rho})}\bigg]^{\frac{15+3\rho-22\sqrt{\rho}}{16\sqrt{\rho}(1-\sqrt{\rho})}}.
\end{eqnarray}
Again the expression is complicated, but is completely explicit in terms of $\xi$ and $\tau$. In fact the dependence of $\psi$, $\psi^{(1)}$, and $\psi^{(2)}$ on $\tau$ is linear, while $L$ in (\ref{s33_R3L}) depends on $\tau$ in an exponential manner, in view of the second equality in (\ref{s33_R3L0}). In region $R_3$ the dominant eigenvalue $\nu_0$, whose four term asymptotic approximation is in (\ref{s2_nuj<1}) with $j=0$, determines the asymptotic behavior of $p_n(t)$. Expression (\ref{s33_R3pn}) corresponds to $e^{-\nu_0t}c_0\phi_0(n)$ where the functions $\psi$, $\psi^{(1)}$, $\psi^{(2)}$ and $L_0$ contain within them the first, second, third and fourth terms in (\ref{s2_nuj<1}), respectively. The factors which do not involve $\tau$ correspond to the large $N$ expansion of $c_0\phi_0(n)$, when $\rho<1$. Thus as we move from region $R_2$ to $R_3$ for a fixed $\xi>0$, the zeroth eigenvalue begins to dominate. 

Next we discuss the transitions from $R_1$ to $R_2$ and from $R_2$ to $R_3$. We thus define the asymptotic transition ranges $T_1$ and $T_2$ by the scalings
\begin{equation}\label{s33_T1scale}
T_1:\;\xi-\xi_0(\tau)=O(N^{-1/2}),\;\tau>0
\end{equation}
and
\begin{equation}\label{s33_T2scale}
T_2:\;\xi-\xi_*(\tau)=O(N^{-1/4}),\;\tau>0.
\end{equation}
This is equivalent to scaling $\tau-\tau_0(\xi)=O(N^{-1/2})$ and $\tau-\tau_*(\xi)=O(N^{-1/4})$, respectively.

In transition region $T_1$ we shall obtain
\begin{eqnarray}\label{s33_T1pn}
p_n(t)&\sim& N^{\frac{\rho-2}{2(1-\rho)}}\xi^{\frac{\rho-2}{2(1-\rho)}}(\xi\rho^2+1-\rho^2)^{\frac{\rho}{2(1-\rho)}}(\xi\rho+1-\rho)\nonumber\\
&&\times\frac{1}{\sqrt{2\pi}}\int_{-\Delta_1}^\infty(y+\Delta_1)^{\frac{\rho}{1-\rho}}\exp\Big[-\frac{1}{2}(1-\rho)^2y^2\Big]dy,
\end{eqnarray}
where
\begin{equation}\label{s33_Delta1}
\Delta_1=\frac{\sqrt{N}}{\sqrt{\xi}\sqrt{\rho^2\xi+1-\rho^2}}\Big[\xi-\frac{1-\rho}{\rho}(e^{\rho\tau}-1)\Big]=O(1).
\end{equation}
The integral in (\ref{s33_T1pn}) may be expressed as a parabolic cylinder function. In this transition region $p_n(t)=O(N^{(\rho-2)/(2-2\rho)})$, which is $o(N^{-1})$ so that $p_n(t)$ is smaller than in region $R_1$. As $\Delta_1\to\infty$ we are approaching $R_1$ and then (\ref{s33_T1pn}) asymptotically matches to $p_n(t)\sim N^{-1}P(\xi,\tau)$ as $\xi\downarrow \xi_0(\tau)$. In the matching region, where $\xi-\xi_0(\tau)\to 0^+$ and $\sqrt{N}\big[\xi-\xi_0(\tau)\big]\to +\infty$, the integral in (\ref{s33_T1pn}) may be approximated by $\sqrt{2\pi}(1-\rho)^{-1}(-\Delta_1)^{\rho/(1-\rho)}$. For $\Delta_1\to-\infty$ we are approaching region $R_2$ and then we have
\begin{eqnarray}\label{s33_T1R2}
&&\int_{-\Delta_1}^\infty(y+\Delta_1)^{\frac{\rho}{1-\rho}}\exp\Big[-\frac{1}{2}(1-\rho)^2y^2\Big]dy\nonumber\\
&\sim&\Gamma\Big(\frac{1}{1-\rho}\Big)(1-\rho)^{-\frac{2}{1-\rho}}(-\Delta_1)^{\frac{1}{1-\rho}}\exp\Big[-\frac{1}{2}(1-\rho)^2\Delta_1^2\Big],\; \Delta_1\to -\infty.
\end{eqnarray}
Using (\ref{s33_T1R2}) in (\ref{s33_T1pn}) we can show that (\ref{s33_T1pn}) as $\Delta_1\to -\infty$ agrees with (\ref{s33_R2pn}) (with (\ref{s33_R2phi})--(\ref{s33_R2K0})) as $\xi\uparrow \xi_0(\tau)$ (then $B\downarrow 0$). Thus (\ref{s33_T1pn}) matches to both of the $R_1$ and $R_2$ results, but for $\Delta_1=O(1)$ the expansion involves a special function.

For the transition region $T_2$ in (\ref{s33_T2scale}) we define $\Delta$ by
\begin{equation}\label{s33_T2Delta}
\tau-\tau_*(\xi)=N^{-1/4}\Delta,\quad \Delta=N^{1/4}\big[\tau-\tau_*(\xi)\big]=O(1).
\end{equation}
Then for $\xi>0$ and $\Delta=O(1)$ we shall obtain
\begin{equation}\label{s33_T2pn}
p_n(t)\sim e^{-(1-\sqrt{\rho})^2t}\rho^{-n/2}N^{-9/8}e^{NF(\xi)}e^{N^{1/4}f(\Delta)}g(\xi,\Delta),
\end{equation}
where $F(\xi)$ can be computed from (\ref{s33_R3psi}) as
\begin{equation}\label{s33_T2F}
F(\xi)=\psi(\xi,\tau)+(1-\sqrt{\rho})^2\tau+\frac{1}{2}\xi\log\rho
\end{equation}
(so that $e^{N\psi(\xi,\tau)}=e^{-(1-\sqrt{\rho})^2t}\rho^{-n/2}e^{NF(\xi)}$), and
\begin{equation}\label{s33_T2f(Delta)}
f(\Delta)=\frac{1}{3}\sqrt{\rho}A\Delta-\frac{2}{3}\int_0^\infty\frac{1}{\sqrt{v}\sqrt{(1-\sqrt{\rho})v^2+1+Av}}dv
\end{equation}
where $A=A(\Delta)$ is defined implicitly by the equation
\begin{equation}\label{s33_T2A}
\int_0^\infty\bigg[\frac{\sqrt{v}}{\sqrt{(1-\sqrt{\rho})v^2+1+Av}}-\frac{1}{\sqrt{(1-\sqrt{\rho})v}}\bigg]dv=2\sqrt{\rho}\Delta.
\end{equation}
Finally,
\begin{equation*}\label{s33_T2g}
g(\xi,\Delta)=g_0(\Delta)J(\xi)\exp\bigg\{-\frac{2(1-\sqrt{\rho})+\sqrt{\rho}\xi}{4(1-\sqrt{\rho})^2\sqrt{\rho\xi^2+4\xi(1-\sqrt{\rho})}}\big[A^2-4(1-\sqrt{\rho})\big]\bigg\}
\end{equation*}
where $J$ is defined in (\ref{s33_R3L}) and 
\begin{equation}\label{s33_T2g0}
g_0(\Delta)=\frac{\sqrt{2}}{(1-\sqrt{\rho})^{5/4}}\exp\Big[\frac{1+\sqrt{\rho}}{2(1-\sqrt{\rho})}\Big]\bigg\{\int_0^\infty\big[(1-\sqrt{\rho})v+1/v+A\big]^{-3/2}dv\bigg\}^{-1/2}.
\end{equation}

Thus in transition region $T_2$ the asymptotic approximation is not very explicit and involves the solution of (\ref{s33_T2A}). We note that the integrand in (\ref{s33_T2A})  is $O(v^{-3/2})$ as $v\to\infty$ and $O(v^{-1/2})$ as $v\to 0^+$ and is thus integrable at both endpoints. We will show in section 5 that $A\to +\infty$ as $\Delta\to -\infty$, while $A\to -2\sqrt{1-\sqrt{\rho}}$ as $\Delta\to +\infty$. More precisely we shall obtain the estimates
\begin{equation}\label{s33_T2Aasym-}
A=\rho(1-\sqrt{\rho})^2\Delta^2+\frac{1}{\rho(1-\sqrt{\rho})\Delta^2}\Big\{\frac{3}{2}-2\log\big[2\sqrt{\rho}(1-\sqrt{\rho})^{3/4}(-\Delta)\big]\Big\}+o(\Delta^{-2}),\; \Delta\to -\infty
\end{equation}
and
\begin{equation}\label{s33_T2Aasym+}
A+2\sqrt{1-\sqrt{\rho}}\sim 64\,e^{-4}\sqrt{1-\sqrt{\rho}}\,\exp\Big\{-2\sqrt{\rho}(1-\sqrt{\rho})^{3/4}\Delta\Big\},\;\Delta\to +\infty,
\end{equation}
which shows that the error term in the relation $A\sim -2\sqrt{1-\sqrt{\rho}}$ is exponentially small for $\Delta\to +\infty$. In the limits $\Delta\to \pm \infty$ we can use (\ref{s33_T2Aasym-}) and (\ref{s33_T2Aasym+}) to considerably simplify $f(\Delta)$ and $g_0(\Delta)$, obtaining
\begin{equation}\label{s33_T2fasym-}
f(\Delta)=\frac{1}{3}\rho^{3/2}(1-\sqrt{\rho})^2\Delta^3+\frac{1+4\log\big[2\sqrt{\rho}(1-\sqrt{\rho})^{3/4}(-\Delta)\big]}{2\sqrt{\rho}(1-\sqrt{\rho})\Delta}+o(\Delta^{-1}),\;\Delta\to -\infty,
\end{equation}
\begin{equation}\label{s33_T2g0asym-}
g_0(\Delta)\sim\frac{\rho^{1/4}}{(1-\sqrt{\rho})^{1/4}}\exp\bigg[\frac{1+\sqrt{\rho}}{2(1-\sqrt{\rho})}\bigg]\sqrt{-\Delta},\;\Delta\to -\infty
\end{equation}
and
\begin{eqnarray}\label{s33_T2fasym+}
f(\Delta)&\sim &-2\sqrt{\rho}\sqrt{1-\sqrt{\rho}}\,\Delta-\frac{8}{3}(1-\sqrt{\rho})^{-1/4}\nonumber\\
&&-32\,e^{-4}(1-\sqrt{\rho})^{-1/4}\,\exp\Big\{-2\sqrt{\rho}(1-\sqrt{\rho})^{3/4}\Delta\Big\},\;\Delta\to +\infty,
\end{eqnarray}
\begin{equation}\label{s33_T2g0asym+}
g_0(\Delta)\sim 8(1-\sqrt{\rho})^{-5/8}e^{-2}\exp\bigg[\frac{1+\sqrt{\rho}}{2(1-\sqrt{\rho})}\bigg]\exp\big[-\sqrt{\rho}(1-\sqrt{\rho})^{3/4}\Delta\big],\;\Delta\to +\infty.
\end{equation}
By using (\ref{s33_T2fasym-}) and (\ref{s33_T2g0asym-}) in (\ref{s33_T2pn}) we can verify that (\ref{s33_T2pn}) asymptotically matches to (\ref{s33_R2pn}), in an intermediate limit where $\tau-\tau_*(\xi)\to 0^-$ and $N^{1/4}\big[\tau-\tau_*(\xi)\big]\to -\infty$. Similarly, using (\ref{s33_T2fasym+}) and (\ref{s33_T2g0asym+}) we find that (\ref{s33_T2pn}) matches to (\ref{s33_R3pn}) in the limit where $\tau-\tau_*(\xi)\to 0^+$ and $N^{1/4}\big[\tau-\tau_*(\xi)\big]\to +\infty$. In the matching region between (\ref{s33_R3pn}) and (\ref{s33_T2pn}) we can also write
\begin{eqnarray}\label{s33_T2R3}
p_n(t)&\sim& 8(1-\sqrt{\rho})^{-5/8}e^{-2}\exp\bigg[\frac{1+\sqrt{\rho}}{2(1-\sqrt{\rho})}\bigg]\rho^{-n/2}N^{-9/8}J(\xi)\nonumber\\
&&\times e^{-(1-\sqrt{\rho})^2t}e^{NF(\xi)}\exp\Big[-2\sqrt{\rho}\sqrt{1-\sqrt{\rho}}\sqrt{N}(\tau-\tau_*(\xi))\Big]\nonumber\\
&&\times\exp\Big[-\sqrt{\rho}(1-\sqrt{\rho})^{3/4}N^{1/4}(\tau-\tau_*(\xi))-\frac{8}{3}(1-\sqrt{\rho})^{-1/4}N^{1/4}\Big]\nonumber\\
&&\times\sum_{j=0}^\infty\frac{e^{-4j}32^j(-1)^j}{(1-\sqrt{\rho})^{j/4}j!}\,e^{-2\sqrt{\rho}(1-\sqrt{\rho})^{3/4}j\widetilde{\Delta}}
\end{eqnarray}
where $\Delta$ and $\widetilde{\Delta}$ are related by
\begin{equation*}\label{s33_T2R3Delta}
\Delta=\frac{1}{8\sqrt{\rho}(1-\sqrt{\rho})^{3/4}}\log N+\widetilde{\Delta}
\end{equation*}
so that 
\begin{equation*}\label{s33_T2R3tau}
\tau=\tau_*(\xi)+\frac{1}{8\sqrt{\rho}(1-\sqrt{\rho})^{3/4}}\frac{\log N}{N^{1/4}}+\frac{\widetilde{\Delta}}{N^{1/4}}.
\end{equation*}
The result in (\ref{s33_T2R3}), which applies for fixed $\widetilde{\Delta}$ (so that $\Delta$ is slightly large and positive) can be interpreted in terms of the spectral expansion in (\ref{s2_pn(t)eigen}), as the $j^\mathrm{th}$ term in the series in (\ref{s33_T2R3}) corresponds to the $j^\mathrm{th}$ eigenvalue $\nu_j$. For $\widetilde{\Delta}=O(1)$ all eigenvalues (with $j=O(1)$) contribute to the expansion of $p_n(t)$, while if $\widetilde{\Delta}\to\infty$ the zeroth eigenvalue dominates. Thus we have shown that the relation $p_n(t)\sim c_0\phi_0(n)e^{-\nu_0t}$, $t\to\infty$ is, for large $N$, valid when time $t$ is such that   
\begin{equation}\label{s33_T2t}
t-N\tau_*(\xi)-\frac{N^{3/4}\log N}{8\sqrt{\rho}(1-\sqrt{\rho})^{1/4}}\gg O(N^{3/4}),\; 0<\xi<1.
\end{equation}
The spectral results in (\ref{s2_nuj<1}) would suggest that $t$ has to be larger than $O(N^{3/4})$ to distinguish $e^{-\nu_0t}$ from the other $e^{-\nu_jt}$, but (\ref{s33_T2t}) is a much more precise criterion.

Next we consider ranges where $\xi\to 0$, so that the fraction of the population being served is small. The various expansions in regions $R_1$, $R_2$, $R_3$, $T_1$ and $T_2$ all have singular behavior as $\xi\to 0$, indicating a lack of uniformity in the asymptotics. We will need to thus consider the following additional space/time scales:
\begin{equation*}
\begin{array}{lll}
\mathrm{(i)} & n=O(1), & t=O(1)\\
\mathrm{(ii)} & n=O(\sqrt{N}), & t=O(N^{3/4})\\
\mathrm{(iii)} & n=O(1), & t=O(N^{3/4})\\
\mathrm{(iv)} & n=O(\sqrt{N}), & t=O(N)\\
\mathrm{(v)} & n=O(1), & t=O(N).
\end{array}
\end{equation*}
We shall use the new space variable $x$ and time variable $\sigma$, with $n=\sqrt{N}x$, $t=N^{3/4}\sigma$ 
and we note that $\xi=x/\sqrt{N}$ and $\tau=N^{-1/4}\sigma$.

For $n,\,t=O(1)$ the result in (\ref{s32_O(1)}) (with (\ref{s32_p^(0)}) and (\ref{s32_phat})) again applies. Note, however, that while $p_n^{(0)}(t)$ had an algebraic tail for $t\to\infty$ if $\rho>1$, for $\rho<1$ it has roughly an exponential tail, with
\begin{eqnarray}\label{s33_O(1)}
p_n^{(0)}(t)&\sim& 2^{2/3}3^{-1/2}\frac{\rho^{-5/12}}{1-\sqrt{\rho}}\exp\Big(\frac{\sqrt{\rho}}{1-\sqrt{\rho}}\Big)\rho^{-n/2}t^{-5/6}e^{-(1-\sqrt{\rho})^2t}\exp\big(-2^{-3/2}3\pi^{2/3}\rho^{1/6}t^{1/3}\big)\nonumber\\
&&\times\frac{1}{2\pi i}\oint\frac{1}{z^{n+1}(1-z)}\exp\Big(\frac{1}{1-z}\Big)dz,\; t\to\infty.
\end{eqnarray}
Here the integral is over a small loop about $z=0$. In (\ref{s33_O(1)}) the contour integral may also be written as
\begin{equation}\label{s33_O(1)int}
\sum_{m=0}^\infty\frac{(m+n)!}{(m!)^2}=(n-1)!\big[(n+1)L_{-n}(1)-L_{-n}^{(1)}(1)\big]
\end{equation}
where $L$ and $L^{(1)}$ are the Laguerre and first generalized Laguerre polynomials, with negative integer indices. For any $n\ge 0$ the series in (\ref{s33_O(1)int}) is an integer multiple of $e$. 

Now consider the large time scale $t=N\tau=O(N)$. Then for $x=n/\sqrt{N}$ fixed with $0<x<\infty$ the expansion is
\begin{eqnarray}\label{s33_xtaupn}
p_n(t) &\sim & 8(1-\sqrt{\rho})^{-3/8}e^{-5/2}\exp\Big(\frac{1}{1-\sqrt{\rho}}\Big)\rho^{-n/2}N^{-3/4}e^{-(1-\sqrt{\rho})^2t}e^{-2\sqrt{\rho}\sqrt{1-\sqrt{\rho}}\,\sqrt{N}\,\tau}\nonumber\\
&&\times\exp\bigg\{N^{1/4}\Big[-\sqrt{\rho}(1-\sqrt{\rho})^{3/4}\tau+2\sqrt{x}-\frac{2}{3}\sqrt{1-\sqrt{\rho}}\,x^{3/2}-\frac{8}{3}(1-\sqrt{\rho})^{-1/4}\Big]\bigg\}\nonumber\\
&&\times\exp\bigg\{-\frac{\sqrt{\rho}(22\sqrt{\rho}-3\rho-15)}{16(1-\sqrt{\rho})}\tau\bigg\}\frac{e^{x^2/4}e^{(1-\sqrt{\rho})^{1/4}\sqrt{x}}}{x^{1/4}\big[1+(1-\sqrt{\rho})^{1/4}\sqrt{x}\big]}.
\end{eqnarray}
For $n=O(1)$ we have
\begin{eqnarray}\label{s33_ntaupn}
p_n(t)&\sim& \frac{16\sqrt{\pi}}{(1-\sqrt{\rho})^{3/8}}e^{-3}\exp\Big(\frac{1}{1-\sqrt{\rho}}\Big)\rho^{-n/2}e^{-(1-\sqrt{\rho})^2t}e^{-2\sqrt{\rho}\sqrt{1-\sqrt{\rho}}\,\sqrt{N}\,\tau}N^{-5/8}\nonumber\\
&&\times\exp\bigg\{N^{1/4}\Big[-\sqrt{\rho}(1-\sqrt{\rho})^{3/4}\tau-\frac{8}{3}(1-\sqrt{\rho})^{-1/4}\Big]\bigg\}\exp\bigg\{-\frac{\sqrt{\rho}(22\sqrt{\rho}-3\rho-15)}{16(1-\sqrt{\rho})}\tau\bigg\}\nonumber\\
&&\times\frac{1}{2\pi i}\oint\frac{1}{z^{n+1}(1-z)}\exp\Big(\frac{1}{1-z}\Big)dz.
\end{eqnarray}

For $n\to\infty$ a standard application of the saddle point method shows that 
\begin{equation}\label{s33_int}
\frac{1}{2\pi i}\oint\frac{1}{z^{n+1}(1-z)}\exp\Big(\frac{1}{1-z}\Big)dz\sim \frac{\sqrt{e}}{2\sqrt{\pi}}n^{-1/4}e^{2\sqrt{n}},\; n\to\infty,
\end{equation}
and then we immediately see that the large $n$ behavior of (\ref{s33_ntaupn}) agrees with the small $x$ behavior of (\ref{s33_xtaupn}), so these expansions asymptotically match. Also, after some calculation we can show that the small $\xi$ expansion of the right side of (\ref{s33_R3pn}) (with (\ref{s33_R3psi})--(\ref{s33_R3L0})) agrees with the large $x$ expansion of (\ref{s33_xtaupn}). For the latter we simply approximate $x^{-1/4}\big[1+(1-\sqrt{\rho})^{1/4}\sqrt{x}\big]^{-1}$ by $x^{-3/4}(1-\sqrt{\rho})^{-1/4}$. Note also that for $\xi\to 0$, in view of (\ref{s33_R3psi}) and (\ref{s33_T2F}), 
\begin{equation}\label{s33_Fxi0}
F(\xi)=-\frac{2}{3}(1-\sqrt{\rho})^{1/2}\xi^{3/2}+\frac{1}{4}\xi^2+O(\xi^{5/2}).
\end{equation} 
In both (\ref{s33_xtaupn}) and (\ref{s33_ntaupn}) the dominant asymptotic factor is $\rho^{-n/2}e^{-(1-\sqrt{\rho})^2t}$, with the remaining factors showing sub-exponential or algebraic dependence on the large parameter $N$.

Finally we consider the time scale $t=N^{3/4}\sigma=O(N^{3/4})$ and it is on this scale where the asymptotics are the most difficult. On the $(x,\sigma)$ scale the expansion has the form
\begin{equation}\label{s33_xsigmapn}
p_n(t)\sim \rho^{-n/2}e^{-(1-\sqrt{\rho})^2t}N^{-3/4}e^{N^{1/4}\eta(x,\sigma)}\gamma(x,\sigma)
\end{equation}
where the functions $\eta$ and $\gamma$ have different forms in 3 regions of the $(x,\sigma)$ plane (see Figure \ref{figure:3}). We call these regions $D_1$, $D_2$ and $D_3$. 

Let $B_1=B_1(x,\sigma)$ be the solution to 
\begin{equation}\label{s33_xsigmaB1}
\int_0^x\big[(1-\sqrt{\rho})v+1/v+B_1\big]^{-1/2}dv=2\sqrt{\rho}\sigma.
\end{equation}
Then region $D_1$ is defined by the conditions $B_1\ge -2\sqrt{1-\sqrt{\rho}}$ and $\eta_x<0$, and there 
\begin{equation}\label{s33_etaD1}
\eta(x,\sigma)=B_1\sqrt{\rho}\,\sigma-\int_0^x\sqrt{(1-\sqrt{\rho})v+1/v+B_1}\,dv.
\end{equation}
Note that the left side of (\ref{s33_xsigmaB1}) is real only for $B_1>-2\sqrt{1-\sqrt{\rho}}$, and when $B_1=-2\sqrt{1-\sqrt{\rho}}$ we have $\sqrt{(1-\sqrt{\rho})v^2+B_1v+1}=|\sqrt{1-\sqrt{\rho}}\,v-1|$ so the integral in (\ref{s33_xsigmaB1}) will not converge at $v=(1-\sqrt{\rho})^{-1/2}$. In region $D_1$ we also have
\begin{eqnarray}\label{s33_gammaD1}
\gamma(x,\sigma)&=&\frac{\sqrt{2}}{1-\sqrt{\rho}}\exp\bigg[\frac{1+\sqrt{\rho}}{2(1-\sqrt{\rho})}\bigg]e^{x^2/4}x^{-1/2}\Big[(1-\sqrt{\rho})x+1/x+B_1\Big]^{-1/4}\nonumber\\
&&\times\bigg\{\int_0^x\big[(1-\sqrt{\rho})v+1/v+B_1\big]^{-3/2}dv\bigg\}^{-1/2}.
\end{eqnarray}
Thus we have expressed the asymptotic result in (\ref{s33_xsigmapn}) in terms of the solution to (\ref{s33_xsigmaB1}), which involves an integral that is expressible also in terms of incomplete elliptic integrals. 

The region $D_2$ is defined by the conditions $B_1<-2\sqrt{1-\sqrt{\rho}}$, $\eta_x<0$ and here $0<x<\alpha$ where
\begin{equation}\label{s33_alpha}
\alpha=\frac{-B_1-\sqrt{B_1^2-4(1-\sqrt{\rho})}}{2(1-\sqrt{\rho})}.
\end{equation}
In $D_2$, $B_1$ is again to be computed from (\ref{s33_xsigmaB1}), and (\ref{s33_etaD1}) and (\ref{s33_gammaD1}) hold. The curve that separates $D_1$ from $D_2$ is 
\begin{equation}\label{s33_D1D2}
2\sqrt{\rho}\sqrt{1-\sqrt{\rho}}\,\sigma=-2\sqrt{x}+\frac{1}{(1-\sqrt{\rho})^{1/4}}\log\bigg(\frac{1+(1-\sqrt{\rho})^{1/4}\sqrt{x}}{1-(1-\sqrt{\rho})^{1/4}\sqrt{x}}\bigg),
\end{equation}
which passes through $(x,\sigma)=(0,0)$ and has the vertical asymptote $x=(1-\sqrt{\rho})^{-1/2}$ in the $(x,\sigma)$ plane.

The region $D_3$ is defined by the conditions $B_1<-2\sqrt{1-\sqrt{\rho}}$ and $0<x<\alpha$, but now $\eta_x>0$. Then within region $D_3$ we define $B_1=B_1(x,\sigma)$ as the solution to
\begin{equation}\label{s33_B1D3}
\int_x^\alpha\Big[(1-\sqrt{\rho})v+1/v+B_1\Big]^{-1/2}dv+\int_0^\alpha\Big[(1-\sqrt{\rho})v+1/v+B_1\Big]^{-1/2}dv=2\sqrt{\rho}\,\sigma.
\end{equation}
Note that when $x=\alpha$, (\ref{s33_B1D3}) and (\ref{s33_xsigmaB1}) agree, and in fact (\ref{s33_B1D3}) simply gives the smooth continuation of $B_1(x,\sigma)$ as we move from region $D_2$ to $D_3$. In $D_3$ we have
\begin{equation}\label{s33_etaD3}
\eta(x,\sigma)=B_1\sqrt{\rho}\,\sigma+\int_0^x\sqrt{(1-\sqrt{\rho})v+1/v+B_1}\,dv-2\int_0^\alpha\sqrt{(1-\sqrt{\rho})v+1/v+B_1}\,dv
\end{equation}
and
\begin{eqnarray}\label{s33_gammaD3}
\gamma(x,\sigma)&=&\frac{\sqrt{2}}{1-\sqrt{\rho}}\exp\bigg[\frac{1+\sqrt{\rho}}{2(1-\sqrt{\rho})}\bigg]e^{x^2/4}x^{-1/2}\Big[(1-\sqrt{\rho})x+1/x+B_1\Big]^{-1/4}\nonumber\\
&&\times\bigg\{\int_0^x\big[(1-\sqrt{\rho})v+1/v+B_1\big]^{-3/2}dv+\gamma_0\bigg\}^{-1/2}
\end{eqnarray}
where $\gamma_0$ is given below in terms of the complete elliptic integrals of the first kind $K(\cdot)$ and of the second kind $E(\cdot)$:
\begin{equation*}\label{s33_gamma0D3}
\gamma_0=-\frac{2\sqrt{2}\sqrt{-B_1+\sqrt{B_1^2-4(1-\sqrt{\rho})}}}{(1-\sqrt{\rho})\big[B_1^2-4(1-\sqrt{\rho})\big]}\bigg[\sqrt{B_1^2-4(1-\sqrt{\rho})}\,K\Big(\sqrt{{\alpha}/{\beta}}\Big)+B_1E\Big(\sqrt{{\alpha}/{\beta}}\Big)\bigg],
\end{equation*}
where
\begin{equation*}\label{s33_beta}
\beta=\frac{-B_1+\sqrt{B_1^2-4(1-\sqrt{\rho})}}{2(1-\sqrt{\rho})}.
\end{equation*}
The functions $B_1$, $\eta$ and $\gamma$ are all smooth at the curve $\eta_x=0$, i.e., 
\begin{equation}\label{s33_D2D3}
\sqrt{\rho}\,(1-\sqrt{\rho})\sigma=\frac{1}{x}\Big[K\Big(\sqrt{(1-\sqrt{\rho})x^3}\,\Big)-E\Big(\sqrt{(1-\sqrt{\rho})x^3}\,\Big)\Big],
\end{equation}
which separates $D_2$ from $D_3$. This completes the summary of our results on the $(x,\sigma)$ scale. In section 5 we discuss the geometric significance of this curve.

The result for $n=O(1)$ and $t=O(N^{3/4})$ (the $(n,\sigma)$ scale) will involve the functions $B_1$, $\gamma_0$, $\alpha$ and $\beta$ evaluated at $x=0$, so we define $B_1^*(\sigma)=B_1(0,\sigma)$ implicitly as the solution to 
\begin{equation}\label{s33_nsgmb1}
\int_0^{\alpha_*(\sigma)}\frac{dv}{\sqrt{(1-\sqrt{\rho})v+1/v+B_1^*}}=\sqrt{\rho}\,\sigma,
\end{equation}
and we also define
\begin{equation}\label{s33_nsgmalpha}
\alpha_*(\sigma)=\alpha(0,\sigma)=\frac{-B_1^*-\sqrt{(B_1^*)^2-4(1-\sqrt{\rho})}}{2(1-\sqrt{\rho})},
\end{equation}
\begin{equation}\label{s33_nsgmbeta}
\beta_*(\sigma)=\beta(0,\sigma)=\frac{-B_1^*+\sqrt{(B_1^*)^2-4(1-\sqrt{\rho})}}{2(1-\sqrt{\rho})},
\end{equation}
and
\begin{equation}\label{s33_nsgmgamma}
\gamma_*(\sigma)=\gamma_0(0,\sigma)=-\frac{2\sqrt{\rho}\,\sigma}{\sqrt{(B_1^*)^2-4(1-\sqrt{\rho})}}+\frac{8\sqrt{\alpha_*}}{(B_1^*)^2-4(1-\sqrt{\rho})}\,E\Big(\sqrt{\frac{\alpha_*}{\beta_*}}\Big).
\end{equation}
The equation (\ref{s33_nsgmb1}) that defines $B_1^*(\sigma)$ may be written in terms of elliptic integrals, as
\begin{equation}\label{s33_nsgmellip}
\frac{2\sqrt{\beta_*}}{\sqrt{1-\sqrt{\rho}}}\bigg[K\Big(\sqrt{\frac{\alpha_*}{\beta_*}}\Big)-E\Big(\sqrt{\frac{\alpha_*}{\beta_*}}\Big)\bigg]=\sqrt{\rho}\,\sigma,
\end{equation}
and we note that $\sqrt{\alpha_*/\beta_*}=\sqrt{1-\sqrt{\rho}}\,\alpha_*$.
Then for $n=O(1)$, $t=N^{3/4}\sigma=O(N^{3/4})$ we have
\begin{eqnarray}\label{s33_nsigmapn}
p_n(t)&\sim& \rho^{-n/2}e^{-(1-\sqrt{\rho})^2t}\frac{2\sqrt{2\pi}}{1-\sqrt{\rho}}\exp\Big(\frac{\sqrt{\rho}}{1-\sqrt{\rho}}\Big)\big[\gamma_*(\sigma)\big]^{-1/2}N^{-5/8}\nonumber\\
&&\times\exp\bigg\{N^{1/4}\bigg[\sqrt{\rho}\,\sigma B_1^*-2\int_0^{\alpha_*}\sqrt{(1-\sqrt{\rho})v+1/v+B_1^*}\,dv\bigg]\bigg\}\nonumber\\
&&\times\frac{1}{2\pi i}\oint\frac{1}{z^{n+1}(1-z)}\exp\Big(\frac{1}{1-z}\Big)dz.
\end{eqnarray}
For $n\to\infty$ we can use (\ref{s33_int}) to show that (\ref{s33_nsigmapn}) matches to (\ref{s33_xsigmapn}), with $\eta$ given by (\ref{s33_etaD3}) and $\gamma$ by (\ref{s33_gammaD3}), as $x\to 0^+$. Note that $x=0$ is contained within region $D_3$. For $\sigma\to 0$, we can easily show that (\ref{s33_nsigmapn}) leads to precisely the tail behavior of the infinite population model in (\ref{s33_O(1)}). This shows that for the finite population model, (\ref{s33_O(1)}) approximates $p_n(t)$ for $n=O(1)$ only for the time ranges $1\ll t\ll N^{3/4}$, and when $t=O(N^{3/4})$ a fundamentally different result is needed. For $\sigma\to \infty$, (\ref{s33_nsigmapn}) matches to (\ref{s33_ntaupn}).

This completes our summary of the 10 space/time scales where $p_n(t)$ has different behaviors, for $N\to\infty$ and $\rho<1$. The results are derived in section 5, where we also discuss in more detail the asymptotic matching between different scales.

\subsection{Unconditional distribution}

When $\rho>1$ and $N\to\infty$ we shall show that 
\begin{equation}\label{s31_rho>1_p}
p(t)=\frac{\rho}{\rho-1}\,\exp\Big(-\frac{\rho\tau}{\rho-1}\Big)N^{-1}+\mathcal{P}_1(\tau)N^{-2}+O(N^{-3})
\end{equation}
where $\tau=t/N$ and the correction term is
\begin{equation}\label{s31_rho>1_p1}
\mathcal{P}_1(\tau)=\frac{\rho^3}{(\rho-1)^5}\Big(\tau+\frac{1-2\rho}{\rho}\Big)\exp\Big(-\frac{\rho\tau}{\rho-1}\Big)+\frac{\rho^4}{(\rho-1)^5}\exp\Big(-\frac{\rho^2\tau}{\rho-1}\Big).
\end{equation}
This approximation is for large times $t=O(N)$, though it holds also for times $t=O(1)$, and breaks down only for $t=O(N^2)$. Note that for $t=O(N^2)$, $\tau=O(N)$ and then $N^{-1}$ becomes comparable to $\tau N^{-2}$, so that the correction term becomes comparable to the leading term. We also note that $\int_0^\infty\mathcal{P}_1(\tau)d\tau=0$ so that the correction term carries zero mass. Unlike the case $\rho<1$ below, there is only one time scale (the $\tau$ scale) important for understanding $p(t)$. 

We can also interpret (\ref{s31_rho>1_p}) in terms of the spectral expansion in (\ref{s2_p(t)eigen}), with the asymptotic result in (\ref{s2_pasym}). We see that the leading term in (\ref{s31_rho>1_p}) involves only the zeroth eigenvalue $\nu_0$, while (\ref{s31_rho>1_p1}) involves a second exponential, and this corresponds to $\nu_1$ ($\sim\frac{\rho^2}{N(\rho-1)}$ as $N\to\infty$). The term that is linear in $\tau$ in (\ref{s31_rho>1_p1}) arises from the $O(N^{-2})$ correction term in (\ref{s2_nuj>1}) for the zeroth eigenvalue $\nu_0$. By computing higher order terms in (\ref{s2_nuj>1}) we could extend the validity of (\ref{s31_rho>1_p}) to times $\tau=O(N)$ ($t=O(N^2)$) and larger, simply by the ``renormalization"
\begin{equation*}\label{s31_renormal}
\frac{\rho}{\rho-1}\exp\Big(-\frac{\rho\tau}{\rho-1}\Big)\Big[1+\frac{\rho^2\tau}{(\rho-1)^4}\,\frac{1}{N}\Big]\sim\frac{\rho}{\rho-1}\exp\bigg\{-\frac{\rho\tau}{\rho-1}\Big[1-\frac{\rho}{(\rho-1)^3}\,\frac{1}{N}\Big]\bigg\},
\end{equation*}
as from the above we can infer the $O(N^{-2})$ term for $j=0$ in (\ref{s2_nuj>1}). Our discussion suggests that the approximation $p(t)\sim d_0e^{-\nu_0t}$ as $t\to\infty$ actually holds for any time $t$, if $N\to\infty$ and $\rho>1$. Thus we do not need the condition $t\gg O(N)$ for the zeroth eigenvalue to dominate. We discuss in section 7 how to obtain (\ref{s31_rho>1_p}) from the spectral results in \cite{ZHE_O12}.

Next we consider $\rho<1$. Now the expansion of $p(t)$ is different on three time scales: $t=O(1)$, $t=N^{3/4}\sigma=O(N^{3/4})$ and $t=N \tau=O(N)$. For times $t=O(1)$ we have
\begin{equation}\label{s31_O(1)}
p(t)=p(t;M/M/1-\textrm{PS})+O(N^{-1}).
\end{equation}
Thus for times $t=O(1)$, we again approximate the finite population model by the infinite population model, with an $O(N^{-1})$ error term that may be computed by analyzing (\ref{s2_recu}) by transform methods.
Some studies of $p(t;M/M/1-\textrm{PS})$ included Morrison \cite{MOR_R}, who derived integral representations and studied asymptotic properties, and Guillemin and Boyer \cite{GUI}, who used a spectral approach. Conditional moments of this model, i.e. $m_{_J}(n)=\int_0^\infty t^Jp_n(t;M/M/1-\textrm{PS})dt$ were studied in \cite{SEN}, in the limit of $n\to\infty$. The spectrum in the infinite population model is purely continuous, in sharp contrast to the present model, which clearly has a purely discrete spectrum, as the matrix $\mathbf{A}$ in (\ref{s2_A}) is finite.

For $t=N^{3/4}\sigma =O(N^{3/4})$ we obtain
\begin{eqnarray}\label{s31_sigma}
p(t)&\sim &2\sqrt{2\pi}\,\frac{1+\sqrt{\rho}}{1-\sqrt{\rho}}\,\exp\Big(\frac{1+\sqrt{\rho}}{1-\sqrt{\rho}}\Big)N^{-5/8}\frac{1}{\sqrt {\gamma_*(\sigma)}}\,e^{-(1-\sqrt{\rho})^2t}\nonumber\\
&&\times\exp\bigg[N^{1/4}\bigg(\sqrt{\rho}\,\sigma B_1^*(\sigma)-2\int_0^{\alpha_*(\sigma)}\sqrt{(1-\sqrt{\rho})v+1/v+B_1^*(\sigma)}\,dv\bigg)\bigg],\end{eqnarray}
where $\alpha_*(\sigma)$, $\beta_*(\sigma)$ and $\gamma_*(\sigma)$ are given in (\ref{s33_nsgmalpha})--(\ref{s33_nsgmgamma}).
The integral in (\ref{s31_sigma}) may also be expressed as a linear combination of the complete elliptic integrals $K(\cdot)$ and $E(\cdot)$. Thus the expression on the $\sigma$-scale is quite complicated, but it can be simplified in the limits $\sigma\to 0$ and $\sigma\to\infty$, as discussed below.

On the time scale $t=N\tau=O(N)$ we have
\begin{eqnarray}\label{s31_tau}
p(t)&\sim& \frac{16\sqrt{\pi}\,(1+\sqrt{\rho})}{(1-\sqrt{\rho})^{3/8}}\,N^{-5/8}\,e^{-3}\,\exp\Big(\frac{2}{1-\sqrt{\rho}}\Big)\,e^{-(1-\sqrt{\rho})^2N\tau}\,
e^{-2\sqrt{\rho}\sqrt{1-\sqrt{\rho}}\sqrt{N}\tau}\,e^{-\sqrt{\rho}(1-\sqrt{\rho})^{3/4}N^{1/4}\tau}\nonumber\\
&&\times\exp\bigg\{-\frac{8}{3}(1-\sqrt{\rho})^{-1/4}N^{1/4}\bigg\}\exp\bigg\{-\frac{\sqrt{\rho}(22\sqrt{\rho}-3\rho-15)}{16(1-\sqrt{\rho})}\,\tau\bigg\}.
\end{eqnarray}
This is an explicit expression in $\tau$, and in fact shows an exponential dependence on the large time scale $\tau$. On the $\tau$ scale the eigenvalue $\nu_0=\nu_0(N,\rho)$ dominates, and in fact (\ref{s31_tau}) corresponds to $d_0e^{-\nu_0t}$ in (\ref{s2_pasym}), with $\nu_0$ having the expansion in (\ref{s2_nuj<1}), for $N\to\infty$ and $\rho<1$. The factors in (\ref{s31_tau}) that do not involve time correspond to the expansion of $d_0=d_0(N,\rho)$ for $N\to\infty$ and $\rho<1$, and we note that this coefficient is exponentially small, of the order $\exp\big[-O(N^{1/4})\big]$, which is the next to last factor in (\ref{s31_tau}). 

We next discuss the asymptotic matching of the $t$, $\sigma$ and $\tau$ time scales. For $\sigma\to 0$, we can show that 
\begin{equation*}\label{s31_B*}
B_1^*(\sigma)\sim -\Big(\frac{\pi}{2\sqrt{\rho}}\Big)^{2/3}\sigma^{-2/3},\quad \sigma\to 0
\end{equation*}
and then
\begin{equation*}\label{s31_B*exp}
B_1^*\sqrt{\rho}\,\sigma-2\int_0^{\alpha_*}\sqrt{(1-\sqrt{\rho})v+1/v+B_1^*}\,dv\sim\Big(\frac{\pi}{2}\Big)^{3/2}\rho^{1/6}\sigma^{1/3},\quad \sigma\to 0.
\end{equation*}
Also, from (\ref{s33_nsgmgamma}) we obtain
$\gamma_*(\sigma)\sim 3\pi^{-2/3}2^{5/3}\rho^{5/6}\sigma^{5/3}$ as $\sigma\to 0.$
Thus for $\sigma\to 0$, the approximation in (\ref{s31_sigma}), written in terms of $t=N^{3/4}\sigma$, becomes
\begin{equation}\label{s31_p(t)sigma}
p(t)\sim \frac{2^{3/2}\pi^{5/6}(1+\sqrt{\rho})}{3^{1/2}\rho^{5/12}(1-\sqrt{\rho})\,t^{5/6}}\,\exp\Big(\frac{1+\sqrt{\rho}}{1-\sqrt{\rho}}\Big)\exp\Big\{-(1-\sqrt{\rho})^2t-2^{-2/3}3\pi^{2/3}\rho^{1/6}t^{{1}/{3}}\Big\}.
\end{equation}
But the above is also the large $t$ expansion of the infinite population model in (\ref{s31_O(1)}). Thus we have shown that (\ref{s31_p(t)sigma}) holds on intermediate time scales where $t\gg 1$ but $tN^{-3/4}=\sigma\ll 1$. 

Next we consider the asymptotic matching between the $\sigma$ and $\tau$ time scales, where $O(N^{3/4})\ll t\ll O(N)$. Expanding (\ref{s31_tau}) as $\tau\to 0$ corresponds simply to approximating the last exponential factor by 1. For $\sigma\to +\infty$ we can show that $B_1^*\to -2\sqrt{1-\sqrt{\rho}}$ and then 
\begin{equation}\label{s31_B*exp2}
B_1^*\sqrt{\rho}\,\sigma-2\int_0^{\alpha_*}\sqrt{(1-\sqrt{\rho})v+1/v+B_1^*}\,dv=-2\sqrt{\rho}\sqrt{1-\sqrt{\rho}}\,\sigma-\frac{8}{3}(1-\sqrt{\rho})^{-1/4}+o(1),\quad \sigma\to \infty,
\end{equation}
and 
\begin{equation}\label{s31_gamma*2}
\gamma_*(\sigma)\sim 2^{-5}e^4(1-\sqrt{\rho})^{-5/4}e^{2\sqrt{\rho}(1-\sqrt{\rho})^{3/4}\sigma},\quad \sigma\to \infty.
\end{equation}
Using (\ref{s31_B*exp2}) and (\ref{s31_gamma*2}) in (\ref{s31_sigma}) we see that the large $\sigma$ expansion of (\ref{s31_sigma}) agrees with the small $\tau$ expansion of (\ref{s31_tau}). In the matching region between the $\sigma$ and $\tau$ scales we can in fact obtain a new asymptotic approximation to $p(t)$. Consider the scaling 
\begin{equation}\label{s31_sigma*scale}
\sigma=\frac{1}{8}\rho^{-1/2}(1-\sqrt{\rho})^{-3/4}\log N+\sigma_*
\end{equation} 
where $\sigma_*=O(1)$. Then by computing explicitly the $o(1)$ error term in (\ref{s31_B*exp2}), which turns out to be exponentially small for large $\sigma$, we obtain
\begin{eqnarray}\label{s31_sigma*}
p(t)&\sim& \frac{16\sqrt{\pi}(1+\sqrt{\rho})}{(1-\sqrt{\rho})^{3/8}}e^{-3}\exp\Big(\frac{2}{1-\sqrt{\rho}}\Big)N^{-5/8}\exp\Big\{-\frac{8}{3}(1-\sqrt{\rho})^{-1/4}N^{1/4}\Big\}\nonumber\\
&&\times \exp\Big\{-(1-\sqrt{\rho})^2N^{3/4}\sigma_*-2\sqrt{\rho}\sqrt{1-\sqrt{\rho}}N^{1/4}\sigma_*-\sqrt{\rho}(1-\sqrt{\rho})^{3/4}\sigma_*\Big\}\nonumber\\
&&\times\exp\Big\{32e^{-4}e^{-2\sqrt{\rho}(1-\sqrt{\rho})^{3/4}\sigma_*}\Big\},
\end{eqnarray}
which applies for $\sigma_*=O(1)$ and gives a double exponential limit law, in view of the last factor in (\ref{s31_sigma*}). We have thus identified more precisely the exact time scale when the eigenvalue $\nu_0$ dominates the behavior of $p(t)$, namely when $\sigma_*\to\infty$ or $\sigma-2^{-3}\rho^{-1/2}(1-\sqrt{\rho})^{-3/4}\log N\to \infty$. By writing the last factor in (\ref{s31_sigma*}) as the series
\begin{equation*}\label{s31_sigma*last}
\sum_{j=0}^\infty\frac{32^je^{-4j}}{j!}\exp\big[-2j\sqrt{\rho}(1-\sqrt{\rho})^{3/4}\sigma_*\big]
\end{equation*}
we can interpret the $j^\textrm{th}$ term in the series as corresponding to the $j^\textrm{th}$ eigenvalue $\nu_j$ (cf. (\ref{s2_nuj<1})). Thus for times $\sigma_*=O(1)$ all of the eigenvalues with $j=O(1)$ contribute roughly equally to the expansion of $p(t)$. In section 7 we shall see how (\ref{s31_sigma*}) follows from the spectral results we obtained in \cite{ZHE_O12}. However, for $\sigma=O(1)$ (then we would have $\sigma_*\to -\infty$) it seems difficult if not impossible to obtain (\ref{s31_sigma}) from the spectral results. At the very least we would need to understand the eigenvalues/eigenvectors of (\ref{s2_A}) of large index $j=O(N^{1/4})$, and then the expansion in (\ref{s2_nuj<1}) breaks down.

We have thus shown that for $\rho<1$ and $N\to\infty$ the asymptotics of $p(t)$ are quite intricate, with the time range $t=O(N^{3/4})$ leading to a complicated transition from the tail behavior in (\ref{s31_p(t)sigma}) for the infinite population model, to the purely exponential behavior in (\ref{s31_tau}). For times slightly larger than $t=O(N^{3/4})$ (cf. (\ref{s31_sigma*scale})) we obtain the double exponential behavior in (\ref{s31_sigma*}).

\section{Asymptotic analysis for $p_n(t)$ if $\rho>1$}

We briefly derive the asymptotic approximations for the sojourn time density with $\rho>1$, recalling that asymptotic expansions are different for the three scales (i) $n,\,t=O(N)$, (ii) $n=O(1)$, $t=O(N)$ and (iii) $n,\,t=O(1)$.

For the large time/space scale where $n,\,t=O(N)$, we use (\ref{s32_scales}) and rewrite $p_n(t)$ as 
\begin{equation}\label{s4_pn(t)}
p_n(t)=\frac{1}{n+1}Q_{n+1}(t)
\end{equation}
with the initial condition $Q_{n+1}(0)=1$. Then $Q_{n+1}(t)$ will satisfy
\begin{equation}\label{s4_Qn(t)}
\frac{d}{dt}Q_{n+1}(t)=\rho\Big(1-\frac{n+1}{N}\Big)\frac{n+1}{n+2}\,Q_{n+2}(t)+Q_{n}(t)-\Big[\rho\Big(1-\frac{n+1}{N}\Big)+1\Big]Q_{n+1}(t).
\end{equation}
By introducing 
\begin{equation}\label{s4_xi1}
\xi_1=\frac{n+1}{N}=\xi+\frac{1}{N},
\end{equation}
we set
\begin{equation}\label{s4_Q}
Q_{n+1}(t)=A(\xi_1,\tau)+\frac{1}{N}A(\xi_1,\tau)C(\xi_1,\tau)+O\Big(\frac{1}{N^2}\Big).
\end{equation}
Then using (\ref{s4_Q}) in (\ref{s4_Qn(t)}) and noting that changing $n$ to $n\pm 1$ corresponds to changing $\xi_1$ to $\xi_1\pm 1/N$, we obtain the following PDE for $A(\xi_1,\tau)$
\begin{equation}\label{s4_APDE}
A_\tau+\big[1-\rho(1-\xi_1)\big]A_{\xi_1}=-\rho\frac{1-\xi_1}{\xi_1}A,
\end{equation}
with the initial condition $A(\xi_1,0)=1$, and for $C(\xi_1, \tau)$ the PDE
\begin{equation}\label{s4_CPDE}
C_{\tau}+\big[1-\rho(1-\xi_1)\big]C_{\xi_1}=\rho\frac{1-\xi_1}{\xi_1^2}-\rho\frac{1-\xi_1}{\xi_1}\frac{A_{\xi_1}}{A}+\frac{1}{2}\big[\rho(1-\xi_1)+1\big]\frac{A_{\xi_1\xi_1}}{A}
\end{equation}
with the initial condition $C(\xi_1,0)=0$. The PDEs (\ref{s4_APDE}) and (\ref{s4_CPDE}) can be solved by the method of characteristics, thus $A(\xi_1,\tau)$ is given by 
\begin{equation*}\label{s4_A}
A(\xi_1,\tau)=\xi_1^{\frac{\rho}{\rho-1}}\exp\Big(-\frac{\rho\tau}{\rho-1}\Big)\Big[\Big(\xi_1-\frac{\rho-1}{\rho}\Big)e^{-\rho\tau}+\frac
{\rho-1}{\rho}\Big]^{-\frac{\rho}{\rho-1}}
\end{equation*}
and $C(\xi_1,\tau)$ is given by (\ref{s32_Cxitau}) with $\xi$ replaced by $\xi_1$. Here we omit the detailed derivation.
We note that by (\ref{s4_xi1}), a Taylor expansion about $\xi_1=\xi$ yields
$$A(\xi_1,\tau)=A(\xi,\tau)+\frac{1}{N}A_{\xi_1}(\xi,\tau)+O\Big(\frac{1}{N^2}\Big)$$
and $C(\xi_1,\tau)=C(\xi,\tau)+O(1/N)$. Thus $Q_{n+1}(t)$ can be expanded as
\begin{equation}\label{s4_Q2}
Q_{n+1}(t)=A(\xi,\tau)+\frac{1}{N}\Big[A_{\xi_1}(\xi,\tau)+A(\xi,\tau)C(\xi,\tau)\Big]+O\Big(\frac{1}{N^2}\Big).
\end{equation}
Using (\ref{s4_xi1}) and (\ref{s4_Q2}) in (\ref{s4_pn(t)}) and letting $P(\xi,\tau)=A(\xi,\tau)/\xi$, we obtain (\ref{s32_xitau}) with (\ref{s32_Pxitau}) and (\ref{s32_Pxitau1}).

We next consider the scale $n=O(1)$, $t=N\tau=O(N)$. From (\ref{s32_xitau}) and (\ref{s32_Pxitau}) it follows that
$$p_n(t)\sim N^{-1}\xi^{\frac{1}{\rho-1}}\Big(\frac{\rho}{\rho-1}\Big)^{\frac{\rho}{\rho-1}}(1-e^{-\rho\tau})^{-\frac{\rho}{\rho-1}}\exp\Big(-\frac{\rho\tau}{\rho-1}\Big)$$
as $\xi\to 0$. Then we expand $p_n(t)$ in the form 
\begin{equation*}\label{s4_pn(t)2}
p_n(t)\sim N^{-\frac{\rho}{\rho-1}}\big(1-e^{-\rho\tau}\big)^{-\frac{\rho}{\rho-1}}\exp\Big(-\frac{\rho\tau}{\rho-1}\Big)\mathcal{P}(n)
\end{equation*}
and use it in (\ref{s2_recu}), to find that $\mathcal{P}(n)$ satisfies the following limiting difference equation
\begin{equation}\label{s4_mathcalP}
\rho\mathcal{P}(n+1)-(\rho+1)\mathcal{P}(n)+\frac{n}{n+1}\mathcal{P}(n-1)=0,
\end{equation}
with the matching condition
$$\mathcal{P}(n)\sim n^{\frac{1}{\rho-1}}\Big(\frac{\rho}{\rho-1}\Big)^{\frac{\rho}{\rho-1}},\quad n\to\infty,$$
and with $\mathcal{P}(-1)$ finite (thus (\ref{s4_mathcalP}) holds for all $0\le n\le N-1$). Solving (\ref{s4_mathcalP}) using generating functions or contour integrals leads to the formula in (\ref{s32_O(N)}).

For the scale $n,\,t=O(1)$ we can approximate the finite population model by the infinite population model which we already discussed in section 3.1. This completes the analysis of $p_n(t)$ for $\rho>1$.

\section{Asymptotic analysis for $p_n(t)$ if $\rho<1$}

We proceed to derive the asymptotic approximations for the conditional sojourn time density with $\rho<1$. In region $R_1$ the expansion (\ref{s32_xitau}) holds and the derivation is the same as in section 4.

\subsection{Region $R_2$}

We first consider region $R_2$ and assume an expansion of the form in (\ref{s33_R2pn}), where $\phi$, $\phi^{(1)}$ and $K$ are to be determined. Then from (\ref{s2_recu}) we obtain the following PDEs for $\phi(\xi,\tau)$, $\phi^{(1)}(\xi,\tau)$ and $K(\xi,\tau)$:
\begin{equation}\label{s51_phiPDE}
\phi_\tau-\Big[\rho(1-\xi)e^{\phi_\xi}+e^{-\phi_\xi}\Big]=-\rho(1-\xi)-1,
\end{equation}
\begin{equation}\label{s51_phi(1)PDE}
\phi^{(1)}_\tau-\Big[\rho(1-\xi)e^{\phi_\xi}-e^{-\phi_\xi}\Big]\phi^{(1)}_\xi=0,
\end{equation}
and
\begin{equation}\label{s51_KPDE}
K_\tau-\Big[\rho(1-\xi)e^{\phi_\xi}-e^{-\phi_\xi}\Big]K_\xi=\bigg\{\rho\Big[\frac{1-\xi}{2}\phi_{\xi\xi}-1\Big]e^{\phi_\xi}+\Big(\frac{\phi_{\xi\xi}}{2}-\frac{1}{\xi}\Big)e^{-\phi_\xi}+\rho\bigg\}K.
\end{equation}
The PDE (\ref{s51_phiPDE}) can be solved by the method of characteristics. The characteristic equations are the five ODEs
\begin{equation}\label{s51_R2rayw}
\frac{d\tau}{dw}=1,\quad \frac{d\xi}{dw}=e^{-\phi_\xi}-\rho(1-\xi)e^{\phi_\xi},
\end{equation}
\begin{equation}\label{s51_R2rayw2}
\frac{d\phi}{dw}=\phi_\tau-\phi_\xi\Big[\rho(1-\xi)e^{\phi_\xi}-e^{-\phi_\xi}\Big],
\end{equation}
and 
\begin{equation}\label{s51_R2rayw3}
\frac{d\phi_\tau}{dw}=0,\quad\frac{d\phi_\xi}{dw}=\rho\Big(1-e^{\phi_\xi}\Big).
\end{equation}
Here $w$ is a parameter that increases along a characteristic curve.
We shall use in $R_2$ the characteristics, also called ``rays", that start from the point $(\xi,\tau)=(0,0)$; these are given by 
\begin{equation}\label{s51_R2ray}
\xi=(1-e^{-\rho\tau})\Big(\frac{e^{\rho\tau}-B}{\rho}+\frac{e^{\rho\tau}}{B-1}\Big),
\end{equation}
where $B$ is a constant on each ray. Solving (\ref{s51_R2ray}) for $B=B(\xi,\tau)$, we have the expression in (\ref{s33_R2B}). Then, solving (\ref{s51_R2rayw2}) with (\ref{s51_R2rayw}) and (\ref{s51_R2rayw3}) leads to $\phi(\xi,\tau)$ in (\ref{s33_R2phi}).

Using the same ray variable $w$ and (\ref{s51_R2rayw}), we can rewrite (\ref{s51_phi(1)PDE}) and (\ref{s51_KPDE}), respectively, as $d\phi^{(1)}/dw=0$ and 
\begin{equation}\label{s51_K/wPDE}
\frac{dK}{dw}=\bigg\{\rho\Big[\frac{1-\xi}{2}\phi_{\xi\xi}-1\Big]e^{\phi_\xi}+\Big(\frac{\phi_{\xi\xi}}{2}-\frac{1}{\xi}\Big)e^{-\phi_\xi}+\rho\bigg\}K.
\end{equation}
Then $\phi^{(1)}(\xi,\tau)$ is a constant along each ray and we set $\phi^{(1)}(\xi,\tau)\equiv \phi_1(B)$. We also get (\ref{s33_R2K}) from (\ref{s51_K/wPDE}) after some calculation and simplification. The functions $\phi_1(B)$ and $K_0(B)$ cannot be determined directly from the above PDEs, so we proceed to use asymptotic matching between $R_2$ and the expansion valid for $n,\,t=O(1)$, cf. (\ref{s32_O(1)}) or (\ref{s32_p^(0)}). We let $\xi$ and $\tau\to 0$ but with $n/t=\xi/\tau=O(1)$ in (\ref{s33_R2phi}), (\ref{s33_R2B}) and (\ref{s33_R2K}), and match with the results in Theorem 2.2 item 3 in \cite{ZHE_O10}, which corresponds to the expansion of (\ref{s32_O(1)}) (or (\ref{s32_p^(0)})) for $n$ and $t\to\infty$ with $t/n$ fixed and $t/n\in(0,(1-\rho)^{-1})$. This ultimately yields (\ref{s33_R2phi(1)}) and (\ref{s33_R2K0}). Here we omit the detailed derivation, but only point out that $B$ in (\ref{s33_R2B}) is asymptotically given by, for $\xi,\,\tau\to 0$,
$$B= 1-\frac{1}{2}\Big(\frac{n}{t}\Big)-\frac{1}{2}\sqrt{\Big(\frac{n}{t}\Big)^2+4\rho}+o(1).$$
This completes the analysis of region $R_2$.

\subsection{Region $R_3$}

We next consider region $R_3$ and assume that $p_n(t)$ has an asymptotic expansion of the form
\begin{equation}\label{s52_R3pn}
p_n(t)\sim N^{\delta_1}e^{N\psi (\xi,\tau)}e^{\sqrt{N}\psi^{(1)}(\xi,\tau)}e^{N^{1/4}\psi^{(2)}(\xi,\tau)}L(\xi,\tau).
\end{equation}
Using (\ref{s52_R3pn}) in (\ref{s2_recu}) yields the PDEs
\begin{equation}\label{s52_psiPDE}
\psi_\tau-\Big[\rho(1-\xi)e^{\psi_\xi}+e^{-\psi_\xi}\Big]=-\rho(1-\xi)-1,
\end{equation}
\begin{equation}\label{s52_psi(1)PDE}
\psi^{(1)}_\tau-\Big[\rho(1-\xi)e^{\psi_\xi}-e^{-\psi_\xi}\Big]\psi^{(1)}_\xi=0,
\end{equation}
\begin{equation}\label{s52_psi(2)PDE}
\psi^{(2)}_\tau-\Big[\rho(1-\xi)e^{\psi_\xi}-e^{-\psi_\xi}\Big]\psi^{(2)}_\xi=0,
\end{equation}
and
\begin{equation}\label{s52_LPDE}
L_\tau-\Big[\rho(1-\xi)e^{\psi_\xi}-e^{-\psi_\xi}\Big]L_\xi=\bigg\{\frac{1}{2}\Big[\rho(1-\xi)e^{\psi_\xi}+e^{-\psi_\xi}\Big]\Big[\psi_{\xi\xi}+\big(\psi_{\xi}^{(1)}\big)^2\Big]-\rho e^{\psi_\xi}-\frac{1}{\xi}e^{-\psi_\xi}+\rho\bigg\}L.
\end{equation}
Using again the method of characteristics and the ray variable $w$, we obtain (\ref{s51_R2rayw})--(\ref{s51_R2rayw3}), with $\phi$ replaced by $\psi$.
Now we must use rays that are tangent to the line $\xi=0$ at $\tau=0$. Omitting the details, this corresponds to the following one parameter family of curves, with $B_0$ indexing the family:
\begin{equation}\label{s52_R3ray}
\xi=\frac{(1-\sqrt{\rho})^2}{\rho B_0}\,e^{\rho\tau}+\frac{B_0}{\rho}\,e^{-\rho\tau}-\frac{2(1-\sqrt{\rho})}{\rho}.
\end{equation}
The rays fill up the region $R_3$ for $0<B_0<1-\sqrt{\rho}$. 
Solving (\ref{s52_R3ray}) for $B_0$ yields (\ref{s33_R3B0}). Then (\ref{s33_R3psi}) is obtained by solving (\ref{s52_psiPDE}) subject also to $\psi(0,0)=0$. The PDEs (\ref{s52_psi(1)PDE}) and (\ref{s52_psi(2)PDE}) imply that $\psi^{(1)}$ and $\psi^{(2)}$ are constant along the rays, so we set $\psi^{(1)}(\xi,\tau)=\psi_1(B_0)$ and $\psi^{(2)}(\xi,\tau)=\psi_2(B_0)$. $L(\xi,\tau)$ can be determined from (\ref{s52_LPDE}) up to a multiplicative constant along the rays, i.e., the general solution is $L(\xi,\tau)=L_0(B_0)J(\xi)$, where $J(\xi)$ is given in (\ref{s33_R3L}). $\psi_1(B_0)$, $\psi_2(B_0)$, $L_0(B_0)$,  as well as the constant $\delta_1$ in (\ref{s52_R3pn}), will be determined in the next section by matching with the boundary layers.

\subsection{Boundary layers}

We consider boundary layer regions where $\xi\to 0$, so that $n=o(N)$. For $n,\,t=O(1)$ the result in (\ref{s32_O(1)}) holds and it has been discussed in section 3.1.

We first consider the $(x,\sigma)$ scale, where $n=\sqrt{N}\,x=O(\sqrt{N})$ and $t=N^{3/4}\sigma=O(N^{3/4})$. Using the expansions in regions $R_2$ and $R_3$, if we let $\tau=\sigma N^{-1/4}$ and $\xi=xN^{-1/2}$ in (\ref{s33_R2phi}) and (\ref{s33_R3psi}), it follows that
\begin{equation*}
e^{N\phi(\xi,\tau)}= \rho^{-n/2}\,\exp\Big\{-(1-\sqrt{\rho})^2\,t+O(N^{1/4})\Big\}
\end{equation*}
and
\begin{equation*}
e^{N\psi(\xi,\tau)}= \rho^{-n/2}\,\exp\Big\{-(1-\sqrt{\rho})^2\,t+O(N^{1/4})\Big\}.
\end{equation*}
Then we assume the following asymptotic expansion on the $(x,\sigma)$ scale:
\begin{equation}\label{s53_xsgmpnt}
p_n(t)\sim\rho^{-n/2}\,e^{-(1-\sqrt{\rho})^2\,t}N^{\delta_2}e^{N^{1/4}\eta(x,\sigma)}\gamma(x,\sigma).
\end{equation}
Using (\ref{s53_xsgmpnt}) in (\ref{s2_recu}) yields at the first two orders the PDEs
\begin{equation}\label{s53_etaPDE}
\eta_\sigma-\sqrt{\rho}\,\eta^2_x=\big(\rho-\sqrt{\rho}\big)x-\frac{\sqrt{\rho}}{x}
\end{equation}
and
\begin{equation}\label{s53_gammaPDE}
\gamma_\sigma-2\sqrt{\rho}\,\eta_x\gamma_x=\sqrt{\rho}\Big(\eta_{xx}-x\eta_x+\frac{1}{x}\eta_x\Big)\gamma.
\end{equation}
We use the method of characteristics to solve (\ref{s53_etaPDE}), with all of the rays starting from the origin $(x,\sigma)=(0,0)$. This is necessary to be able to match to the $R_2$ and $R_3$ expansions. We find that the geometry of the rays naturally defines three regions in the $(x,\sigma)$ plane (as shown in Figure \ref{figure:3}). The first region corresponds to $\eta_x<0$, where the rays $B_1\geq -2\sqrt{1-\sqrt{\rho}}$ satisfy (\ref{s33_xsigmaB1}). When $B_1=-2\sqrt{1-\sqrt{\rho}}$, we have the curve that separates $D_1$ from $D_2$ in (\ref{s33_D1D2}), which is the ray denoted by the dashed curve in Figure \ref{figure:3}. The second region corresponds to $\eta_x<0$ and the rays $B_1<-2\sqrt{1-\sqrt{\rho}}$ satisfy (\ref{s33_xsigmaB1}) with $0<x<\alpha$, where $\alpha$ is defined in (\ref{s33_alpha}). If we let $x\to \alpha$ in (\ref{s33_xsigmaB1}), we have $\eta_x=0$, which is the curve in (\ref{s33_D2D3}). It is not a ray, and is denoted by the dotted curve in Figure \ref{figure:3}. This curve corresponds to the locus of the maximum values of $x$ achieved along those rays that start from $(0,0)$ and return to $x=0$ at some later $\sigma>0$. The third region corresponds to $\eta_x>0$ and  $B_1<-2\sqrt{1-\sqrt{\rho}}$ satisfying (\ref{s33_B1D3}) with $0<x<\alpha$. In all of these three regions, $\eta_\sigma=\sqrt{\rho}\,B_1$ is a constant along a ray. Thus, (\ref{s33_etaD1}) and (\ref{s33_etaD3}) are obtained from the corresponding regions $D_j$.

The function $\gamma(x,\sigma)$ can also be determined from (\ref{s53_gammaPDE}) by the method of characteristics, up to a multiplicative function of $B_1$, so we write
\begin{equation}\label{s53_gamma1}
\gamma(x,\sigma)=\frac{\gamma_1(B_1)e^{x^2/4}}{x^{1/2}\Big[(1-\sqrt{\rho})x+1/x+B_1\Big]^{1/4}\sqrt{\int_0^x\big[(1-\sqrt{\rho})v+1/v+B_1\big]^{-3/2}dv}}\quad \mathrm{in }\; D_1,
\end{equation}
\begin{equation*}
\gamma(x,\sigma)=\frac{\gamma_2(B_1)e^{x^2/4}}{x^{1/2}\Big[(1-\sqrt{\rho})x+1/x+B_1\Big]^{1/4}\sqrt{\int_0^x\big[(1-\sqrt{\rho})v+1/v+B_1\big]^{-3/2}dv}}\quad \mathrm{in }\; D_2,
\end{equation*}
and 
\begin{equation*}
\gamma(x,\sigma)=\frac{\gamma_3(B_1)e^{x^2/4}}{x^{1/2}\Big[(1-\sqrt{\rho})x+1/x+B_1\Big]^{1/4}\sqrt{\int_0^x\big[(1-\sqrt{\rho})v+1/v+B_1\big]^{-3/2}dv+\gamma_0}}\quad \mathrm{in }\; D_3,
\end{equation*}
where $\gamma_1(B_1)$, $\gamma_2(B_1)$ and $\gamma_3(B_1)$ will be determined by matching with results for $n,\,t=O(1)$. We let $x,\,\sigma\to 0$ but keep $a=n\,t^{-2/3}=x\,\sigma^{-2/3}=O(1)$ fixed and match with the results in Theorem 2.2 item 4 in \cite{ZHE_O10}. In region $D_1$, we have $B_1\to+\infty$ with $B_1=O(1/x)$. Then (\ref{s33_xsigmaB1}) is asymptotically approximated as
\begin{eqnarray}\label{s53_xsgmRay1}
2\sqrt{\rho}\,\sigma
&=&\int_0^x\big[(1-\sqrt{\rho})v+1/v+B_1\big]^{-1/2}dv\nonumber\\
&\sim&\int_0^x\frac{\sqrt{v}}{\sqrt{B_1v+1}}dv\nonumber\\
&=& B_1^{-3/2}\Big[\sqrt{B_1x(B_1x+1)}-\sinh^{-1}\Big(\sqrt{B_1x}\Big)\Big].
\end{eqnarray}
Setting $B_1x=A$, we notice that (\ref{s53_xsgmRay1}) is equivalent to the rays in (2.19) in \cite{ZHE_O10}. Then we can easily verify that $N^{1/4}\eta(x,\sigma)\sim\Phi(n,t)$, where $\Phi(n,t)$ is given in (2.18) in \cite{ZHE_O10}. The integral in (\ref{s53_gamma1}) is for $x\to 0$ asymptotically given by
\begin{equation}\label{s53_xsgmint1}
\int_0^x\big[(1-\sqrt{\rho})v+1/v+B_1\big]^{-3/2}dv \sim  \int_0^x\frac{dv}{\big[B_1+1/v\big]^{3/2}}=\frac{2x^{5/2}a^{-3/2}}{A\sqrt{1+A}}\Big[3\sqrt{\rho}\sqrt{1+A}-a^{3/2}\Big].
\end{equation}
Thus, from (\ref{s53_gamma1}) and (\ref{s53_xsgmint1}) we have 
\begin{equation}\label{s53_xsgmNgam}
N^{\delta_2}\gamma(x,\sigma)\sim N^{\delta_2+3/4}\frac{\gamma_1(B_1)\sqrt{A}}{\sqrt{2}\,a^{3/4}\big(3\sqrt{\rho}\sqrt{1+A}-a^{3/2}\big)^{1/2}t}.
\end{equation}
Comparing (\ref{s53_xsgmNgam}) with $\Lambda(n,t)$ in (2.17) in \cite{ZHE_O10}, we obtain $\delta_2=-3/4$ and 
\begin{equation}\label{s53_xsgmconstant}
\gamma_1(B_1)=\frac{\sqrt{2}}{1-\sqrt{\rho}}\exp\bigg[\frac{1+\sqrt{\rho}}{2(1-\sqrt{\rho})}\bigg].
\end{equation}
Then (\ref{s53_gamma1}) and (\ref{s53_xsgmconstant}) lead to (\ref{s33_gammaD1}). Using similar matching arguments for regions $D_2$ and $D_3$, along with corresponding results in \cite{ZHE_O10}, we find that $\gamma_2(B_1)=\gamma_3(B_1)=\gamma_1(B_1)$ and this completes the derivation of the $(x,\sigma)$ scale.

We next consider the $(n,\sigma)$ scale. From the results in region $D_3$ on the $(x,\sigma)$ scale, if we let $x\to 0$, it follows that
$$\eta(x,\sigma)= \sqrt{\rho}\,\sigma B_1^*-2\int_0^{\alpha_*}\sqrt{(1-\sqrt{\rho})v+1/v+B_1^*}\,dv+O(\sqrt{x}),$$
and
$\gamma(x,\sigma)\sim \gamma_1(B_1)x^{-1/4}\big[\gamma_*(\sigma)\big]^{-1/2},$ where $B_1^*$, $\alpha_*$ and $\gamma_*(\sigma)$ are given by (\ref{s33_nsgmb1})--(\ref{s33_nsgmgamma}). Thus, we assume that $p_n(t)$ has an expansion in the following form
\begin{equation}\label{s53_nsgmpnt}
p_n(t)\sim \frac{N^{\delta_3}\mathcal{Q}(n)}{\sqrt{\gamma_*(\sigma)}}\rho^{-n/2}e^{-(1-\sqrt{\rho})^2t}\exp\bigg\{N^{1/4}\bigg[\sqrt{\rho}\,\sigma B_1^*-2\int_0^{\alpha_*}\sqrt{(1-\sqrt{\rho})v+\frac{1}{v}+B_1^*}\,dv\bigg]\bigg\}.
\end{equation}
Using (\ref{s53_nsgmpnt}) in (\ref{s2_recu}), we obtain the following limiting difference equation for $\mathcal{Q}(n)$:
\begin{equation}\label{s53_nsgmDEQ}
\mathcal{Q}(n+1)+\frac{n}{n+1}\mathcal{Q}(n-1)-2\mathcal{Q}(n)=0
\end{equation}
with $\mathcal{Q}(-1)$ finite. Solving (\ref{s53_nsgmDEQ}) yields 
\begin{equation}\label{s53_nsgmQ}
\mathcal{Q}(n)=\frac{\mathcal{Q}(0)e^{-1}}{2\pi i}\oint\frac{1}{z^{n+1}(1-z)}\exp\Big(\frac{1}{1-z}\Big)dz.
\end{equation}
Here the integral is over a small loop about $z=0$. $\mathcal{Q}(0)$ is a constant and can be determined by matching with the results on the $(x,\sigma)$ scale, using (\ref{s33_int}), which leads to 
\begin{equation}\label{s53_nsgmQ0}
\mathcal{Q}(0)=\frac{2\sqrt{2\pi}}{1-\sqrt{\rho}}\exp\bigg(\frac{1}{1-\sqrt{\rho}}\bigg).
\end{equation}
The matching also yields $\delta_3=-5/8$. Using (\ref{s53_nsgmQ}) and (\ref{s53_nsgmQ0}) in (\ref{s53_nsgmpnt}), we have established (\ref{s33_nsigmapn}).

In the rest of the section, we consider $t=N\tau=O(\tau)$ and obtain the undetermined functions $\psi_1(B_0)$, $\psi_2(B_0)$, $L_0(B_0)$, and the constant $\delta_1$ in region $R_3$. 

We found in \cite{ZHE_O12} that the eigenvectors in the spectral representation in (\ref{s2_pn(t)eigen}) are concentrated on the $y$ scale, where
\begin{equation*}\label{s53_yscale}
n=\frac{\sqrt{N}}{\sqrt{1-\sqrt{\rho}}}+N^{3/8}y,\quad y=O(1).
\end{equation*}
In the asymptotic analysis of $p_n(t)$ in this paper, the $(y,\tau)$ scale is contained within the $(x,\tau)$ scale, but its analysis is crucial for deriving the undetermined functions in $R_3$. Thus we consider the $(y,\tau)$ scale first and we assume that $p_n(t)$ has the asymptotic expansion in the following form:
\begin{eqnarray}\label{s53_ypnt}
p_n(t)&\sim &N^{\delta^*}\rho^{-n/2}e^{-(1-\sqrt{\rho})^2t}e^{\sqrt{N}\Psi^{(1)}(\tau)}e^{N^{1/4}\Psi^{(2)}(\tau)}\nonumber\\
&&\times\Big[Z(y,\tau)+N^{-1/8}Z^{(1)}(y,\tau)+N^{-1/4}Z^{(2)}(y,\tau)\Big].
\end{eqnarray}
Using (\ref{s53_ypnt}) in (\ref{s2_recu}), we obtain to leading order the limiting ODE $\Psi^{(1)}_\tau=-2\sqrt{\rho}\sqrt{1-\sqrt{\rho}}$, whose solution is
\begin{equation*}\label{s53_yPsi(1)}
\Psi^{(1)}(\tau)=-2\sqrt{\rho}\sqrt{1-\sqrt{\rho}}\,\tau+d_1^*,
\end{equation*}
where $d_1^*$ is a constant. At the next order we obtain the following PDE for $Z(y,\tau)$
\begin{equation}\label{s53_yZPDE}
Z_{yy}+\bigg[-\frac{\Psi^{(2)}_\tau}{\sqrt{\rho}}-(1-\sqrt{\rho})^{3/2}y^2\bigg]Z=0.
\end{equation}
This is a Hermite or parabolic cylinder equation.
The only acceptable solution (which is positive and decays as $y\to\pm\infty$) corresponds to $\Psi^{(2)}_\tau=-\sqrt{\rho}(1-\sqrt{\rho})^{3/4}$, from which we have
\begin{equation*}\label{s53_yPsi(2)}
\Psi^{(2)}(\tau)=-\sqrt{\rho}(1-\sqrt{\rho})^{3/4}\,\tau+d_2^*,
\end{equation*}
where $d_2^*$ is another constant.
Then (\ref{s53_yZPDE}) has the solution
\begin{equation}\label{s53_yZ}
Z(y,\tau)=u(\tau)\exp\bigg\{-\frac{(1-\sqrt{\rho})^{3/4}}{2}y^2\bigg\},
\end{equation}
where $u(\tau)$ is also to be determined. 
Letting $\mathcal{L}$ be the differential operator $\mathcal{L}\{f(y,\tau)\}=f_{yy}+\big[(1-\sqrt{\rho})^{3/4}-(1-\sqrt{\rho})^{3/2}y^2\big]f$, the PDEs for $Z^{(1)}$ and $Z^{(2)}$ in (\ref{s53_ypnt}) turn out to be
\begin{equation*}\label{s53_yZ1PDE}
\mathcal{L}\{Z^{(1)}\}=\Big[-(1-\sqrt{\rho})^2y^3-\sqrt{\rho}(1-\sqrt{\rho})^{1/4}y\Big]Z
\end{equation*}
and
\begin{equation}\label{s53_yZ2PDE}
\mathcal{L}\{Z^{(2)}\}=\frac{Z_\tau}{\sqrt{\rho}}+(2-\sqrt{\rho})yZ_y+(1-\sqrt{\rho})^{5/2}y^2Z+\frac{\sqrt{\rho}}{\sqrt{1-\sqrt{\rho}}}Z^{(1)}_y-(1-\sqrt{\rho})^2y^3Z^{(1)}.
\end{equation}
To compute $u(\tau)$ in (\ref{s53_yZ}), we solve for $Z^{(1)}$, use the result in (\ref{s53_yZ2PDE}), and then use a solvability condition for (\ref{s53_yZ2PDE}). We omit the detailed derivation. This ultimately leads to a simple ODE for $u(\tau)$ whose general solution is 
\begin{equation*}\label{s53_yu}
u(\tau)=d_3^*\exp\bigg\{\frac{\sqrt{\rho}(15-22\sqrt{\rho}+3\rho)}{16(1-\sqrt{\rho})}\,\tau\bigg\}.
\end{equation*}
Here $d_3^*$ is also a constant. We shall determine $\delta^*$ and $d_j^*$ (for $j=1,2,3$) later.

We next consider $n=O(1)$ and $t=N\tau=O(N)$. The expansion on the $(y,\tau)$ scale in (\ref{s53_ypnt}) suggests that we assume that $p_n(t)$ has the following asymptotic expansion:
\begin{eqnarray}\label{s53_ntaupnt}
p_n(t)&\sim& N^{\delta_4}\rho^{-n/2}\,\mathcal{R}_n\,e^{-(1-\sqrt{\rho})^2t}\exp\Big\{\sqrt{N}\big(-2\sqrt{\rho}\sqrt{1-\sqrt{\rho}}\,\tau+d_1\big)\Big\}\nonumber\\
&&\times\exp\Big\{N^{1/4}\big(-\sqrt{\rho}(1-\sqrt{\rho})^{3/4}\,\tau+d_2\big)\Big\}\exp\bigg\{\frac{\sqrt{\rho}(15-22\sqrt{\rho}+3\rho)}{16(1-\sqrt{\rho})}\,\tau\bigg\}.
\end{eqnarray}
Here $\delta_4$, $d_1$, and $d_2$ are constants to be determined later.
Similarly to the analysis of the $(n,\sigma)$ scale, we use (\ref{s53_ntaupnt}) in (\ref{s2_recu}) and find that $\mathcal{R}_{n+1}+\frac{n}{n+1}\mathcal{R}_{n-1}-2\mathcal{R}_n=0$, which leads to the contour integral 
\begin{equation}\label{s53_ntauR}
\mathcal{R}_n=\frac{\mathcal{R}_0\,e^{-1}}{2\pi i}\oint\frac{1}{z^{n+1}(1-z)}\exp\Big(\frac{1}{1-z}\Big)dz.
\end{equation}
We determine $d_1$, $d_2$ and $\mathcal{R}_0$ by matching (\ref{s53_ntaupnt}) with the $(n,\sigma)$ scale, for $\tau\to 0$ but $\sigma=N^{1/4}\tau\to\infty$. Letting $\sigma\to\infty$ in (\ref{s33_nsgmellip}) with the help of (\ref{s33_nsgmalpha}) and (\ref{s33_nsgmbeta}), we approximate the left-hand side of (\ref{s33_nsgmellip}) as
\begin{equation*}\label{s53_}
2(1-\sqrt{\rho})^{-3/4}\bigg[-\frac{1}{4}\log\Big(-B_1^*-2\sqrt{1-\sqrt{\rho}}\,\Big)+\log\Big(2\sqrt{2}(1-\sqrt{\rho})^{1/8}\Big)-1+o(1)\bigg],
\end{equation*}
which leads to the following approximation for $B_1^*(\sigma)$
\begin{equation}\label{s53_B1*match}
B_1^*(\sigma)\sim -2\sqrt{1-\sqrt{\rho}}-64\sqrt{1-\sqrt{\rho}}\,\exp\Big\{-4-2\sqrt{\rho}(1-\sqrt{\rho})^{3/4}\sigma\Big\},\quad \sigma\to \infty.
\end{equation}
Thus $B_1^*+2\sqrt{1-\sqrt{\rho}}$ is exponentially small for $\sigma\to\infty$.
Using (\ref{s53_B1*match}) in (\ref{s33_nsgmgamma}) and letting $\sigma\to\infty$, we have
\begin{equation}\label{s53_gamma*mch}
\gamma_*(\sigma)\sim 2^{-5}e^4(1-\sqrt{\rho})^{-5/4}\exp\Big\{2\sqrt{\rho}(1-\sqrt{\rho})^{3/4}\sigma\Big\},
\end{equation}
and we also find that
\begin{equation}\label{s53_intmch}
B_1^*\sqrt{\rho}\,\sigma -2\int_0^{\alpha_*}\sqrt{(1-\sqrt{\rho})v+1/v+B_1^*}\,dv=-2\sqrt{\rho}\sqrt{1-\sqrt{\rho}}\,\sigma-\frac{8}{3}(1-\sqrt{\rho})^{-1/4}+o(1).
\end{equation}
Using (\ref{s53_gamma*mch}) and (\ref{s53_intmch}) in (\ref{s33_nsigmapn}) and recalling that $\sigma=N^{1/4}\tau$, we compare the result to (\ref{s53_ntaupnt}), with $\tau\to 0$. If these two scales are to match in an intermediate region where $O(N^{3/4})\ll t\ll O(N)$, it follows that $\delta_4=-5/8$, 
\begin{equation}\label{s53_d1d2}
\quad d_1=0,\quad d_2=-\frac{8}{3}(1-\sqrt{\rho})^{-1/4},
\end{equation}
and
\begin{equation}\label{s53_nsgmmch}
\mathcal{R}_0=\frac{16\sqrt{\pi}}{(1-\sqrt{\rho})^{3/8}}\exp\Big(\frac{1}{1-\sqrt{\rho}}-2\Big).
\end{equation}
Thus, (\ref{s53_ntaupnt}) with (\ref{s53_ntauR}), (\ref{s53_d1d2}) and (\ref{s53_nsgmmch}) yields (\ref{s33_ntaupn}). This completes the analysis of the $(n,\tau)$ scale.

Now we consider the $(x,\tau)$ scale, where $n=\sqrt{N}\,x=O(\sqrt{N})$ and $t=N\tau=O(N)$. The $(n,\tau)$ scale we just analyzed and (\ref{s33_int}) suggests that now $p_n(t)$ has the following asymptotic expansion
\begin{eqnarray}\label{s53_xtaupnt}
p_n(t)&\sim& N^{\delta_5}\rho^{-n/2}\,e^{-(1-\sqrt{\rho})^2t}\exp\Big\{-2\sqrt{\rho}\sqrt{1-\sqrt{\rho}}\sqrt{N}\,\tau\Big\}\,e^{N^{1/4}\mathcal{F}(x,\tau)}\mathcal{G}(x,\tau)\nonumber\\
&&\times\exp\Big\{N^{1/4}\big(-\sqrt{\rho}(1-\sqrt{\rho})^{3/4}\,\tau-\frac{8}{3}(1-\sqrt{\rho})^{-1/4}\big)\Big\}.
\end{eqnarray}
Using (\ref{s53_xtaupnt}) in (\ref{s2_recu}) yields the PDEs
\begin{equation}\label{s53_xtauFPDE}
\mathcal{F}^2_x+(\sqrt{\rho}-1)x-\frac{1}{x}+2\sqrt{1-\sqrt{\rho}}=0
\end{equation}
and
\begin{equation}\label{s53_xtauGPDE}
\mathcal{G}_x+\frac{\mathcal{F}_{xx}+\big(1/x-x\big)\mathcal{F}_x-\mathcal{F}_{\tau}/\sqrt{\rho}+(1-\sqrt{\rho})^{3/4}}{2\mathcal{F}_x}\,\mathcal{G}=0.
\end{equation}
The solution to (\ref{s53_xtauFPDE}) involves an additive function of $\tau$. But this function must be zero if (\ref{s53_xtaupnt}) is to match to the $(\xi,\tau)$ result in $R_3$, hence
\begin{equation}\label{s53_xtauF}
\mathcal{F}(x,\tau)=\mathcal{F}(x)=2\sqrt{x}-\frac{2}{3}\sqrt{1-\sqrt{\rho}}\,x^{3/2}.
\end{equation}
Solving (\ref{s53_xtauGPDE}) with the help of (\ref{s53_xtauF}), we obtain
\begin{equation}\label{s53_xtauG}
\mathcal{G}(x,\tau)=\mathcal{G}_0(\tau)\frac{e^{x^2/4}e^{(1-\sqrt{\rho})^{1/4}\sqrt{x}}}{x^{1/4}\big[1+(1-\sqrt{\rho})^{1/4}\sqrt{x}\big]}.
\end{equation}
We determine $\mathcal{G}_0(\tau)$ through matching with the $(n,\tau)$ scale. In (\ref{s33_ntaupn}) we let $n=\sqrt{N}\,x\to\infty$ and use (\ref{s33_int}), and then we let $x\to 0$ in (\ref{s53_xtaupnt}). $\mathcal{G}(x,\tau)$ can be approximated for $x\to 0$ by $\mathcal{G}_0(\tau)x^{-1/4}$. Thus, if these two scales are to match in some intermediate limit where $O(1)\ll n\ll O(\sqrt{N})$ and $t=O(N)$, it follows that $\delta_5=-3/4$ and
\begin{equation}\label{s53_xtauG0}
\mathcal{G}_0(\tau)=\frac{8}{(1-\sqrt{\rho})^{3/8}}\exp\Big(\frac{1}{1-\sqrt{\rho}}-\frac{5}{2}\Big)\exp\bigg\{\frac{\sqrt{\rho}(15-22\sqrt{\rho}+3\rho)}{16(1-\sqrt{\rho})}\,\tau\bigg\}.
\end{equation}
Using (\ref{s53_xtauF})--(\ref{s53_xtauG0}) in (\ref{s53_xtaupnt}) yields (\ref{s33_xtaupn}).

Finally, we consider the matching between the $(\xi,\tau)$ (the region $R_3$) and $(x,\tau)$ scales. 
If we let $\xi\to 0$ in region $R_3$, (\ref{s33_R3B0}) leads to
\begin{equation}\label{s53_R3B0aym}
B_0\,e^{-\rho\tau}=(1-\sqrt{\rho})-\sqrt{\rho}\sqrt{1-\sqrt{\rho}}\sqrt{\xi}+O(\xi),\quad \xi\to 0,
\end{equation}
which is equivalent to
\begin{equation}\label{s53_B0xi0}
\frac{1}{\rho}\log\Big(\frac{B_0}{1-\sqrt{\rho}}\Big)=\tau+O(\sqrt{\xi}).
\end{equation}
Letting $\xi\to 0$ in (\ref{s33_R3psi}), (\ref{s33_R3psi(1)}), and (\ref{s33_R3psi(2)}) with the help of (\ref{s53_R3B0aym}), and using
\begin{equation*}\label{s53_psiasym}
\psi(\xi,\tau)=-(1-\sqrt{\rho})^2\tau-\frac{\xi}{2}\log\rho-\frac{2}{3}\sqrt{1-\sqrt{\rho}}\,\xi^{3/2}+\frac{\xi^2}{4}+O(\xi^{5/2}),
\end{equation*}
we see that for $\xi\to 0$ (\ref{s52_R3pn}) becomes
\begin{eqnarray}\label{s53_*mtch}
&&N^{\delta_1}\xi^{-3/4}L_0\Big(e^{\rho\tau}(1-\sqrt{\rho})\Big)\exp\left\{N\left[-(1-\sqrt{\rho})^2\tau-\frac{\xi}{2}\log\rho-\frac{2}{3}\sqrt{1-\sqrt{\rho}}\,\xi^{3/2}+\frac{\xi^2}{4}\right]\right\}\nonumber\\
&&\times\exp\left\{N^{1/2}\psi_1\Big(e^{\rho\tau}(1-\sqrt{\rho})\Big)-N^{1/2}\sqrt{\rho}\sqrt{1-\sqrt{\rho}}\sqrt{\xi}\,e^{\rho\tau}\psi_1'\Big(e^{\rho\tau}(1-\sqrt{\rho})\Big)\right\}\nonumber\\
&&\times\exp\left\{N^{1/4}\psi_2\Big(e^{\rho\tau}(1-\sqrt{\rho})\Big)-N^{1/4}\sqrt{\rho}\sqrt{1-\sqrt{\rho}}\sqrt{\xi}\,e^{\rho\tau}\psi_2'\Big(e^{\rho\tau}(1-\sqrt{\rho})\Big)\right\}\nonumber\\
&=&N^{\delta_1+3/8}\rho^{-n/2}x^{-3/4}L_0\Big(e^{\rho\tau}(1-\sqrt{\rho})\Big)e^{-(1-\sqrt{\rho})^2t}\exp\left\{N^{1/2}\psi_1\Big(e^{\rho\tau}(1-\sqrt{\rho})\Big)\right\}\nonumber\\
&&\times\exp\left\{N^{1/4}\left[-\frac{2}{3}\sqrt{1-\sqrt{\rho}}\,x^{3/2}-\sqrt{\rho}\sqrt{1-\sqrt{\rho}}\sqrt{x}\,e^{\rho\tau}\psi_1'\Big(e^{\rho\tau}(1-\sqrt{\rho})\Big)+\psi_2\Big(e^{\rho\tau}(1-\sqrt{\rho})\Big)\right]\right\}\nonumber\\
&&\times\exp\left\{-\sqrt{\rho}\sqrt{1-\sqrt{\rho}}\sqrt{x}\,e^{\rho\tau}\psi_2'\Big(e^{\rho\tau}(1-\sqrt{\rho})\Big)+\frac{x^2}{4}\right\}.
\end{eqnarray}
Comparing (\ref{s53_*mtch}) to large $x$ behavior of (\ref{s53_xtaupnt}) with the help of (\ref{s53_B0xi0}), we conclude that $\delta_1=\delta_5-3/8=-9/8$, 
\begin{equation*}\label{s53_R3psi1d1}
\psi_1(B_0)=-\frac{2\sqrt{1-\sqrt{\rho}}}{\sqrt{\rho}}\log\Big(\frac{B_0}{1-\sqrt{\rho}}\Big),
\end{equation*}
\begin{equation*}\label{s53_R3psi2d2}
\psi_2(B_0)=-\frac{(1-\sqrt{\rho})^{3/4}}{\sqrt{\rho}}\log\Big(\frac{B_0}{1-\sqrt{\rho}}\Big)-\frac{8}{3}(1-\sqrt{\rho})^{-1/4},
\end{equation*}
and
$$L_0(B_0)=(1-\sqrt{\rho})^{-1/4}\mathcal{G}_0\left(\frac{1}{\rho}\log\Big(\frac{B_0}{1-\sqrt{\rho}}\Big)\right).$$
Then we have derived (\ref{s33_R3pn}).
Note that $-\sqrt{\rho}\sqrt{1-\sqrt{\rho}}\sqrt{x}\,e^{\rho\tau}\psi_1'\Big(e^{\rho\tau}(1-\sqrt{\rho})\Big)=2\sqrt{x}$, which matches to the first term in (\ref{s53_xtauF}), while $-\sqrt{\rho}\sqrt{1-\sqrt{\rho}}\sqrt{x}\,e^{\rho\tau}\psi_2'\Big(e^{\rho\tau}(1-\sqrt{\rho})\Big)=(1-\sqrt{\rho})^{1/4}\sqrt{x}$, which matches to the second exponential factor in (\ref{s53_xtauG}).

Finally we note that the $(y,\tau)$ scale result in (\ref{s53_ypnt}) becomes a special case of the $(x,\tau)$ scale result, since $x=(1-\sqrt{\rho})^{-1/2}+yN^{-1/8}$ and for $y=O(1)$
$$\mathcal{F}(x)=\frac{4}{3}(1-\sqrt{\rho})^{-1/4}-N^{-1/4}\frac{(1-\sqrt{\rho})^{3/4}}{2}y^2+o(N^{-1/4}).$$
It follows that $\delta^*=\delta_5=-3/4$, $d_1^*=0$, $d_2^*=-\frac{4}{3}(1-\sqrt{\rho})^{-1/4}$, and
$$d_3^*=\frac{1}{2}\mathcal{G}_0(0)(1-\sqrt{\rho})^{1/8}\exp\left(\frac{1}{4(1-\sqrt{\rho})}+1\right).$$

This concludes the analysis of region $R_3$ and all the boundary layers. We have now completely determined the expansions in region $R_3$ and those on the $(n,\tau)$, $(y,\tau)$ and $(x,\tau)$ scales.

\subsection{Transition region $T_1$}

In this section we discuss the transition from $R_1$ to $R_2$, by considering the region $T_1$ defined in (\ref{s33_T1scale}). We scale $\xi$ as 
\begin{equation}\label{s54_theta}
\xi=\xi_0(\tau)+\frac{\theta}{\sqrt{N}}=\big(e^{\rho\tau}-1\big)\frac{1-\rho}{\rho}+\frac{\theta}{\sqrt{N}},\; \theta=O(1).
\end{equation} 
and then use it in the leading term approximation of $p_n(t)$ in $R_1$, which is given in (\ref{s32_xitau}) and (\ref{s32_Pxitau}) with $\rho<1$. It follows that, in the matching region between $R_1$ and $T_1$, where $\theta\to +\infty$ but $\xi-\xi_0(\tau)\to 0$, we should have 
\begin{equation}\label{s54_mchR1}
p_n(t)\sim N^{-\frac{2-\rho}{2(1-\rho)}}\,\frac{\rho\xi+1-\rho}{1-\rho}\,\xi^{-\frac{1}{1-\rho}}\,\theta^{\frac{\rho}{1-\rho}}.
\end{equation}
Next we use (\ref{s54_theta}) in the region $R_2$ result given by (\ref{s33_R2pn})--(\ref{s33_R2K0}), where we need $\theta<0$. Then (\ref{s33_R2B}) yields
\begin{equation}\label{s54_mchB}
B=\frac{1}{\sqrt{N}}\frac{\rho\xi+1-\rho}{\rho^2\xi+1-\rho^2}\frac{1-\rho}{\xi}(-\theta)+O\Big(\frac{1}{N}\Big).
\end{equation}
Using (\ref{s54_theta}) and (\ref{s54_mchB}) we have
\begin{equation*}\label{s54_mchphi}
N\phi(\xi,\tau)\sim-\frac{(1-\rho)^2\theta^2}{2\xi(\rho^2\xi+1-\rho^2)},
\end{equation*}
\begin{equation*}\label{s54_mchphi(1)}
\phi^{(1)}(\xi,\tau)\sim-\frac{3-\rho}{2(1-\rho)},
\end{equation*}
and
\begin{equation*}\label{s54_mchK}
K(\xi,\tau)\sim N^{\frac{1}{2(1-\rho)}}\frac{1}{\sqrt{2\pi}\,}(1-\rho)^{-\frac{2}{1-\rho}}\Gamma\Big(\frac{1}{1-\rho}\Big)\frac{\rho\xi+1-\rho}{\sqrt{\xi}}\,(\rho^2\xi+1-\rho^2)^{\frac{1+\rho}{2(1-\rho)}}\,(-\theta)^{-\frac{1}{1-\rho}}.
\end{equation*}
Thus, in the matching region between $R_2$ and $T_1$, where $\theta\to -\infty$ and $\xi-\xi_0(\tau)\to 0$, we should have
\begin{eqnarray}\label{s54_mchR2}
p_n(t)&\sim& N^{-\frac{2-\rho}{2(1-\rho)}}\frac{1}{\sqrt{2\pi}\,}(1-\rho)^{-\frac{2}{1-\rho}}\Gamma\Big(\frac{1}{1-\rho}\Big)\xi^{-\frac{2-\rho}{2(1-\rho)}}\big(\rho\xi+1-\rho\big)\,\big(\rho^2\xi+1-\rho^2\big)^{\frac{\rho}{2(1-\rho)}}\nonumber\\
&&\times\bigg(\frac{-\theta}{\sqrt{\xi}\sqrt{\rho^2\xi+1-\rho^2}}\bigg)^{-\frac{1}{1-\rho}}\exp\bigg\{-\frac{(1-\rho)^2}{2}\bigg(\frac{-\theta}{\sqrt{\xi}\sqrt{\rho^2\xi+1-\rho^2}}\bigg)^{2}\bigg\}.
\end{eqnarray}
The limits (\ref{s54_mchR1}) and (\ref{s54_mchR2}) in the matching regions suggest that in the transition region $T_1$ we expand $p_n(t)$ in the form 
\begin{equation}\label{s54_T1pnt}
p_n(t)\sim N^{-\frac{2-\rho}{2(1-\rho)}}V(\xi)H(\xi,\theta)
\end{equation}
where
\begin{equation}\label{s54_Vxi}
V(\xi)=\xi^{-\frac{2-\rho}{2(1-\rho)}}\big(\rho\xi+1-\rho\big)\,\big(\rho^2\xi+1-\rho^2\big)^{\frac{\rho}{2(1-\rho)}}.
\end{equation}
Using (\ref{s54_T1pnt}) in (\ref{s2_recu}) and noting that 
\begin{equation*}\label{s54_dt}
\frac{\partial }{\partial t}\Big(V(\xi)H(\xi,\theta)\Big)=\frac{\theta}{N} VH-\frac{\rho\xi+1-\rho}{\sqrt{N}}VH_\theta,
\end{equation*}
\begin{equation*}\label{s54_dxi}
\frac{\partial }{\partial \xi}\Big(V(\xi)H(\xi,\theta)\Big)=\sqrt{N}\,VH_\theta+\big(V'H+VH_\xi\big),
\end{equation*}
and
\begin{equation*}\label{s54_dxixi}
\frac{\partial ^2}{\partial \xi^2}\Big(V(\xi)H(\xi,\theta)\Big)=NVH_{\theta\theta}+\sqrt{N}\big(2VH_{\xi\theta}+V'H_\theta\big)+\big(V''H+2V'H_\xi+VH_{\xi\xi}\big),
\end{equation*}
we obtain the following limiting parabolic PDE for $H(\xi,\theta)$:
\begin{equation}\label{s54_PDEH}
\big(1+\rho-\rho\xi\big)H_{\theta\theta}-2\rho\theta H_\theta-2(\rho\xi+1-\rho)H_\xi=\frac{\rho(1-\rho)(1+\rho-\rho\xi)}{\xi(\rho^2\xi+1-\rho^2)}H.
\end{equation}
The matching results (\ref{s54_mchR1}) and (\ref{s54_mchR2}) also suggest that we seek a solution of (\ref{s54_PDEH}) in terms of the similarity variable $\Delta_1$, where 
\begin{equation}\label{s54_Delta1}
\Delta_1=\frac{\theta}{\sqrt{\xi}\sqrt{\rho^2\xi+1-\rho^2}},
\end{equation}
which is equivalent to (\ref{s33_Delta1}) after we solve (\ref{s54_theta}) for $\theta$ and using it in (\ref{s54_Delta1}). Thus, we let $H(\xi,\theta)=\mathcal{H}(\Delta_1)$, so that
$$H_\xi(\xi,\theta)=-\frac{\xi+\rho^2}{2\xi}\,\Delta_1\mathcal{H}'(\Delta_1),$$
and
$$H_\theta(\xi,\theta)=\frac{1}{\sqrt{\xi}\sqrt{\rho^2\xi+1-\rho^2}}\mathcal{H}'(\Delta_1),\quad H_{\theta\theta}(\xi,\theta)=\frac{1}{\xi(\rho^2\xi+1-\rho^2)}\mathcal{H}''(\Delta_1).$$
Then (\ref{s54_PDEH}) becomes a second order ODE for $\mathcal{H}(\Delta_1)$:
\begin{equation}\label{s54_ODEH}
\mathcal{H}''+(1-\rho)^2\Delta_1\mathcal{H}'-\rho(1-\rho)\mathcal{H}=0.
\end{equation}
If we rewrite (\ref{s54_mchR1}) and (\ref{s54_mchR2}) using (\ref{s54_Vxi}) and (\ref{s54_Delta1}), we have, in the respective matching regions, 
\begin{equation}\label{s54_mchR1Delta}
p_n(t)\sim N^{-\frac{2-\rho}{2(1-\rho)}}\frac{V(\xi)}{1-\rho}\,\Delta_1^{\frac{\rho}{1-\rho}},\quad \Delta_1\to +\infty
\end{equation}
and for $\Delta_1\to -\infty$
\begin{equation}\label{s54_mchR2Delta}
p_n(t)\sim N^{-\frac{2-\rho}{2(1-\rho)}}\frac{1}{\sqrt{2\pi}\,}(1-\rho)^{-\frac{2}{1-\rho}}\Gamma\Big(\frac{1}{1-\rho}\Big)V(\xi)\big(-\Delta_1\big)^{-\frac{1}{1-\rho}}\exp\bigg\{-\frac{(1-\rho)^2}{2}\Delta_1^{2}\bigg\}.
\end{equation}
Thus (\ref{s54_mchR1Delta}) and (\ref{s54_mchR2Delta}) provide two matching conditions for $\mathcal{H}(\Delta_1)$, which are
\begin{equation*}\label{s54_ODEcon1}
\mathcal{H}(\Delta_1)\sim \frac{1}{1-\rho}\,\Delta_1^{\frac{\rho}{1-\rho}},\; \Delta_1\to \infty
\end{equation*}
and
\begin{equation*}\label{s54_ODEcon2}
\mathcal{H}(\Delta_1)\sim \frac{1}{\sqrt{2\pi}\,}(1-\rho)^{-\frac{2}{1-\rho}}\Gamma\Big(\frac{1}{1-\rho}\Big)\big(-\Delta_1\big)^{-\frac{1}{1-\rho}}\exp\bigg\{-\frac{(1-\rho)^2}{2}\Delta_1^{2}\bigg\},\; \Delta_1\to -\infty.
\end{equation*}
The solution of (\ref{s54_ODEH}) that satisfies both matching conditions is given by  
\begin{equation}\label{s54_H(Del)}
\mathcal{H}(\Delta_1)=\frac{1}{\sqrt{2\pi}}\int_{-\Delta_1}^\infty\big(y+\Delta_1\big)^{\frac{\rho}{1-\rho}}\exp\bigg\{-\frac{(1-\rho)^2}{2}y^{2}\bigg\}\,dy.
\end{equation}
Using (\ref{s54_T1pnt}), (\ref{s54_Vxi}) and (\ref{s54_H(Del)}) leads to (\ref{s33_T1pn}). This completes the analysis of the first transition region.

\subsection{Transition region $T_2$}

Now we consider the transition from $R_2$ to $R_3$, with the second transition region $T_2$ defined in (\ref{s33_T2scale}). By defining $\Delta$ as in (\ref{s33_T2Delta}), the approximations in $R_3$ (cf. (\ref{s33_R3pn})) and on the $(x,\sigma)$ scale (cf. (\ref{s33_xsigmapn})) suggest that we assume $p_n(t)$ in the following form:
\begin{equation}\label{s55_pnt}
p_n(t)\sim N^{\delta_6}\rho^{-n/2}e^{-(1-\sqrt{\rho})^2t}e^{NF(\xi)}e^{N^{1/4}f(\xi,\Delta)}g(\xi,\Delta),
\end{equation}
where $F(\xi)$ is given in (\ref{s33_T2F}). 
Using (\ref{s55_pnt}) in (\ref{s2_recu}), the perturbation method leads to $f_\xi(\xi,\Delta)=0$, which indicates that $f$ is a function of $\Delta$ only, so we set $f(\xi,\Delta)=f(\Delta)$. The next-order in the perturbation expansion gives a PDE for $g(\xi,\Delta)$:
\begin{equation*}\label{s55_gPDE}
\Big[e^{-F'}-(1-\xi)e^{F'}\Big]g_\xi=\bigg\{\frac{1}{2}\Big[(1-\xi)e^{F'}+e^{-F'}\Big]\Big[F''+\big(f'(\Delta)\,\tau_*'(\xi)\big)^2\Big]-e^{F'}-\frac{e^{-F'}}{\xi}+\sqrt{\rho}\bigg\}g,
\end{equation*}
where $(\tau_*')^2=\big[\rho^2\xi^2+4\rho\xi(1-\sqrt{\rho})\big]^{-1}$, and $e^{-F'}-(1-\xi)e^{F'}=\sqrt{\rho\xi^2+4\xi(1-\sqrt{\rho})}$.
Thus, after some calculation and simplification, we obtain $g(\xi,\Delta)$ as
\begin{equation}\label{s55_g}
g(\xi,\Delta)=g_0(\Delta)J(\xi)\exp\bigg\{\frac{\big[4\rho(1-\sqrt{\rho})-(f'
(\Delta))^2\big]\big[2(1-\sqrt{\rho})+\sqrt{\rho}\,\xi\big]}{4\sqrt{\rho}(1-\sqrt{\rho})^2\sqrt{\rho^2\xi^2+4\rho\xi(1-\sqrt{\rho})}}\bigg\},
\end{equation}
where $J(\xi)$ is given in (\ref{s33_R3L}) and $g_0(\Delta)$ will be determined by matching.

We determine $f(\Delta)$ and $g_0(\Delta)$ by matching with the $(x,\sigma)$ scale. If we let $\xi\to 0$ in (\ref{s33_tau*}), $\tau_*(\xi)$ is asymptotically given by 
\begin{equation}\label{s55_tau*}
\tau_*(\xi)=\frac{\sqrt{\xi}}{\sqrt{\rho}\sqrt{1-\sqrt{\rho}}}+O(\xi^{3/2}).
\end{equation}
Using (\ref{s55_tau*}) in (\ref{s33_T2Delta}) and letting $\tau=\sigma N^{-1/4}$ and $\xi=xN^{-1/2}$ we have
\begin{equation}\label{s55_Delta}
\Delta=\Big(\sigma-\frac{\sqrt{x}}{\sqrt{\rho}\sqrt{1-\sqrt{\rho}}}\Big)+\frac{\sqrt{\rho}}{24(1-\sqrt{\rho})^{3/2}}x^{3/2}N^{-1/2}+O(N^{-1}).
\end{equation}
Thus, to match between the $(x,\sigma)$ scale and $T_2$, we need $\sigma,\,x\to\infty$ but keep 
\begin{equation*}\label{s55_Z}
Z\equiv\sigma-\frac{\sqrt{x}}{\sqrt{\rho}\sqrt{1-\sqrt{\rho}}}
\end{equation*}
fixed, so in the matching region we have $\Delta\sim Z$. Then (\ref{s55_Delta}) implies that, $f(\Delta)\sim f(Z)$, or more precisely, by a Taylor expansion, 
\begin{equation*}\label{s55_f}
f(\Delta)=f(Z)+f'(Z)\frac{\sqrt{\rho}}{24(1-\sqrt{\rho})^{3/2}}x^{3/2}N^{-1/2}+O(N^{-1}).
\end{equation*} 

We consider region $D_1$ on the $(x,\sigma)$ scale with $B_1(x,\sigma)$ given by (\ref{s33_xsigmaB1}). If we rewrite (\ref{s33_xsigmaB1}) as 
$$\left(\int_0^\infty-\int_x^\infty\right)\left[\frac{1}{\sqrt{(1-\sqrt{\rho})v+1/v+B_1}}-\frac{1}{\sqrt{(1-\sqrt{\rho})v}}\right]dv=2\sqrt{\rho}\sigma-\frac{2\sqrt{x}}{\sqrt{1-\sqrt{\rho}}}= 2\sqrt{\rho}Z$$
and approximate the second integral (integrating from $x$ to $\infty$) as $x\to\infty$, the above becomes
\begin{equation}\label{s55_intB1}
\int_0^\infty\left[\frac{1}{\sqrt{(1-\sqrt{\rho})v+1/v+B_1}}-\frac{1}{\sqrt{(1-\sqrt{\rho})v}}\right]dv=2\sqrt{\rho}Z-\frac{B_1}{(1-\sqrt{\rho})^{3/2}\sqrt{x}}+O(x^{-3/2}).
\end{equation}
From (\ref{s55_intB1}) it follows that, $B_1(x,\sigma)$ has, for $x,\,\sigma\to \infty$ with $Z=O(1)$, the expansion 
\begin{equation}\label{s55_B1}
B_1(x,\sigma)=A(Z)+\frac{1}{\sqrt{x}}A_1(Z)+O(x^{-1}).
\end{equation}
From (\ref{s55_intB1}), we find that the leading term $A(Z)$ will be defined implicitly by the integral
\begin{equation}\label{s55_A}
\int_0^\infty\bigg[\frac{\sqrt{v}}{\sqrt{(1-\sqrt{\rho})v^2+1+Av}}-\frac{1}{\sqrt{(1-\sqrt{\rho})v}}\bigg]dv=2\sqrt{\rho}Z,
\end{equation}
and $A_1(Z)$ is given in terms of $A(Z)$ by
\begin{equation}\label{s55_A1}
A_1=2A(1-\sqrt{\rho})^{-3/2}\left[\int_0^\infty\Big((1-\sqrt{\rho})v+1/v+A\Big)^{-3/2}dv\right]^{-1}.
\end{equation}
Then we expand $\eta(x,\sigma)$ in (\ref{s33_etaD1}) and $\gamma(x,\sigma)$ in (\ref{s33_gammaD1}) for $x\to\infty$, using (\ref{s55_B1}). After some calculation, including integration by parts, we obtain
\begin{eqnarray}\label{s55_eta}
\eta(x,\sigma)&=&-\frac{2}{3}\sqrt{1-\sqrt{\rho}}\,x^{3/2}+\frac{1}{3}\sqrt{\rho}AZ-\frac{2}{3}\int_0^\infty\frac{dv}{v\sqrt{(1-\sqrt{\rho})v+1/v+A}}\\
&&+\frac{2}{3\sqrt{x}}\left[\frac{\sqrt{\rho}}{2}A_1Z+\frac{A^2+12(1-\sqrt{\rho})}{8(1-\sqrt{\rho})^{3/2}}+\frac{A_1}{2}\int_0^\infty\frac{dv}{v\big[(1-\sqrt{\rho})v+1/v+A\big]^{3/2}}\right]+O(x^{-1})\nonumber
\end{eqnarray}
and 
\begin{equation*}\label{s55_gamma}
\gamma(x,\sigma)\sim\frac{\sqrt{2}}{(1-\sqrt{\rho})^{5/4}}\exp\left[\frac{1+\sqrt{\rho}}{2(1-\sqrt{\rho})}\right]e^{x^2/4}x^{-3/4}\left[\int_0^\infty\frac{dv}{\big[(1-\sqrt{\rho})v+1/v+A\big]^{3/2}}\right]^{-1/2}.
\end{equation*}
We note that the $O(x^{-1/2})$ term in (\ref{s55_eta}) can be further simplified using (\ref{s55_A}) and (\ref{s55_A1}) as
\begin{equation*}\label{s55_Oterm}
\left(\frac{1}{\sqrt{1-\sqrt{\rho}}}-\frac{A^2}{4(1-\sqrt{\rho})^{3/2}}\right)\frac{1}{\sqrt{x}}.
\end{equation*}
Thus, in the matching region between $T_2$ and the $(x,\sigma)$ scale, we should have the asymptotic approximation of $p_n(t)$ as
\begin{eqnarray}\label{s55_xsgmasym}
p_n(t)&\sim & N^{-3/4}\frac{\sqrt{2}}{(1-\sqrt{\rho})^{5/4}}\exp\left[\frac{1+\sqrt{\rho}}{2(1-\sqrt{\rho})}\right]\rho^{-n/2}e^{-(1-\sqrt{\rho})^2t}e^{x^2/4}x^{-3/4}\\
&&\times\exp\left\{\bigg(-\frac{2}{3}\sqrt{1-\sqrt{\rho}}\,x^{3/2}+\frac{1}{3}\sqrt{\rho}AZ-\frac{2}{3}\int_0^\infty\frac{dv}{v\sqrt{(1-\sqrt{\rho})v+1/v+A}}\bigg)N^{1/4}\right\}\nonumber\\
&&\times\exp\left\{\left(\frac{1}{\sqrt{1-\sqrt{\rho}}}-\frac{A^2}{4(1-\sqrt{\rho})^{3/2}}\right)\frac{N^{1/4}}{\sqrt{x}}\right\}\left[\int_0^\infty\frac{dv}{\big[(1-\sqrt{\rho})v+1/v+A\big]^{3/2}}\right]^{-1/2}.\nonumber
\end{eqnarray}

Now we let $\xi\to 0$ in (\ref{s55_g}), which leads to 
\begin{equation}\label{s55_gasym}
g(\xi,\Delta)\sim g_0(Z)\xi^{-3/4}\exp\left\{\left[\frac{1}{\sqrt{1-\sqrt{\rho}}}-\frac{\big(f'(Z)\big)^2}{4\rho(1-\sqrt{\rho})^{3/2}}\right]\frac{1}{\sqrt{\xi}}\right\}.
\end{equation}
The approximation of $F(\xi)$ as $\xi\to 0$ is given in (\ref{s33_Fxi0}). Then using (\ref{s33_Fxi0}) and (\ref{s55_gasym}) in (\ref{s55_pnt}), we have another approximation for $p_n(t)$ in the matching region,
\begin{eqnarray}\label{s55_T2asym}
p_n(t) &\sim& N^{\delta_6}\rho^{-n/2}e^{-(1-\sqrt{\rho})^2t}\exp\left\{N\left(-\frac{2}{3}\sqrt{1-\sqrt{\rho}}\,\xi^{3/2}+\frac{\xi^2}{4}\right)\right\} e^{N^{1/4}f(Z)}\nonumber\\
&&\times g_0(Z)\xi^{-3/4}\exp\left\{\left[\frac{1}{\sqrt{1-\sqrt{\rho}}}-\frac{\big(f'(Z)\big)^2}{4\rho(1-\sqrt{\rho})^{3/2}}\right]\frac{1}{\sqrt{\xi}}\right\},
\end{eqnarray}
and this must agree with (\ref{s55_xsgmasym}).
If we let $\xi=xN^{-1/2}$ in (\ref{s55_T2asym}) and compare it with (\ref{s55_xsgmasym}), it follows that $\delta_6=-9/8$, 
\begin{equation}\label{s55_f(D)}
f(Z)=\frac{1}{3}\sqrt{\rho}AZ-\frac{2}{3}\int_0^\infty\frac{1}{\sqrt{v}\sqrt{(1-\sqrt{\rho})v^2+1+Av}}dv
\end{equation}
and
\begin{equation}\label{s55_g0(D)}
g_0(Z)=\frac{\sqrt{2}}{(1-\sqrt{\rho})^{5/4}}\exp\Big[\frac{1+\sqrt{\rho}}{2(1-\sqrt{\rho})}\Big]\bigg[\int_0^\infty\big[(1-\sqrt{\rho})v+1/v+A\big]^{-3/2}dv\bigg]^{-1/2}.
\end{equation}
We note that $f'(Z)=\sqrt{\rho}A$, which can be easily verified from (\ref{s55_f(D)}) using the chain rule. Using the fact that $\Delta\sim Z$ in the matching region, (\ref{s55_A}), (\ref{s55_f(D)}) and (\ref{s55_g0(D)}) lead to (\ref{s33_T2A}), (\ref{s33_T2f(Delta)}) and (\ref{s33_T2g0}). Hence (\ref{s55_pnt}) is now completely determined, since $f(\xi,\Delta)=f(\Delta)$ and $g(\xi,\Delta)$ follows from (\ref{s55_g}), (\ref{s55_g0(D)}) and (\ref{s33_R3L}).

In the rest of this section, we verify the matching between the transition region $T_2$ and regions $R_2$ and $R_3$. Now we fix $\xi>0$ and consider the limits $\Delta\to \pm\infty$ in the $T_2$ result. To simplify the calculation, we set 
\begin{equation}\label{s55_A*Z*}
A_*=\frac{1}{\sqrt{1-\sqrt{\rho}}}A\quad \textrm{and}\quad Z_*=2\sqrt{\rho}(1-\sqrt{\rho})^{3/4}\Delta
\end{equation}
in (\ref{s33_T2A}). Then after an integration by parts, (\ref{s33_T2A}) becomes
\begin{equation}\label{s55_Z*int}
Z_*=-\int_0^\infty\frac{A_*+2/y}{\big[y+A_*+1/y\big]^{3/2}}dy.
\end{equation}

We first let $\Delta\to -\infty$, i.e., $Z_*\to -\infty$, in order to match with region $R_2$. Now we will have $A_*\to\infty$. By the substitution $y=A_*u$, ({\ref{s55_Z*int}) can be rewritten as
\begin{equation}\label{s55_Z*intR2}
\frac{|Z_*|}{\sqrt{A_*}}=\left(\int_0^\epsilon+\int_\epsilon^\infty\right)\frac{1+2/(uA_*^2)}{\big[u+1+1/(uA_*^2)\big]^{3/2}}du\equiv \textrm{I}+\textrm{II},
\end{equation}
where $0< \epsilon\ll 1$ and $\epsilon$ will be chosen appropriately below. In the first integral, we let $uA_*^2=v$ and assume that $\epsilon A_*^2\to \infty$, which yields
\begin{eqnarray}\label{s55_I}
\textrm{I}&=&\frac{1}{A^{2}_*}\int_{1/(\epsilon A_*^2)}^\infty\frac{1+2v}{v^2\big[v+1+1/(vA_*^2)\big]^{3/2}}dv\nonumber\\
&=& \frac{1}{A^{2}_*}\int_{1/(\epsilon A_*^2)}^\infty\frac{1+2v}{(1+v)^{3/2}v^2}dv+O(\epsilon^2)\nonumber\\
&\sim& \frac{1}{A^{2}_*}\Big[\epsilon A_*^2-\frac{3}{2}+\log 2+\frac{1}{2}\log(\epsilon A_*^2)\Big].
\end{eqnarray}
In the second integral, by expanding the integrand as $A_*\to\infty$, we have
\begin{eqnarray}\label{s55_II}
\textrm{II}&=&\int_\epsilon^\infty\frac{1}{(u+1)^{3/2}}\Big[1+\frac{4u+1}{2u(u+1)A_*^2}+O\Big(\frac{1}{u^2A_*^4}\Big)\Big]du\nonumber\\
&\sim& \frac{2}{\sqrt{1+\epsilon}}+\frac{1}{A^{2}_*}\Big[-\frac{1}{2}\log \epsilon+\log 2\Big].
\end{eqnarray}
Using (\ref{s55_I}) and (\ref{s55_II}) in (\ref{s55_Z*intR2}) with $A_*^{-2}\ll \epsilon\ll A_*^{-1}$, we have
\begin{equation*}\label{s55_Z*asym}
\frac{|Z_*|}{\sqrt{A_*}}\sim 2+\frac{\log A_*}{A_*^2}+\frac{2\log 2-3/2}{A_*^2},
\end{equation*}
which can be solved asymptotically for $A_*=A_*(Z_*)$ to give
\begin{equation}\label{s55_A*asym}
A_*\sim\frac{1}{4}Z_*^2+\frac{1}{Z_*^2}\big(6-4\log Z_*^2\big).
\end{equation}
Using (\ref{s55_A*Z*}) in (\ref{s55_A*asym}) yields (\ref{s33_T2Aasym-}). 

Now we expand $f(\Delta)$ in (\ref{s33_T2f(Delta)}) for $\Delta\to -\infty$. By setting $v=uA_*/\sqrt{1-\sqrt{\rho}}$, we have
\begin{equation}\label{s55_f(D)R2}
f(\Delta)=(1-\sqrt{\rho})^{-1/4}\left[\frac{1}{6}A_*Z_*-\frac{2}{3\sqrt{A_*}}\int_0^\infty\frac{du}{u\sqrt{u+1+1/(uA_*^2)}}\right].
\end{equation}
Similarly as in (\ref{s55_Z*intR2}), we split the integral in (\ref{s55_f(D)R2}) into two parts and approximate these as
$$\int_0^\epsilon\frac{du}{u\sqrt{u+1+1/(uA_*^2)}}\sim \log\epsilon+2\log A_*+2\log 2-\frac{\epsilon}{2}$$
and 
$$\int_\epsilon^\infty\frac{du}{u\sqrt{u+1+1/(uA_*^2)}}\sim -\log\epsilon+2\log 2.$$
Thus, the asymptotic approximation of $f(\Delta)$ is given by
\begin{equation}\label{s55_fR2asym}
f(\Delta)\sim (1-\sqrt{\rho})^{-1/4}\Big[\frac{Z_*^3}{24}+\frac{1}{Z_*}(1+2\log Z_*^2)\Big]+O(Z_*^{-5}).
\end{equation}
Using (\ref{s55_A*Z*}) in (\ref{s55_fR2asym}) yields (\ref{s33_T2fasym-}). (\ref{s33_T2g0asym-}) is obtained by first approximating the integral in (\ref{s33_T2g0}) by $2(1-\sqrt{\rho})^{-1}A^{-1/2}$ and then using (\ref{s33_T2Aasym-}).

Finally we let $\Delta\to \infty$ to verify the matching with region $R_3$. From (\ref{s55_Z*int}) we notice that if $A_*=-2$, the integrand has a double pole at $y=1$. We let $A_*=-2+2\delta$ with $0<\delta\ll 1$, and make the substitution $y=1-\delta+\sqrt{2\delta-\delta^2}\,\xi$; then (\ref{s55_Z*int}) can be rewritten as
\begin{equation}\label{s55_Z*R3}
Z_*=-2(1+\lambda^2)^{-1/4}\int_{-\lambda}^\infty\frac{\sqrt{\lambda+\xi}\,(1-\lambda\xi)}{(1+\xi^2)^{3/2}}\,d\xi,
\end{equation}
where
\begin{equation}\label{s55_lambda}
\lambda=\frac{1-\delta}{\sqrt{2\delta-\delta^2}}=\frac{1}{\sqrt{2\delta}}+O(\sqrt{\delta})\quad \textrm{as}\quad \delta\to 0.
\end{equation}
Then letting $\lambda\to\infty$ in (\ref{s55_Z*R3}), the integral can be approximated in the following two parts
$$\int_{0}^\infty\frac{\sqrt{\lambda+\xi}\,(1-\lambda\xi)}{(1+\xi^2)^{3/2}}\,d\xi\sim -\lambda^{3/2}-\frac{\sqrt{\lambda}}{2}\log\lambda+\sqrt{\lambda}\Big(1-\frac{3}{2}\log 2\Big)$$
and
$$\int_{-\lambda}^0\frac{\sqrt{\lambda+\xi}\,(1-\lambda\xi)}{(1+\xi^2)^{3/2}}\,d\xi\sim \lambda^{3/2}-\frac{\sqrt{\lambda}}{2}\log \lambda+\sqrt{\lambda}\Big(1-\frac{3}{2}\log 2\Big).$$
Thus, for $A_*\to -2$ (or $\lambda\to\infty$) $Z_*$ is asymptotically given by 
\begin{equation}\label{s55_Z*R3asym}
Z_*=2\log\lambda+(6\log 2-4)+o(1).
\end{equation}
Since $2\delta=A_*+2$, using (\ref{s55_lambda}) in (\ref{s55_Z*R3asym}) yields
$Z_*\sim -\log(A_*+2)+6\log 2-4$, and this in turn yields (\ref{s33_T2Aasym+}). (\ref{s33_T2fasym+}) and (\ref{s33_T2g0asym+}) follow by a similar asymptotic analysis, whose details we omit, noting only that the integral in (\ref{s33_T2f(Delta)}) can be approximated as
$$\int_0^\infty\frac{1}{\sqrt{v}\sqrt{(1-\sqrt{\rho})v^2+1+Av}}dv\sim (1-\sqrt{\rho})^{-1/4}\big(6\log 2+2\log\lambda\big),$$
and the integral in (\ref{s33_T2g0}) is asymptotically equal to
$$\int_0^\infty\big[(1-\sqrt{\rho})v+1/v+A\big]^{-3/2}dv\sim 2(1-\sqrt{\rho})^{-5/4}\lambda^2.$$
This completes the analysis of the transition region $T_2$.

\section{Asymptotic analysis for $p(t)$}

We now derive the asymptotic approximations for the unconditional density $p(t)$ in (\ref{s2_p(t)}). We first consider $\rho>1$ and then there is only one important time scale, $t=N\tau=O(N)$. We use (\ref{s32_xitau}) in (\ref{s2_p(t)}) and let $N\to\infty$. In view of (\ref{s2_pin}), the summation in (\ref{s2_p(t)}) is concentrated about $n=(1-1/\rho)N$. Letting $\xi^*=1-1/\rho$, we expand (\ref{s32_xitau}) about $\xi=\xi^*$ as follows
\begin{equation}\label{s6_pn}
p_n(t)=\frac{1}{N}P(\xi^*,\tau)+\frac{1}{N}(\xi-\xi^*)P_\xi(\xi^*,\tau)+\frac{1}{2N}(\xi-\xi^*)^2P_{\xi\xi}(\xi^*,\tau)+\frac{1}{N^2}P^{(1)}(\xi^*,\tau)+o(N^{-2}).
\end{equation}
Using the relation $\int_0^\infty v^n e^{-v}dv=n!$, we define
\begin{equation}\label{s6_D0}
\mathcal{D}_0\equiv\sum_{n=0}^{N-1}\frac{(N-1)!}{(N-n-1)!}\Big(\frac{\rho}{N}\Big)^n=\int_0^\infty\Big(1+\frac{\rho}{N}v\Big)^{N-1}e^{-v}dv.
\end{equation}
Then by the substitution $v=Nu$, we rewrite $\mathcal{D}_0$ as
\begin{equation*}\label{s6_D0u}
\mathcal{D}_0=N\int_0^\infty\frac{1}{1+\rho u}e^{N[-u+\log(1+\rho u)]}du.
\end{equation*}
Using the Laplace method, the above integrand is concentrated near the point $u=u_*=1-1/\rho$, where $\frac{d}{du}\big[-u+\log(1+\rho u)\big]=0$. Then we obtain a two-term asymptotic approximation of $\mathcal{D}_0$ as $N\to\infty$:
\begin{equation}\label{s6_D0asym}
\mathcal{D}_0=\sqrt{N}\sqrt{2\pi}\,e^{N(-1+1/\rho+\log\rho)}\bigg[\frac{1}{\rho}+\frac{1}{12\rho N}+O\Big(\frac{1}{N^2}\Big)\bigg].
\end{equation}
Now we differentiate both sides of (\ref{s6_D0}) with respect to $\rho$ twice, and after integration by parts and some simplification we obtain
\begin{equation}\label{s6_D1}
\mathcal{D}_1\equiv\sum_{n=0}^{N-1}\frac{(N-1)!}{(N-n-1)!}\Big(\frac{\rho}{N}\Big)^nn=\int_0^\infty(v-1)\Big(1+\frac{\rho}{N}v\Big)^{N-1}e^{-v}dv
\end{equation}
and
\begin{equation}\label{s6_D2}
\mathcal{D}_2\equiv\sum_{n=0}^{N-1}\frac{(N-1)!}{(N-n-1)!}\Big(\frac{\rho}{N}\Big)^nn^2=\int_0^\infty(v^2-3v+1)\Big(1+\frac{\rho}{N}v\Big)^{N-1}e^{-v}dv.
\end{equation}
Again by the Laplace method we have
\begin{equation}\label{s6_D1asym}
\int_0^\infty v\Big(1+\frac{\rho}{N}v\Big)^{N-1}e^{-v}dv=N^{3/2}\sqrt{2\pi}\,e^{N(-1+1/\rho+\log\rho)}\bigg[\frac{\rho-1}{\rho^2}+\frac{\rho-1}{12\rho^2N}+O\Big(\frac{1}{N^2}\Big)\bigg]
\end{equation}
and
\begin{equation}\label{s6_D2asym}
\int_0^\infty v^2\Big(1+\frac{\rho}{N}v\Big)^{N-1}e^{-v}dv=N^{5/2}\sqrt{2\pi}\,e^{N(-1+1/\rho+\log\rho)}\bigg[\frac{(\rho-1)^2}{\rho^3}+O\Big(\frac{1}{N}\Big)\bigg].
\end{equation}
Thus, using (\ref{s6_pn}) in (\ref{s2_p(t)}), we have
\begin{eqnarray*}\label{s6_p}
p(t)&=& \frac{1}{N}P(\xi^*,\tau)+\frac{1}{N^2}P_\xi(\xi^*,\tau)\frac{\mathcal{D}_1-N\xi^*\mathcal{D}_0}{\mathcal{D}_0}+\frac{1}{N^2}P^{(1)}(\xi^*,\tau)\nonumber\\
&&+\frac{1}{2N^3}P_{\xi\xi}(\xi^*,\tau)\frac{\mathcal{D}_2-2N\xi^*\mathcal{D}_1+(N\xi^*)^2\mathcal{D}_0}{\mathcal{D}_0}+o(N^{-2}),
\end{eqnarray*}
which yields, with the help of (\ref{s6_D0asym})--(\ref{s6_D2asym}),
\begin{equation*}\label{s6_pasym}
p(t)=\frac{1}{N}P(\xi^*,\tau)+\frac{1}{N^2}\Big[-P_\xi(\xi^*,\tau)+\frac{1}{2\rho}P_{\xi\xi}(\xi^*,\tau)+P^{(1)}(\xi^*,\tau)\Big]+O(N^{-3}).
\end{equation*}
Evaluating (\ref{s32_Pxitau}) and (\ref{s32_Pxitau1}) at $\xi=\xi^*$ and using 
$$P_\xi(\xi^*,\tau)=\frac{\rho^2}{(\rho-1)^3}\exp\Big(-\frac{\rho\tau}{\rho-1}\Big)-\frac{\rho^3}{(\rho-1)^3}\exp\Big(-\frac{\rho^2\tau}{\rho-1}\Big)$$
and
$$P_{\xi\xi}(\xi^*,\tau)=\frac{\rho^3}{(\rho-1)^5}\exp\Big(-\frac{\rho\tau}{\rho-1}\Big)\Big[2-\rho-2\rho e^{-\rho\tau}+\rho(2\rho-1) e^{-2\rho\tau}\Big],$$
we obtain (\ref{s31_rho>1_p}) and (\ref{s31_rho>1_p1}).

Next we consider $\rho<1$. Since $0<\rho<1$ and $N\to\infty$, $\pi_n$ in (\ref{s2_pin}) suggests that the summation in (\ref{s2_p(t)}) is concentrated for $n=O(1)$. We remove the condition on $n$ for the three time scales, $t=O(1)$, $t=O(N^{3/4})$ and $t=O(N)$. For times $t=O(1)$, we again find that we can approximate $p(t)$ by the infinite population model.

We now consider $t=O(N^{3/4})$ (the $\sigma$ scale) and use (\ref{s33_nsigmapn}) in (\ref{s2_p(t)}). Considering the factors involving $n$ in (\ref{s33_nsigmapn}) and taking the summation in (\ref{s2_p(t)}), we have
\begin{equation*}\label{s6_N/D0}
\sum_{n=0}^{N-1}\pi_n\rho^{-n/2}\frac{1}{2\pi i}\oint\frac{1}{z^{n+1}(1-z)}\exp\Big(\frac{1}{1-z}\Big)dz\equiv \frac{\mathcal{N}}{\mathcal{D}_0},
\end{equation*}
where $\mathcal{D}_0$ is as in (\ref{s6_D0}) and
\begin{equation}\label{s6_N}
\mathcal{N}=\frac{1}{2\pi i}\oint\frac{1}{z(1-z)}\exp\Big(\frac{1}{1-z}\Big)\sum_{n=0}^{N-1}\frac{(N-1)!}{(N-n-1)!}\Big(\frac{\rho}{zN}\Big)^ndz.
\end{equation}
Here the contour is over a small loop about $z=0$ with $|z|<\sqrt{\rho}$. Since $0<\rho<1$ and $n=O(1)$, we have $(N-1)!/(N-n-1)!\sim N^n$, which yields
\begin{equation*}\label{s6_<1D0asym}
\mathcal{D}_0\sim\frac{1}{1-\rho}\quad \textrm{as}\quad N\to\infty.
\end{equation*}
Rewriting the summation in (\ref{s6_N}) in integral form as
\begin{equation}\label{s6_sumint}
\sum_{n=0}^{N-1}\frac{(N-1)!}{(N-n-1)!}\Big(\frac{\rho}{zN}\Big)^n=\int_0^\infty\Big(1+\frac{\sqrt{\rho}}{zN}v\Big)^{N-1}e^{-v}dv,
\end{equation}
we notice that the integral in (\ref{s6_sumint}) is the same as in (\ref{s6_D0}) with $\rho$ replaced by $\sqrt{\rho}/z$. Thus, using the Laplace method, we obtain the leading term approximation
\begin{equation}\label{s6_sumasym}
\sum_{n=0}^{N-1}\frac{(N-1)!}{(N-n-1)!}\Big(\frac{\rho}{zN}\Big)^n\sim\sqrt{N}\frac{\sqrt{2\pi}}{\sqrt{\rho}}\,z\exp\Big\{N\Big[-1+\frac{z}{\sqrt{\rho}}+\log\frac{\sqrt{\rho}}{z}\Big]\Big\}.
\end{equation}
Using (\ref{s6_sumasym}) in (\ref{s6_N}), we have
\begin{equation}\label{s6_Nasym}
\mathcal{N}\sim \sqrt{N}\sqrt{2\pi}(\sqrt{\rho})^{N-1}e^{-N}\frac{1}{2\pi i}\oint\frac{1}{z^{N}(1-z)}\exp\Big(\frac{1}{1-z}+\frac{N}{\sqrt{\rho}}z\Big)dz.
\end{equation}
Rewriting the integral in (\ref{s6_Nasym}) as
\begin{equation*}\label{s6_intsadl}
\frac{1}{2\pi i}\oint\frac{1}{1-z}\exp\Big(\frac{1}{1-z}\Big)\exp\Big\{N\Big[\frac{z}{\sqrt{\rho}}-\log z\Big]\Big\}dz,
\end{equation*}
we use the saddle point method, with a saddle point at $z=\sqrt{\rho}$, which leads to the leading-term approximation of the above integral as
\begin{equation}\label{s6_intasym}
\frac{1}{\sqrt{N}\sqrt{2\pi}(1-\sqrt{\rho})}(\sqrt{\rho})^{1-N}e^N\exp\Big(\frac{1}{1-\sqrt{\rho}}\Big).
\end{equation}
Using (\ref{s6_intasym}) in (\ref{s6_Nasym}), we have
\begin{equation*}\label{s6_<1Nasym}
\mathcal{N}\sim\frac{1}{1-\sqrt{\rho}}\exp\Big(\frac{1}{1-\sqrt{\rho}}\Big),
\end{equation*}
which yields
\begin{equation}\label{s6_N/D0asym}
\frac{\mathcal{N}}{\mathcal{D}_0}\sim (1+\sqrt{\rho})\exp\Big(\frac{1}{1-\sqrt{\rho}}\Big).
\end{equation}
Thus, using (\ref{s33_nsigmapn}) in (\ref{s2_p(t)}) with (\ref{s6_N/D0asym}) leads to (\ref{s31_sigma}).

Finally we consider $t=N\tau=O(N)$. We notice that the factors containing $n$ in (\ref{s33_ntaupn}) are exactly the same as those in the $\sigma$-scale. Thus using (\ref{s33_ntaupn}) in (\ref{s2_p(t)}), (\ref{s6_N/D0asym}) still holds and therefore we obtain (\ref{s31_tau}). This completes the analysis of the unconditional density.

\section{Alternate approach to the asymptotics}

We use results in our previous paper \cite{ZHE_O12} on the eigenvalue problem (\ref{s2_recu_eigen}) to obtain the asymptotics of $p_n(t)$ from the spectral expansion in (\ref{s2_pn(t)eigen}), with (\ref{s2_cj}). However, as we show below, this type of analysis is only useful in the tail, where time is sufficiently large so that the asymptotic relation $p_n(t)\sim c_0\phi_0(n)e^{-\nu_0t}$ holds, or in somewhat shorter time ranges where all eigenvalues with $j=O(1)$ in (\ref{s2_recu_eigen}) contribute equally to the expansion of $p_n(t)$. The appropriate time ranges where these statements hold must be determined from the analysis, and we already identified these ranges from the singular perturbation analysis of sections 4 and 5. We shall mostly consider the case $\rho<1$, which is much more difficult than $\rho>1$. 

For $\rho<1$ and $N\to\infty$ a four-term approximation to the eigenvalues $\nu_j$ is given by (\ref{s2_nuj<1}). In \cite{ZHE_O12} we obtained the corresponding approximations to the eigenvectors $\phi_j(n)$. Then it was necessary to consider separately the four spatial scales: $n=O(N)$, $n=O(\sqrt{N})$, $n=\sqrt{N}/\sqrt{1-\sqrt{\rho}}+O(N^{3/8})$, and $n=O(1)$. On the third scale we showed that 
\begin{equation}\label{s7_y_phi}
\phi_j(n)\sim k_0\rho^{-n/2}e^{-z^2/4}\mathrm{He}_j(z),
\end{equation}
where
\begin{equation}\label{s7_y}
n=\frac{\sqrt{N}}{\sqrt{1-\sqrt{\rho}}}+\frac{N^{3/8}}{\sqrt{2}(1-\sqrt{\rho})^{3/8}}\,z,\quad z=O(1),
\end{equation}
and $\mathrm{He}_j(\cdot)$ is the $j^\mathrm{th}$ Hermite polynomial, so that $\mathrm{He}_0(z)=1$ and $\mathrm{He}_1(z)=z$. Also, $k_0$ is a normalizing constant, which could depend upon $j$. For $n=\sqrt{N}x=O(\sqrt{N})$, but $x\ne 1/\sqrt{1-\sqrt{\rho}}$ we showed that
\begin{equation}\label{s7_x_phi}
\phi_j(n)\sim k_1\rho^{-n/2}g_j(x)\exp\bigg\{N^{1/4}\Big[2\sqrt{x}-\frac{2}{3}\sqrt{1-\sqrt{\rho}}\,x^{3/2}\Big]\bigg\},
\end{equation}
where
\begin{equation}\label{s7_x_g}
g_j(x)=\frac{(1-\sqrt{1-\sqrt{\rho}}\,x)^j\,e^{x^2/4}}{x^{1/4}\big[1+(1-\sqrt{\rho})^{1/4}\sqrt{x}\big]^{2j+1}}\exp\Big[(2j+1)(1-\sqrt{\rho})^{1/4}\sqrt{x}\Big]
\end{equation}
and the constants $k_0$ and $k_1$ in (\ref{s7_y_phi}) and (\ref{s7_x_phi}) are related by
\begin{equation}\label{s7_k0k1}
k_0=k_1N^{-j/8}(-1)^j2^{-1-5j/2}(1-\sqrt{\rho})^{(j+1)/8}\exp\bigg[2j+1+\frac{1}{4(1-\sqrt{\rho})}+\frac{4}{3(1-\sqrt{\rho})^{1/4}}N^{1/4}\bigg].
\end{equation}
By expanding (\ref{s7_x_phi}) and (\ref{s7_x_g}) about $x=1/\sqrt{1-\sqrt{\rho}}$ we can easily show that (\ref{s7_y_phi}) is a special case of (\ref{s7_x_phi}) if $j=0$ and $j=1$, but this is not true for $j\ge 2$. In (\ref{s7_x_g}) $g_j(x)$ has a zero of order $j$ at $x=1/\sqrt{1-\sqrt{\rho}}$, while (\ref{s7_y_phi}) has simple zeros at the roots of the Hermite polynomials. 

On the scale $\xi=n/N\in (0,1)$ we found that 
\begin{equation}\label{s7_xi_phi}
\phi_j(n)\sim k_2\rho^{-n/2}G(\xi,j)\exp\Big\{NF(\xi)+\sqrt{N}F_1(\xi)+N^{1/4}F_2(\xi,j)\Big\}
\end{equation}
where $F(\xi)$ is as in (\ref{s33_T2F}), 
\begin{equation*}\label{s7_xi_F1}
F_1(\xi)=\frac{2\sqrt{1-\sqrt{\rho}}}{\sqrt{\rho}}\log\bigg[\frac{2(1-\sqrt{\rho})+\rho\xi+\sqrt{\rho^2\xi^2+4\rho\xi(1-\sqrt{\rho})}}{2(1-\sqrt{\rho})}\bigg],
\end{equation*}
\begin{equation*}\label{s7_xi_F2}
F_2(\xi,j)=(2j+1)\frac{(1-\sqrt{\rho})^{3/4}}{\sqrt{\rho}}\log\bigg[\frac{2(1-\sqrt{\rho})+\rho\xi+\sqrt{\rho^2\xi^2+4\rho\xi(1-\sqrt{\rho})}}{2(1-\sqrt{\rho})}\bigg],
\end{equation*}
and 
\begin{equation*}\label{s7_xi_G}
G(\xi,j)=\xi^{-3/4}\exp\bigg\{\int_0^\xi\Big[\frac{3}{4v}-H(v,j)\Big]dv\bigg\}
\end{equation*}
where
\begin{equation*}\label{s7_xi_H}
H(\xi,j)=\frac{1}{(1-\xi)e^{F'}-e^{-F'}}\bigg\{\frac{\widetilde{c}(j)+\rho}{\sqrt{\rho}}-e^{F'}-\frac{1}{\xi}e^{-F'}+\Big(1-\frac{1}{2}\sqrt{\rho}\,\xi\Big)\big[F''+(F_1')^2\big]\bigg\},
\end{equation*}
\begin{equation*}\label{s7_xi_c}
\widetilde{c}(j)=\frac{\sqrt{\rho}(22\sqrt{\rho}-3\rho-15)}{16(1-\sqrt{\rho})}-\frac{3}{8}\sqrt{\rho}(1-\sqrt{\rho})j(j+1).
\end{equation*}
We can show that 
\begin{equation*}\label{s7_G=J}
G(\xi,0)=J(\xi),
\end{equation*}
where $J(\cdot)$ is defined in (\ref{s33_R3L}). The constant $k_2$ in (\ref{s7_xi_phi}) is related to $k_1$ in (\ref{s7_x_phi}) by
\begin{equation*}\label{s7_k1k2}
k_2=(-1)^j(1-\sqrt{\rho})^{-1/4}N^{-3/8}k_1.
\end{equation*}
Finally, for $n=O(1)$ we have
\begin{equation}\label{s7_O(1)_phi}
\phi_j(n)\sim k_3\rho^{-n/2}\frac{1}{2\pi i}\oint\frac{1}{z^{n+1}(1-z)}\exp\Big(\frac{1}{1-z}\Big)dz,
\end{equation}
with
\begin{equation}\label{s7_k1k3}
k_3=k_1\frac{2\sqrt{\pi}}{\sqrt{e}}N^{1/8}.
\end{equation}

We use the above asymptotic results for $\phi_j(n)$ to evaluate the spectral coefficients $c_j$ in (\ref{s2_cj}). Let us normalize the eigenvectors by setting $k_3=1$, so that $k_1=N^{-1/8}\sqrt{e}/(2\sqrt{\pi})$. Consider the sum in the denominator in (\ref{s2_cj}). The summand concentrates on the scale where $z=O(1)$ (and thus $n\sim \sqrt{N}/\sqrt{1-\sqrt{\rho}}$) in (\ref{s7_y}). We use the asymptotic approximations
\begin{equation*}\label{s7_phiHe}
\rho^n\phi_j^2(n)\sim k_0^2e^{-z^2/2}\mathrm{He}_j^2(z)
\end{equation*}
and
\begin{equation*}\label{s7_Napro}
N^{-n}\frac{N!}{(N-n-1)!}\sim Ne^{-x^2/2}\sim N\exp\Big[-\frac{1}{2(1-\sqrt{\rho})}\Big].
\end{equation*}
Then, using the Euler-Maclaurin formula, we have
\begin{eqnarray}\label{s7_sum}
&&\sum_{n=0}^{N-1}\frac{N!}{(N-n-1)!}\Big(\frac{\rho}{N}\Big)^n(n+1)\phi_j^2(n)\nonumber\\
&\sim& k_0^2\,N^{15/8}\frac{1}{\sqrt{2}(1-\sqrt{\rho})^{7/8}}\exp\Big[-\frac{1}{2(1-\sqrt{\rho})}\Big]\int_{-\infty}^\infty e^{-z^2/2}\mathrm{He}_j^2(z)dz,
\end{eqnarray}
where
\begin{equation}\label{s7_k0^2}
k_0^2=\frac{e^{4j+3}}{2^{5j+4}\,\pi}N^{-j/4}\exp\bigg[\frac{1}{2(1-\sqrt{\rho})}+\frac{8}{3(1-\sqrt{\rho})^{1/4}}N^{1/4}\bigg],
\end{equation}
which follows from (\ref{s7_k0k1}), (\ref{s7_k1k3}) and the fact that we set $k_3=1$. The integral in (\ref{s7_sum}) is the well-know normalization constant for the Hermite polynomials, and its explicit value is
\begin{equation}\label{s7_normalHe}
\int_{-\infty}^\infty e^{-z^2/2}\mathrm{He}_j^2(z)dz=\sqrt{2\pi}\,j!.
\end{equation}

In contrast, the sum in the numerator in (\ref{s2_cj}) is concentrated in the range $n=O(1)$. Then we use the approximation in (\ref{s7_O(1)_phi}) for $\phi_j(n)$, and now note that $N^{-n}N!/(N-n-1)!\sim N$ so that
\begin{eqnarray}\label{s7_sumappx}
\sum_{n=0}^{N-1}\frac{N!}{(N-n-1)!}\Big(\frac{\rho}{N}\Big)^n\phi_j(n)&\sim &N\sum_{n=0}^\infty \rho^{n/2}\frac{1}{2\pi i}\oint\frac{1}{z^{n+1}(1-z)}\exp\Big(\frac{1}{1-z}\Big)dz\nonumber\\
&=&\frac{N}{2\pi i}\oint_{\sqrt{\rho}<|z|<1}\frac{1}{(1-z)(z-\sqrt{\rho})}\exp\Big(\frac{1}{1-z}\Big)dz\nonumber\\
&=&\frac{N}{1-\sqrt{\rho}}\exp\Big(\frac{1}{1-\sqrt{\rho}}\Big).
\end{eqnarray}
Using (\ref{s7_sum})--(\ref{s7_sumappx}) in (\ref{s2_cj}) we thus have 
\begin{equation}\label{s7_cjapprx}
c_j\sim\sqrt{\pi}\,\frac{2^{5j+4}}{j!}e^{-4j-3}N^{j/4-5/8}(1-\sqrt{\rho})^{-j/4-3/8}\exp\bigg[\frac{1}{1-\sqrt{\rho}}-\frac{8}{3(1-\sqrt{\rho})^{1/4}}N^{1/4}\bigg],\quad N\to\infty.
\end{equation}
This shows that, under the normalization $k_3=1$, the spectral coefficients are sub-exponentially small in $N$, with roughly $c_j=\exp\big[-O(N^{1/4})\big]$. The result in (\ref{s7_cjapprx}) applies for $j=O(1)$ only, and shows that $c_{j+1}/c_j=O(N^{1/4})$ for $N$ large, so that the asymptotic order of magnitude increases with the eigenvalue index $j$. 

In order for the result in (\ref{s7_cjapprx}) to be useful for evaluating $p_n(t)$ we must consider time ranges where the large eigenvalues (with $j\to\infty$) are not important asymptotically. This will also depend on the spatial range of $n$. For example, consider $n=O(\sqrt{N})$ with $n=\sqrt{N}\,x$, where (\ref{s7_x_phi}) applies. For times $t=N^{3/4}\sigma=O(N^{3/4})$ we furthermore scale $\sigma$ to be slightly large, with
\begin{equation*}\label{s7_sigma*}
\sigma=\frac{1}{8\sqrt{\rho}(1-\sqrt{\rho})^{3/4}}\log N+\sigma_*,\quad \sigma_*=O(1).
\end{equation*}
For fixed $\sigma_*$ we have
\begin{equation}\label{s7_sigma*appx}
e^{-\nu_jt}N^{j/4}\sim N^{-1/8}e^{-(1-\sqrt{\rho})^2t}\exp\Big[-2\sqrt{\rho}\sqrt{1-\sqrt{\rho}}\,N^{1/4}\sigma-(2j+1)\sqrt{\rho}(1-\sqrt{\rho})^{3/4}\sigma_*\Big].
\end{equation}
Using (\ref{s7_sigma*appx}) along with (\ref{s7_cjapprx}) we obtain from (\ref{s2_pn(t)eigen}) the leading term of the expansion of $p_n(t)$ on the $(x,\sigma_*)$ scale as
\begin{eqnarray}\label{s7_sigma*pn(t)}
p_n(t) &\sim & 8e^{-5/2}(1-\sqrt{\rho})^{-3/8}N^{-7/8}\exp\Big(\frac{1}{1-\sqrt{\rho}}\Big)\rho^{-n/2}e^{-(1-\sqrt{\rho})^2t}\nonumber\\
&&\times\exp\bigg[N^{1/4}\bigg(2\sqrt{x}-\frac{2}{3}\sqrt{1-\sqrt{\rho}}\,x^{3/2}-2\sqrt{\rho}\sqrt{1-\sqrt{\rho}}\,\sigma-\frac{8}{3(1-\sqrt{\rho})^{1/4}}\bigg)\bigg]\nonumber\\
&&\times\frac{x^{-1/4}e^{x^2/4}}{1+(1-\sqrt{\rho})^{1/4}\sqrt{x}}\exp\Big[(1-\sqrt{\rho})^{1/4}\sqrt{x}-\sqrt{\rho}(1-\sqrt{\rho})^{3/4}\sigma_*\Big]\nonumber\\
&&\times \exp\left\{\frac{32e^{-4}\big[1-(1-\sqrt{\rho})x\big]}{\big[1+(1-\sqrt{\rho})^{1/4}\sqrt{x}\big]^2}\,e^{2(1-\sqrt{\rho})^{1/4}\big(\sqrt{x}-\sqrt{\rho}\sqrt{1-\sqrt{\rho}}\,\sigma_*\big)}\right\}.
\end{eqnarray}
But the expansion in (\ref{s7_sigma*pn(t)}) can also be obtained by expanding (\ref{s33_xsigmapn}) for $\sigma\to\infty$ and $x=O(1)$. For times large enough so that $\sigma_*\to\infty$, the last (double-exponential) factor in (\ref{s7_sigma*pn(t)}) can be approximated by unity, and then (\ref{s7_sigma*pn(t)}) agrees with the small $\tau$ expansion of the result in (\ref{s33_xtaupn}), which applies on the $(x,\tau)$ scale. We have thus shown that the spectral results in \cite{ZHE_O12} can also be used to obtain the expansion of $p_n(t)$, but only for sufficiently large time ranges. For the scale $x=O(1)$, the expression in (\ref{s7_sigma*pn(t)}) applies in the asymptotic matching region between the $\sigma$ ($t=O(N^{3/4})$) and $\tau$ ($t=O(N)$) time scales. 

By considering $n=O(1)$ and using (\ref{s7_O(1)_phi}) to approximate the eigenvectors, we can obtain an analogous approximation to $p_n(t)$ on the $(n,\sigma_*)$ scale. We omit that analysis, but note that from this result we can immediately remove the condition on $n$ to get (\ref{s31_sigma*}). Also, by considering $n=N\xi=O(N)$ and using (\ref{s7_xi_phi}) we can easily recover the approximation in (\ref{s33_T2R3}), on the $(\xi,\widetilde{\Delta})$ scale. Thus it is possible to obtain the results for $p_n(t)$ in Region 3, where the zeroth eigenvalue $\nu_0$ dominates, and also in the matching regions between the $(x,\sigma)$ and $(x,\tau)$ scales, and between the $T_2$ and $R_3$ ranges, where all the eigenvalues $\nu_j$ with $j=O(1)$ contribute equally to the expansion of $p_n(t)$. However, it is not feasible to obtain, for example, the $R_1$, $R_2$ and $T_1$ results using spectral asymptotics. At the very least this would require analyzing the eigenvalues $\nu_j$ for $N$ and $j$ simultaneously large. 

We briefly discuss the case $\rho>1$. Now the $\nu_j$ are all $O(N^{-1})$, with (\ref{s2_nuj>1}) holding, and in \cite{ZHE_O12} we showed that $\phi_j(n)$ is concentrated in the range $n=N(1-1/\rho)+O(\sqrt{N})$, with
\begin{equation}\label{s7_rho>1_phi}
\phi_j(n)\sim k_0\,\mathrm{He}_j(\sqrt{\rho}\,X),\quad X=\frac{n-N(1-1/\rho)}{\sqrt{N}}=O(1).
\end{equation}
On the scale $\xi=O(1)$, $\xi\ne 1-1/\rho$, a different result applies, with
\begin{equation*}\label{s7_rho>1_phixi}
\phi_j(n)\sim k_0\,N^{j/2}\rho^{j/2}\Big(1-\frac{1}{\rho}\Big)^{1/(1-\rho)}\Big(\xi-1+\frac{1}{\rho}\Big)^{j}\xi^{1/(\rho-1)}.
\end{equation*}
For $X=O(1)$ in (\ref{s7_rho>1_phi}), an application of Stirling's formula shows that
\begin{equation}\label{s7_rho>1_stirling}
\frac{N!}{(N-n-1)!}\Big(\frac{\rho}{N}\Big)^n\sim\frac{N}{\sqrt{\rho}}\,e^{N(\log \rho-1+1/\rho)}e^{-\rho X^2/2}.
\end{equation}
Then using (\ref{s7_rho>1_stirling}) along with (\ref{s7_rho>1_phi}) we asymptotically evaluate the sum in the denominator in (\ref{s2_cj}) as
\begin{eqnarray}\label{s7_rho>1_cjdeno}
\sum_{n=0}^{N-1}\frac{N!}{(N-n-1)!}\Big(\frac{\rho}{N}\Big)^n(n+1)\phi_j^2(n)&\sim &k_0^2\frac{N^{5/2}}{\sqrt{\rho}}\Big(1-\frac{1}{\rho}\Big)e^{N(\log\rho-1+1/\rho)}\int_{-\infty}^\infty e^{-\rho X^2/2}\mathrm{He}_j^2(\sqrt{\rho}\,X) dX\nonumber\\
&=&k_0^2\,N^{5/2}\frac{\sqrt{2\pi}}{\rho}\Big(1-\frac{1}{\rho}\Big)j!\,e^{N(\log\rho-1+1/\rho)}.
\end{eqnarray}
An analogous calculation shows that the sum in the numerator in (\ref{s2_cj}) can be approximated by
\begin{equation}\label{s7_rho>1_cjnume}
\sum_{n=0}^{N-1}\frac{N!}{(N-n-1)!}\Big(\frac{\rho}{N}\Big)^n\phi_j(n)\sim k_0\,\frac{N^{3/2}}{\sqrt{\rho}}\,e^{N(\log\rho-1+1/\rho)}\int_{-\infty}^\infty e^{-\rho X^2/2}\mathrm{He}_j(\sqrt{\rho}\,X)dX.
\end{equation}
But, in view of the orthogonality of the Hermite polynomials the integral in (\ref{s7_rho>1_cjnume}) vanishes for all $j\ge 1$. Thus (\ref{s7_rho>1_cjnume}) yields the leading term only when $j=0$, and shows that the sum in (\ref{s7_rho>1_cjnume}) is $o\big(N^{3/2}\exp[N(\log\rho-1+1/\rho)]\big)$ for $j\ge 1$, in which case obtaining a more precise estimate would necessitate the calculation of some of the higher order terms in the asymptotic relation in (\ref{s7_rho>1_phi}). From (\ref{s2_cj}), (\ref{s7_rho>1_cjdeno}) and (\ref{s7_rho>1_cjnume}) we thus conclude that 
\begin{equation}\label{s7_rho>1_c0}
c_0\sim\frac{1}{N}\frac{\rho}{\rho-1}\frac{1}{k_0},\quad c_j=o(N^{-1}),\quad j\ge 1.
\end{equation}

Now consider $p_n(t)$ in (\ref{s2_pn(t)eigen}). On the scale $0<\xi<1$, $\xi\ne 1-1/\rho$, we cannot estimate the conditional density, as we do not have the $c_j$ for $j\ge 1$. But, using the perturbation method we obtained the approximation in (\ref{s32_xitau}), with (\ref{s32_Pxitau}). By expanding (\ref{s32_Pxitau}) as a binomial series involving powers of $e^{-\rho\tau}$, and noting that, in view of (\ref{s2_nuj>1}),
$$e^{-\nu_jt}\sim\exp\Big(-\frac{\rho}{\rho-1}\,\tau\Big)e^{-\rho j\tau},$$
we can infer that $c_j=O(N^{-j/2-1})$, and indeed also get the constant in this relation, since
\begin{eqnarray*}\label{s7_rho>1_cjappx}
&&\xi^{1/(\rho-1)}\exp\Big(-\frac{\rho\tau}{\rho-1}\Big)\Big(\frac{\rho-1}{\rho}\Big)^{\rho/(1-\rho)}\bigg[1+\Big(\frac{\rho}{\rho-1}\xi-1\Big)e^{-\rho\tau}\bigg]^{-\rho/(\rho-1)}\nonumber\\
&=& \xi^{1/(\rho-1)}\Big(1-\frac{1}{\rho}\Big)^{1/(1-\rho)}\sum_{j=0}^\infty\widetilde{c}_j\rho^{j/2}\Big(\xi-1+\frac{1}{\rho}\Big)^j\exp\Big[-\Big(\frac{\rho}{\rho-1}+\rho j\Big)\tau\Big].
\end{eqnarray*}
Here $\widetilde{c}_j\equiv \lim_{N\to\infty}(c_jN^{1+j/2})$ and we can set $k_0=1$ for each $j$. Thus in this case the perturbation method is a much more efficient approach for getting the asymptotics of $p_n(t)$. However, on the $X$-scale, we have $\phi_j(n)=O(1)$ for each $j$ and then (\ref{s7_rho>1_phi}) and (\ref{s7_rho>1_c0}) shows that
\begin{equation}\label{s7_rho>1_pn(t)}
p_n(t)\sim c_0\phi_0(n)e^{-\nu_0t}\sim\frac{1}{N}\frac{\rho}{\rho-1}\exp\Big(-\frac{\rho}{\rho-1}\tau\Big),\quad X=O(1).
\end{equation}
This is the conditional result, but the unconditional result is identical (and agrees with the leading term in (\ref{s31_rho>1_p})) since removing the condition on $n$ would asymptotically be the same as applying $\int_{-\infty}^\infty e^{-\rho X^2/2}\sqrt{\rho/(2\pi)}\,(\cdot)dX$ to the right side of (\ref{s7_rho>1_pn(t)}), which amounts to doing nothing as the expression is independent of $X$.

Our discussion also shows that when $\rho>1$, the zeroth eigenvalue $\nu_0$ is dominant even for times $t=O(1)$, as long as $N$ is large. Calculating the correction term $\mathcal{P}_1(\tau)$ in (\ref{s31_rho>1_p}) would require that we precisely estimate $c_1$ and also get the $O(N^{-2})$ correction term in (\ref{s2_nuj>1}) for $\nu_0$.

\section{Large $\rho$, fixed $N$ limit}

We examine the eigenvalue problem for (\ref{s2_recuMatrx}) (with (\ref{s2_A})), for $\rho\to\infty$ with $N=O(1)$. We shall thus verify the conjectures we made in section 2, after we considered the cases $N=2,\,3$. Consider (\ref{s2_recu_eigen}) and expand the eigenvalues and eigenvectors as
\begin{equation}\label{s8_eigen}
\begin{array}{l}
\displaystyle \nu_j=\rho\nu_j^{(0)}+\nu_j^{(1)}+\rho^{-1}\nu_j^{(2)}+O(\rho^{-2}),\\
\\
\displaystyle \phi_j(n)=\phi_j^{(0)}(n)+\rho^{-1}\phi_j^{(1)}(n)+O(\rho^{-2}).
\end{array} 
\end{equation}
Using (\ref{s8_eigen}) in (\ref{s2_recu_eigen}) leads to the difference equations, for $0\le n\le N-1$,
\begin{equation}\label{s8_phi0}
-\nu_j^{(0)}\phi_j^{(0)}(n)=\Big(1-\frac{n+1}{N}\Big)\Big[\phi_j^{(0)}(n+1)-\phi_j^{(0)}(n)\Big]
\end{equation}
\begin{equation}\label{s8_phi1}
-\nu_j^{(0)}\phi_j^{(1)}(n)-\nu_j^{(1)}\phi_j^{(0)}(n)=\Big(1-\frac{n+1}{N}\Big)\Big[\phi_j^{(1)}(n+1)-\phi_j^{(1)}(n)\Big]+\frac{n}{n+1}\phi_j^{(0)}(n-1)-\phi_j^{(0)}(n).
\end{equation}
We also note that the matrix $\mathbf{A}$ in (\ref{s2_A}) may be written as $\mathbf{A}=\rho \mathbf{A_0}+\mathbf{A_1}$, where
\begin{equation}\label{s8_A0}
\mathbf{A_0}=\left[\begin{array}{cccccc}
-\frac{N-1}{N} & \frac{N-1}{N} &  0 & \cdots & 0 & 0\\
0 & -\frac{N-2}{N} & \frac{N-2}{N}  & \cdots & 0 & 0 \\
\vdots  & \vdots  & \ddots  & \ddots  & \vdots  & \vdots \\
0 & 0 & 0  & \cdots & -\frac{1}{N}  & 0\\
0 & 0 & 0  & \cdots & 0  & 0\\
\end{array}\right]
\end{equation}
and
\begin{equation}\label{s8_A1}
\mathbf{A_1}=\left[\begin{array}{cccccc}
-1 & 0 &  0 & \cdots & 0 & 0\\
1/2 & -1 & 0  & \cdots & 0 & 0 \\
0 & 2/3 & -1  & \cdots & 0 & 0 \\
\vdots  & \vdots  & \ddots  & \ddots  & \vdots  & \vdots \\
0 & 0 & 0  & \cdots & -1  & 0\\
0 & 0 & 0  & \cdots & \frac{N-1}{N}  & -1\\
\end{array}\right],
\end{equation}
so that $\mathbf{A_0}$ (resp. $\mathbf{A_1}$) is an upper (resp. lower) bi-diagonal matrix. Using (\ref{s8_A0}) and (\ref{s8_A1}) we can also recast the expansion in (\ref{s8_eigen}) as an expansion of a perturbed eigenvalue problem, with the leading term satisfying $\mathbf{A_0}\overrightarrow{\phi}_j^{(0)}=-\nu_j^{(0)}\overrightarrow{\phi}_j^{(0)}$, where $\overrightarrow{\phi}_j^{(0)}$ is the column vector $(\phi_j^{(0)}(0),\,\phi_j^{(0)}(1),\cdots,\,\phi_j^{(0)}(N-1))^T$. The correction term corresponds to $\mathbf{A_0}\overrightarrow{\phi}_j^{(1)}+\mathbf{A_1}\overrightarrow{\phi}_j^{(0)}
=-\nu_j^{(0)}\overrightarrow{\phi}_j^{(1)}-\nu_j^{(1)}\overrightarrow{\phi}_j^{(0)}$.

From (\ref{s8_phi0}) we see that $\nu_j^{(0)}=0$ is a possible solution (which corresponds to the zeroth eigenvalue), with $\phi_0^{(0)}(n)=$ constant, and we set $\phi_0^{(0)}(n)=1$. Then rearranging (\ref{s8_phi1}) leads to
\begin{equation}\label{s8_rearrang83}
-\nu_0^{(1)}+\frac{1}{n+1}=\Big(1-\frac{n+1}{N}\Big)\Big[\phi_0^{(1)}(n+1)-\phi_0^{(1)}(n)\Big],\; 0\le n\le N-1.
\end{equation}
The solvability condition for (\ref{s8_rearrang83}) is obtained by setting $n=N-1$ from which we conclude that $\nu_0^{(1)}=1/N$. Then (\ref{s8_rearrang83}) becomes $\phi_0^{(1)}(n+1)-\phi_0^{(1)}(n)=1/(n+1)$ so that, up to an additive constant, $\phi_0^{(1)}(n)=H_n=1+1/2+\cdots+1/n$ is the harmonic number. We have thus found that for $j=0$
\begin{equation*}\label{s8_j=0}
\nu_0\sim \frac{1}{N},\quad \phi_0(n)=1+\frac{1}{\rho}H_n+O(\rho^{-2}).
\end{equation*}
Here $\phi_0(n)$ is the zeroth unnormalized eigenvector. 

To obtain the remaining eigenvalues we could simply examine the diagonal entries in (\ref{s8_A0}), or we can introduce the generating function
\begin{equation}\label{s8_generating}
\mathcal{F}_j(z)=\sum_{n=0}^{N-1}\phi_j^{(0)}(n)z^n.
\end{equation}
Using (\ref{s8_generating}) in (\ref{s8_phi0}) leads, after some calculation, to the ODE
\begin{equation}\label{s8_ODE}
\mathcal{F}'_j(z)\Big(\frac{z-1}{N}\Big)+\mathcal{F}_j(z)\Big(\frac{1}{N}-1+\frac{1}{z}+\nu_j^{(0)}\Big)=\frac{1}{z}\phi_j^{(0)}(0).
\end{equation}
Now, $\mathcal{F}_j(z)$ is a polynomial and hence an entire function of $z$. The only solution of (\ref{s8_ODE}) that is analytic at $z=1$ is given by
\begin{equation}\label{s8_ODEsln}
\mathcal{F}_j(z)=\phi_j^{(0)}(0)\frac{N z^N}{(z-1)^{1+N\nu_j^{(0)}}}\int_1^z\frac{(u-1)^{N\nu_j^{(0)}}}{u^{N+1}}\,du.
\end{equation}

We can clearly normalize the $j^\mathrm{th}$ eigenvector by setting its first component $\phi_j^{(0)}(0)=1$. Now, (\ref{s8_ODEsln}) must also be analytic at $z=0$. Note that as $z\to 0$, $\mathcal{F}_j(z)\to \phi_j^{(0)}(0)=1$. By expanding $(u-1)^{N\nu_j^{(0)}}$ as a binomial series we see that if the coefficient of $u^N$ in this series is non-zero, the expansion of (\ref{s8_ODEsln}) about $z=0$ will contain a term proportional to $z^N\log z$. To avoid this lack of analyticity we must have $N\nu_j^{(0)}$ to be a positive integer that is $<N$. We thus have $\nu_j^{(0)}=j/N$, for $1\le j\le N-1$ ($j=0$ will regain the zeroth eigenvalue).

To identify the $\phi_j^{(0)}(n)$ we expand the right side of (\ref{s8_ODEsln}) as
\begin{eqnarray*}\label{s8_Fexpan}
\mathcal{F}_j(z)&=&\frac{Nz^N}{(z-1)^{j+1}}\int_1^z\frac{(u-1)^j}{u^{N+1}}du\nonumber\\
&=&\frac{Nz^N}{(z-1)^{j+1}}\sum_{l=0}^j\int_1^z{{j}\choose{l}}(-1)^{l-j}u^{l-N-1}du\nonumber\\
&=&\frac{N}{(1-z)^{j+1}}\sum_{l=0}^{j}(-1)^l{{j}\choose{l}}\frac{z^l-z^N}{N-l}.
\end{eqnarray*}
Then using a binomial expansion for $(1-z)^{-j-1}$ and letting $\big[z^n\big]\mathcal{F}_j(z)$ denote the coefficient of $z^n$ in the Taylor expansion of $\mathcal{F}_j(z)$, we obtain
\begin{equation}\label{s8_[zn]}
\phi_j^{(0)}(n)=\big[z^n\big]\mathcal{F}_j(z)=\sum_{\rho=\max\{n-j,0\}}^n\frac{N(j+\rho)!(-1)^{n-\rho}}{(N+\rho-n)\rho!(n-\rho)!(j-n+\rho)!}.
\end{equation}
An alternate form of the eigenvector may be obtained by simply setting $\nu_j^{(0)}=j/N$ in (\ref{s8_phi0}) and solving the elementary difference equation, which yields
\begin{equation}\label{s8_phi0sln}
\phi_j^{(0)}(n)=
\left\{ \begin{array}{ll}
\displaystyle \frac{(N-j-1)!(N-n-1)!}{(N-n-j-1)!(N-1)!}, & 0\le n\le N-j-1\\
\\
0, & N-j\le n\le N-1.
\end{array} \right.
\end{equation}
Hence, only the first $N-j$ components, of the leading term in the large $\rho$ approximation to the eigenvectors, are non-zero. To summarize, we have shown that for $1\le j\le N-1$, $\nu_j\sim \rho j/N$ and $\phi_j(n)\sim \phi_j^{(0)}(n)$, as given by (\ref{s8_[zn]}) or (\ref{s8_phi0sln}). This verifies our conjectures in section 2.

It is interesting to note that $\phi_j^{(0)}(n)$ are all of one sign, and this would seem to be inconsistent with the orthogonality relation in (\ref{s2_orth}). But in view of (\ref{s2_pin}) the orthogonality relation gives the largest weight, for $\rho\to \infty$, to the $n=N-1$ components of the eigenvectors. Indeed (\ref{s2_orth}) may be asymptotically rewritten as, for $j\ne k$,
\begin{equation}\label{s8_orth}
\rho^{N-1}N\phi_j(N-1)\phi_k(N-1)+N(N-1)\rho^{N-2}\phi_j(N-2)\phi_k(N-2)+O(\rho^{N-3})=0,
\end{equation}
where we assume all eigenvectors are normalized to be $O(1)$ as $\rho\to\infty$. But if $j,\,k\ge 1$ we have $\phi_j^{(0)}(N-1)=0$ so that the higher order terms become important in evaluating (\ref{s8_orth}). Even for $j=0$ and $k=1$ we have $\phi_0(N-1)\sim \rho^{-1}\phi_0^{(1)}(N-1)$, $\phi_1(N-2)\sim 1/(N-1)$ and (\ref{s8_orth}) becomes
\begin{equation*}\label{s8_orth2}
\rho^{N-2}\big[N\phi_1^{(1)}(N-1)+N+O(\rho^{-1})\big]=0,
\end{equation*}
which can be used to infer the value of $\phi_1^{(1)}(N-1)$, which we did not previously calculate. We can also conclude that $\phi_1^{(1)}(N-1)=-1$ by setting $j=1$ and $n=N-1$ in (\ref{s8_phi1}), and using the fact that $\phi_1^{(0)}(n)=(N-n-1)/(N-1)$.

\appendix
\section{Appendix}

We briefly discuss using generating functions to solve (\ref{s2_recu}), with (\ref{s2_ic}). Using an ordinary generating function
\begin{equation*}\label{A1}
\mathcal{G}(t,z)=\sum_{n=0}^{N-1} z^np_n(t)
\end{equation*}
in (\ref{s2_recu}) leads to, after some calculation, 
\begin{equation}\label{A2}
\mathcal{G}_t+z\mathcal{G}_{zt}=\frac{\rho}{N}(z^2-z)\mathcal{G}_{zz}+\Big[\frac{\rho}{N}(2z-1)+(1-z)(\rho-z)\Big]\mathcal{G}_z+(z-\rho-1)\mathcal{G}-Nz^Np_{_{N-1}}(t),
\end{equation}
and the initial condition in (\ref{s2_ic}) yields
\begin{equation}\label{A3}
\mathcal{G}(0,z)=\sum_{n=0}^{N-1} \frac{z^n}{n+1}.
\end{equation}
But (\ref{A2}) is a complicated PDE which is also a functional equation, due to the last term in the right-hand side. We have not been able to analyze this equation.

\newpage

\begin{center}
\includegraphics[width=0.6\textwidth]{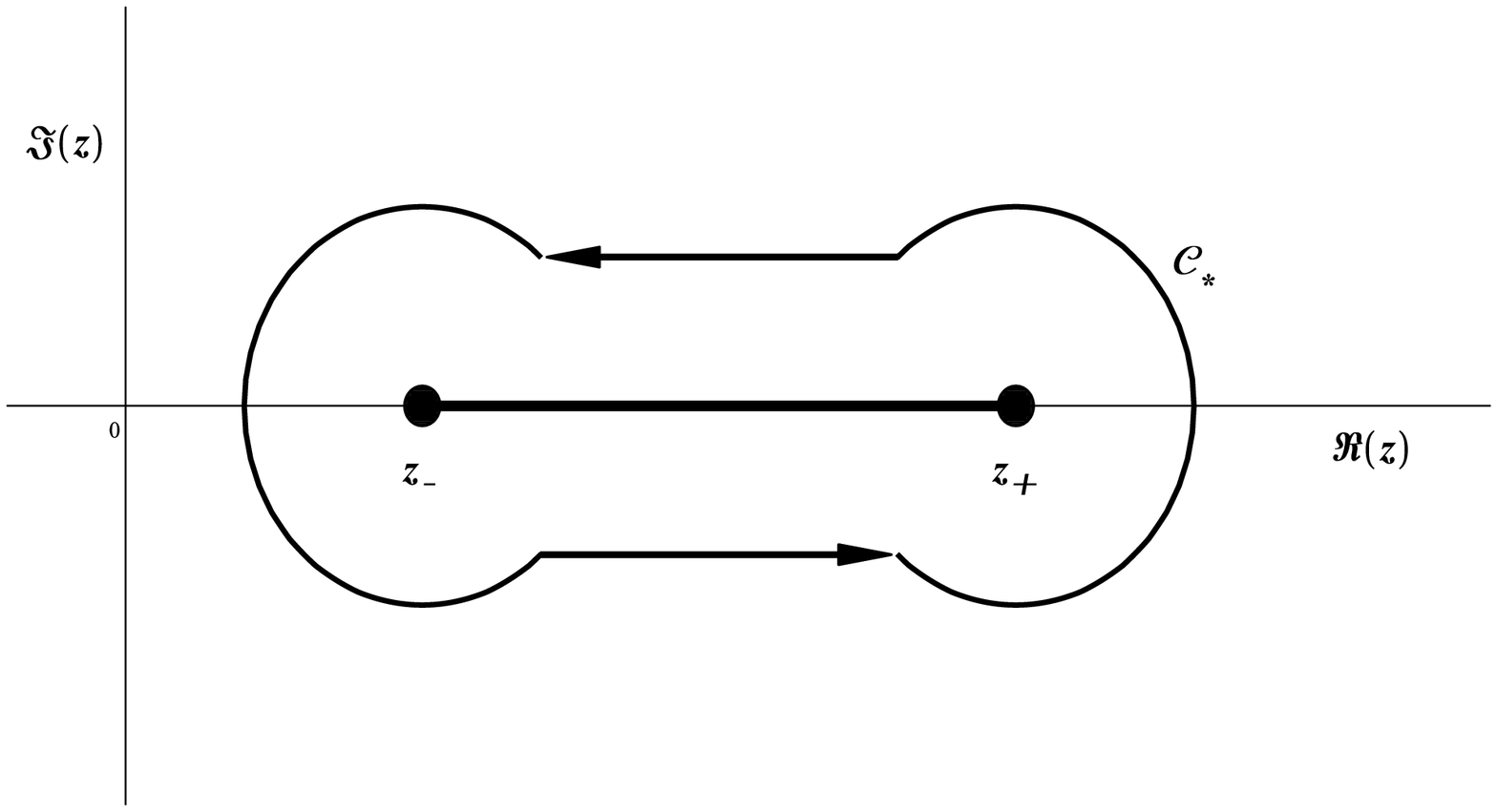}
\caption{The contour $\mathcal{C}_*$.} \label{figure:1}\end{center}

\begin{center}
\includegraphics[width=0.6\textwidth]{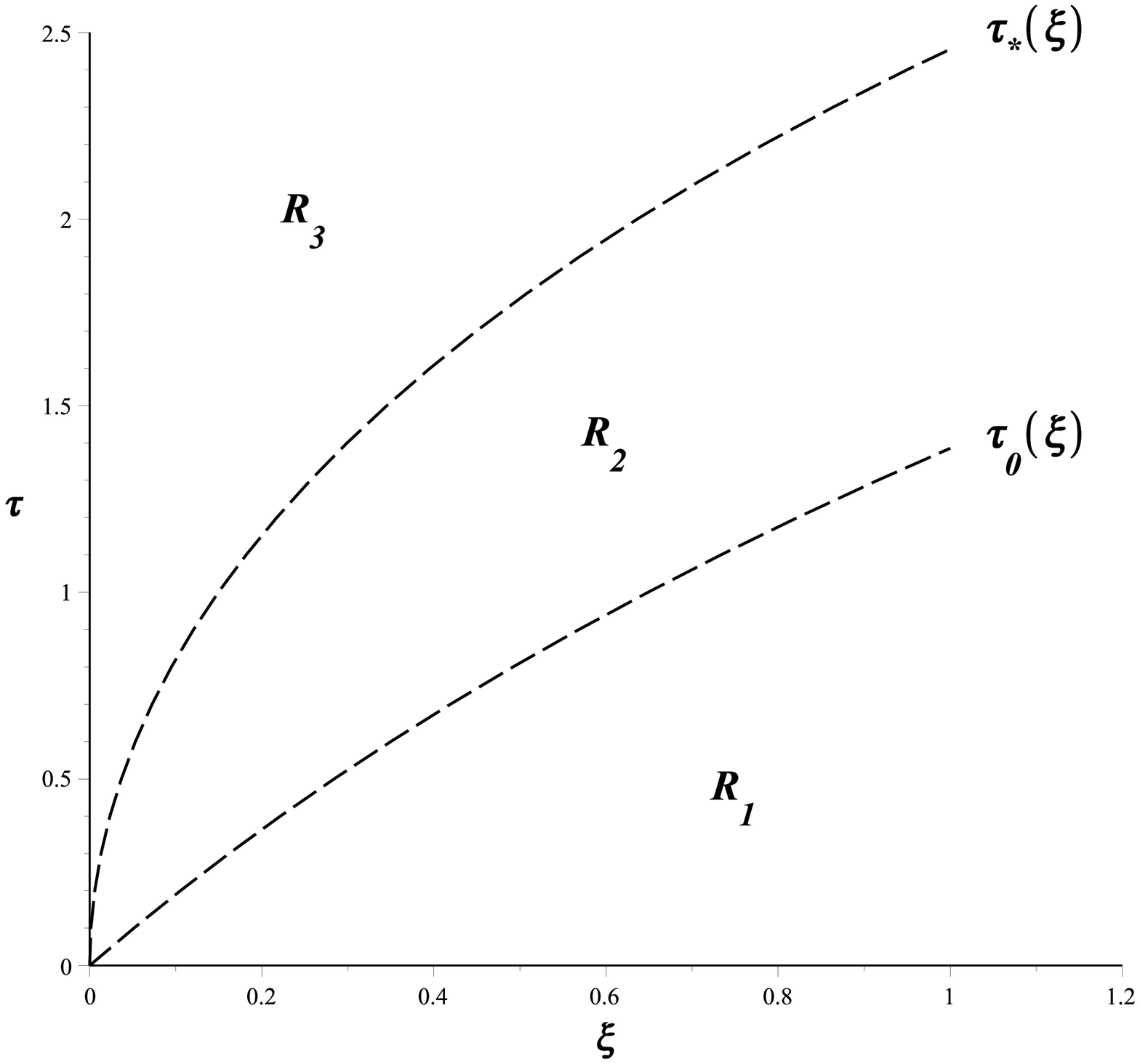}
\caption{The regions $R_1$, $R_2$, $R_3$ in the $(\xi,\tau)$ plane ($\rho=0.5$).} \label{figure:2}
\end{center}

\begin{center}
\includegraphics[width=0.6\textwidth]{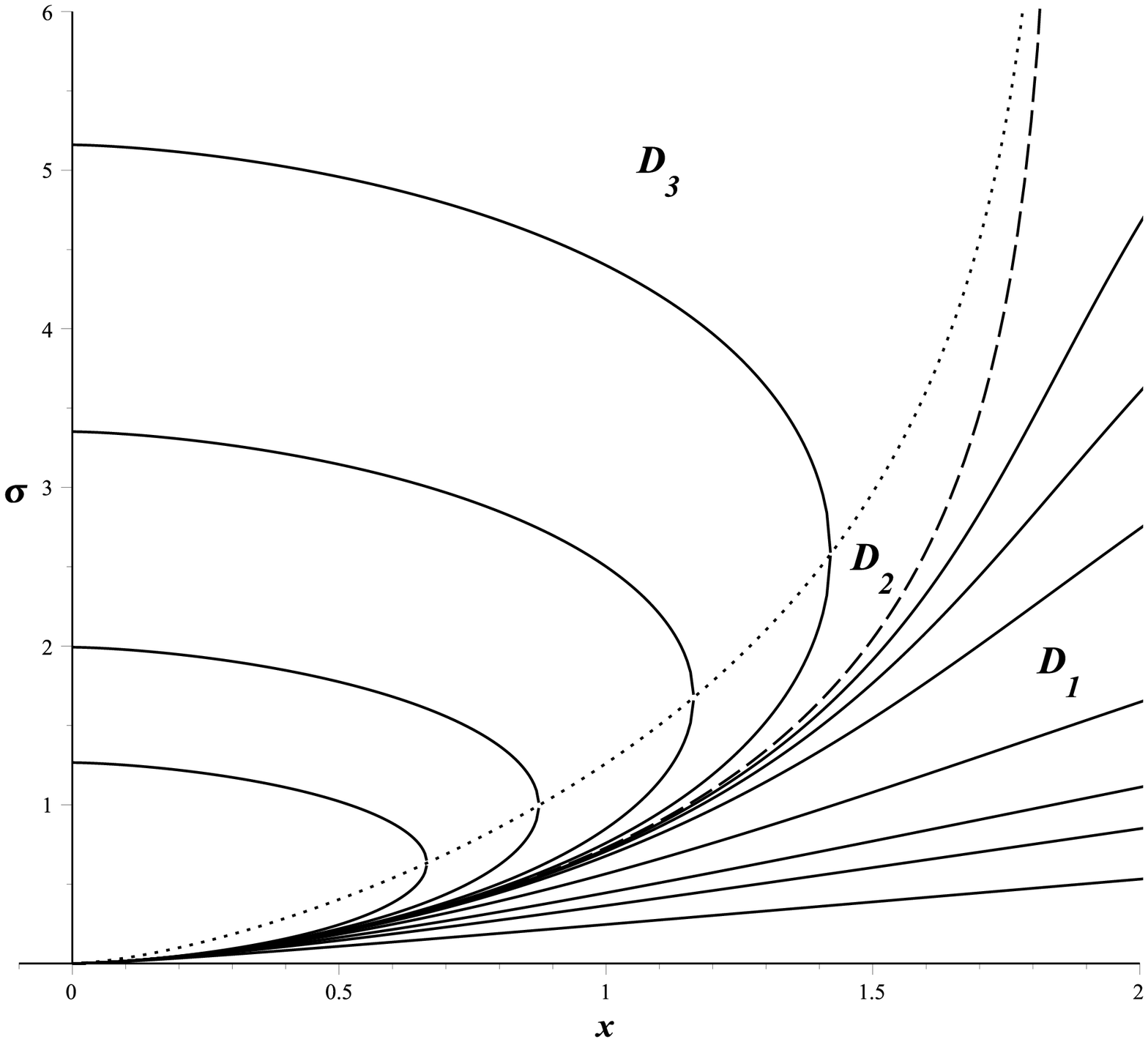}
\caption{The regions $D_1$, $D_2$, $D_3$ in the $(x,\sigma)$ plane ($\rho=0.5$).} \label{figure:3}
\end{center}

\end{document}